\documentclass[12pt]{article}
\usepackage{amssymb,amsmath}

\def\tr {{\rm tr}}
\def\Ker {{\rm Ker}}
\def\Im {{\rm Im}}
\def\dim {{\rm dim}}
\def\div {{\rm div}}
\def\exp {{\rm exp}}
\def\grad {{\rm grad}}

\begin{document}

\title{The space of associated metrics on a symplectic manifold}
\author{N.~K.~Smolentsev
\thanks{Kemerovo State University, Kemerovo, 650043, RUSSIA.
Email smolen@kemsu.ru, }
} \maketitle

\begin{abstract}
In this work the spaces of Riemannian metrics on a closed
manifold $M$ are studied. On the space ${\mathcal M}$ of all
Riemannian metrics on $M$ the various weak Riemannian structures
are defined and the corresponding connections are studied. The
space ${\mathcal AM}$ of associated metrics on a symplectic
manifold $M,\omega$ is considered in more detail. A natural
parametrization of the space ${\mathcal AM}$ is defined. It is
shown, that ${\mathcal AM}$ is a complex manifold. A curvature of
the space ${\mathcal AM}$ and quotient space ${\mathcal
AM}/{\mathcal D}_{\omega}$ is found. The spaces ${\mathcal AM}$
in cases when $M$ is a two-dimensional sphere and two-dimensional
torus are considered as application of general results. The
critical metrics of the functional of the scalar curvature on
${\mathcal AM}$ are considered. The finite dimensionality of the
space of associated metrics of a constant scalar curvature with
Hermitian Ricci tensor is shown.
\end{abstract}

\newpage

\vspace {10mm}

\centerline {\bf Contents.}

\vspace {8mm}

{\bf 1.} Preliminaries. \dotfill 3

\vspace {2mm}

\centerline{\vbox{\hsize=5.1truein \noindent 1.1. The topology in
spaces of tensor fields. \ 1.2. Differential operators $\alpha_g$
and $\delta_g$ and the decomposition of Berger-Ebin.\ 1.3. The
space of Riemannian metrics. \ 1.4. Decomposition of space
$\mathcal{M}$. \ 1.5. Action of group of diffeomorphisms on the
space of metrics.\ 1.6. ILH-manifolds.\ 1.7. Some final notation.
}}

\vspace {2mm}

{\bf 2.} Natural weak Riemannian structures on space of the
Riemannian metrics.  \dotfill 10

\vspace {2mm}

\centerline{\vbox{\hsize=5.1truein  \noindent 2.1. A flat
structure. \ 2.2. A conformally flat structure.\ 2.3. A
homogeneous structure. \ 2.4. General canonical structure.\ 2.5.
Non-Riemannian connection.}}

\vspace{2mm}

{\bf 3.} The space of associated Riemannian metrics on a
symplectic manifold. \dotfill 19

\vspace {2mm}

\centerline{\vbox{\hsize=5.1truein \noindent 3.1. The spaces of
associated metrics and almost complex structures.\ 3.2. A
parametrization of the spaces $\mathcal{A}_\omega$ and
$\mathcal{AM}$.\ 3.3. A complex structure of the space
$\mathcal{AM}$.\ 3.4. Local expressions. The Beltrami equation.}}

\vspace{2mm}

{\bf 4.} Decomposition of the  space of Riemannian metrics on a
symplectic manifold. \dotfill 33

\vspace {2mm}

{\bf 5.} Curvature of the space of associated metrics. \dotfill 37

\vspace{2mm}

{\bf 6.} Orthogonal decompositions of the space of symmetric tensors \\
on an almost K\"{a}hler manifold. \dotfill 41

\vspace{2mm}

{\bf 7.} A curvature of a quotient space $\mathcal
{AM}/\mathcal{D}_{\omega}$. \dotfill 49

\vspace{2mm}

{\bf 8.} The Spaces of the associated metrics on sphere and torus.
\dotfill 52

\vspace {2mm}

{\bf 8.1.} Associated metrics on a sphere $S^2$. \dotfill 52

\vspace {2mm}

\centerline{\vbox{\hsize=5.1truein \noindent 8.1.1. A view of
associated metrics.\ 8.1.2. A curvature of the space
$\mathcal{AM}(S^2)$.}}

\vspace{2mm}

{\bf 8.2.} Associated metrics on a torus $T^2$. \dotfill 60

\vspace {2mm}

\centerline{\vbox{\hsize=5.1truein \noindent 8.2.1. A view of
associated metrics.\ 8.2.2. A curvature of the space
$\mathcal{AM}(T^2)$.\ 8.2.3. A curvature of the space
$\mathcal{AM}/\mathcal{G}$.}}

\vspace{2mm}

{\bf 9.} Critical associated metrics on symplectic manifold.
\dotfill 69

References. \dotfill 79

\newpage

\vspace{10mm}

\centerline{\bf \S 1. Preliminaries.}

\vspace{8mm}

{\bf 1.1. The topology in spaces of tensor fields.} Let $M$ be a
smooth (i.e. $C^{\infty}$), closed, orientable $n$-manifold, $TM$
is its tangent bundle, $T^p_qM$ the bundle of $(p,q)$-tensors and
$S_2M$ the bundle of symmetric 2-forms on $M$.

The symbol $\Gamma (T^p_q)$ will denote space of all smooth
$(p,q)$-tensor fields on $M$. In particular:

$\Gamma (TM)$ is the space of all smooth vector fields on $M$,

$S_2=\Gamma (S_2M)$ is the space of all smooth symmetric 2-forms
on $M$.

Define topology in these spaces. For this purpose we fix some
smooth Riemannian metric $g$ on $M$ and denote its covariant
derivative as $\nabla$. In local coordinates $x^1,\dots,x^n$ on
$M$ let: $g_{ij}(x)$ are components of the metric tensor $g$,
$g^{ij}(x)$ are components of inverse matrix to $g_{ij}(x)$ and
$\mu_g =\left(\mbox{det}g_{ij}\right)^{1/2}dx^1\wedge\dots\wedge
dx^n$ is Riemannian volume form.

The metric $g$ defines an inner product on the space
$\Gamma(T_q^pM)$ of tensor fields of any type $(p,q)$. If $T_1$
and $T_2$ are $(p,q)$-tensor fields, then their inner product is
set by the formula:
$$
(T_1,T_2)_g=\int_Mg^{i_1k_1}\dots g^{i_qk_q}g_{j_1l_1} \dots
g_{j_pl_p}\ T_{1\ i_1\dots i_q}^{j_1\dots j_p} T_{2\ k_1\dots
k_q}^{l_1\dots l_p}\ d\mu_g.\eqno{(1.1)}
$$

The metric $g$ allows to define a stronger inner product on the
space $\Gamma (T_q^pM)$ of $(p,q)$-tensor fields. Let $s$ is an
integer non-negative number and $T_1$ and $T_2$ are
$(p,q)$-tensor fields. Then
$$
(T_1,T_2)_g^s=\sum^{s}_{i=0}(\nabla^iT_1,\nabla^iT_2)_g,\eqno
{(1.2)}
$$
where $\nabla^i=\nabla\circ\dots\circ\nabla$ is $i$-th degree of
the covariant derivative and $(\nabla^iT_1,\nabla^iT_2)$ is the
inner product (1.1).

Denote as $H^s(T_q^pM)$ the completion of space $\Gamma (T_q^pM)$
with respect to the topology in $\Gamma (T_q^pM)$ given by the
inner product (1.2). The space $H^s(T_q^pM)$ is called a space of
$(p,q)$-tensor fields Sobolev class $H^s$. It is a Hilbert space.
Denote corresponding norm as $\|.\|_s$. The fundamental property
of $H^s$ spaces is the Sobolev embedding theorem which states: If
$s\geq\frac n2+1+k$, then $H^s\subset C^k$, and the inclusion is
a continuous linear mapping \cite{Pal1}. If $s\geq\frac n2+1+k$,
then every tensor field $T$ of class $H^s$ is of a differentiable
class $C^k$. The further restrictions on $s$ connect with
necessity to ensure an appropriate class of a smoothness of
tensor fields from $H^s(T_q^pM)$.

Define topology in the space $\Gamma (T_q^pM)$ of smooth
$(p,q)$-tensor fields on $M$ by a set of norms $\{\|.\|_s,\ s\geq
0\}$. Then $\Gamma (T_q^pM)$ is the Frechet space.

In further it is supposed, that the space $\Gamma(TM)$ of smooth
vector fields on $M$ and the space $S_2=\Gamma (S_2M)$ of smooth
symmetric 2-forms on $M$ are endowed just this topology and,
thus, are infinite dimensional Frechet spaces.

The symbol $S_2^s=H^s(S_2M)$ will denote the Hilbert space of
symmetric 2-forms of a class $H^s$, $s>\frac n2+2$.

{\bf Remark.} One can show, that the topology on spaces $\Gamma
(T_q^pM)$ defined above  does not depend on a choice of the
metric $g$ on $M$.

\vspace {5mm}

{\bf 1.2. Differential operators $\alpha_g$ and $\delta_g$ and the
decomposition of Berger-Ebin.}

Let $g$ is some smooth Riemannian structure on manifold $M$ and
$\nabla$ is its covariant derivative of Riemannian connection.
Let $\nabla_i$ is a covariant derivative along the vector field
$\frac{\partial}{\partial x^i}$ in local coordinates
$x^1,\dots,x^n$ on $M$.

The symbol $g^{-1}$ will denote operation of lifting of the first
index of tensor. If $a$ is a symmetric 2-form on $M$, then
$A=g^{-1}a=g^{ik}a_{kj}$ is an endomorphism of the tangent bundle
$TM$.

The trace (with respect to $g$) of symmetric 2-form $a$ is
defined by the formula $\mbox{tr}_g\ a =\tr\ A=g^{ij}a_{ij}$.

The covariant divergence $\delta_ga$ of the symmetric 2-form $a$
on $M$ is a vector field on $M$ defined by equality $(\delta_g
a)^i=-\nabla_j a^{ij}$, where $a^{ij}=g^{ik}g^{jl}a_{kl}$.

Thus, the covariant divergence is a differential operator of the
1-st order,
$$
\delta_g:S_2\longrightarrow\Gamma (TM),\qquad
(\delta_ga)^i=-\nabla_ja^{ij}.
$$
Enter one more differential operator
$$
\alpha_g:\Gamma (TM) \longrightarrow S_2,\qquad \alpha_g(X)=\frac
12L_Xg,
$$
where $L_Xg$ is Lie derivative along the vector field $X$ on $M$,
$L_Xg=\nabla_iX_{j}+\nabla_jX_{i}$.

The spaces $\Gamma (TM)$ of vector fields on $M$ and $S_2$ of the
symmetric 2-forms have natural inner products:
$$
(X,Y)_g=\int_Mg(X,Y)d\mu_g,\qquad X,Y \in\Gamma (TM), \eqno{(1.3)}
$$
$$
(a,b)_g=\int_Mg(a,b)d\mu_g=\int_Mg^{ik}g^{jl}a_{ij}b_{kl}d\mu_g,
\qquad a,b\in S_2. \eqno{(1.4)}
$$
where $\mu_g=\left(\mbox{det}g_{ij}\right)^{1/2}\
dx^1\wedge\dots\wedge dx^n$ is Riemannian volume form of the
metric $g$.

It just follows from the Stokes theorem, that operator $\alpha_g$
is the adjoint of $\delta_g$ \cite{Bes2}: for any $X\in\Gamma
(TM)$ and $a\in S_2$,
$$
(\alpha_g(X),a)=(X,\delta_ga).
$$
Since $\alpha_g$ is operator with an injective symbol, there is
following orthogonal Berger-Ebin decomposition \cite{Ber-Eb} of
the space $S_2$:
$$
S_2=S_2^0\oplus\alpha_g(\Gamma (TM)),
$$
Where $S_2^0=\mbox{ker}\delta_g=\{a\in S_2;\ \delta_g\ a=0\}$ is
the space of divergence-free symmetric 2-forms. Then each 2-form
$a\in S_2$ has unique representation:
$$
a=a^0\ +\ L_Xg,
$$
where $\delta_ga^0=0$. The components $a^0$ and $L_Xg$ are
orthogonal and they are defined in a unique way.

There is one more (pointwise) orthogonal decomposition of the
space $S_2$:
$$
S_2=S_2^T\oplus S_2^C,
$$
where $S_2^T=\{a\in S_2;\ \mbox{tr}_g\ a=0\}$ is the space of
traceless symmetric 2-forms and $S_2^C=\{a\in S_2;\ a=\sigma g,\
\sigma\in C^{\infty}(M,\bf{R}) \}$. Each 2-form is represented in
a unique way as:
$$
a=\left(a-\frac 1n(\mbox{tr}\ a)g\right)+\frac 1n(\mbox{tr}\ a)g.
$$

\vspace {5mm}

{\bf 1.3. The space of Riemannian metrics.} Let $\mathcal{M}$ is
the space of Riemannian structures on a manifold $M$. The space
$\mathcal{M}$ is an open convex cone in the Frechet space $S_2$
(the space $S_2$ consists of all smooth symmetric 2-forms on $M$
and $\mathcal{M}$ of all positive defined symmetric 2-forms). So
the space $\mathcal{M}$ is the Frechet manifold and for any
$g\in\mathcal{M}$, the tangent space $T_g\mathcal{M}$ is
naturally identified with the space $S_2$.

The manifold $\mathcal{M}$ has a canonical weak Riemannian
structure. Namely, if $a,b\in T_g\mathcal{M}=S_2$ are two smooth
symmetric 2-forms on $M$, which represent elements of tangent
space $T_g\mathcal{M}$, then their inner product is defined by
the formula:
$$
(a,b)_g=\int_M\tr(g^{-1}ag^{-1}b)d\mu_g=\int_Mg^{ik}
g^{jl}a_{ij}b_{kl}d\mu_g. \eqno {(1.4)}
$$
This structure on $\mathcal{M}$ is called weak because the inner
product in the tangent space $T_g\mathcal{M}=S_2$ defines weaker
topology, than the topology of the Frechet space. Detail research
of the manifold $\mathcal{M}$ with a canonical Riemannian
structure (1.4) is in work of D. Ebin \cite{Ebi1}.

The curvature and geodesics of the space $\mathcal{M}$ were
studied in the work \cite{Fre-Gro} and \cite{Gil-Mic}. Other weak
Riemannian structures on $\mathcal{M}$ have been considered in
the work of author \cite {Smo23}. Recall main facts about the
space $\mathcal{M}$ obtained in the works \cite {Ebi1}, \cite
{Fre-Gro} and \cite {Gil-Mic}.

Let elements $a,b,c\in T_g\mathcal{M}=S_2$ represent constant
(parallel in $S_2$) vector fields on $\mathcal{M}$ and $A,B,C$
are their endomorphisms, $A=g^{-1}a$, $B=g^{-1}b$, $C=g^{-1}c$.
The covariant derivative $\nabla^0$ of Riemannian connection on
$\mathcal{M}$ of corresponding weak Riemannian structure (1.4) is
obtained in \cite {Ebi1}:
$$
\nabla^0_ab=-\frac 12(aB+bA)+\frac
14\left(\tr(A)b+\tr(B)a-\tr(AB)g\right), \eqno {(1.5)}
$$
where $aB=a_{ik}b^k_j$.

The curvature tensor has been found in works \cite {Fre-Gro} and
\cite {Gil-Mic} (see also \S 2 of this work):
$$
R^0(a,b)=-\frac 14g\left[[A,B],C\right]-\frac{1}{16}\tr
(C)\left(\tr(A)b-\tr(B) a\right)+\frac{n}{16}\left(\tr
(AC)b^T-\tr (BC)a^T\right), \eqno {(1.6)}
$$
where $[A,B]=AB-BA$ and $a^T=a-\frac 1n\tr (A)g$ is traceless part
of tensor $a$.

Geodesics on $\mathcal{M}$ have been found in works \cite
{Fre-Gro} and \cite{Gil-Mic}. At first, we will show them in the
form of \cite {Gil-Mic}

{\bf Theorem 1.1} \cite {Gil-Mic} {\it Let $g_0\in\mathcal{M}$ and
$a\in T_{g_0}\mathcal{M}$. Then geodesic on $\mathcal{M}$, going
out from $g_0$ in direction $a$ is curve
$$
G(t)=g_0\ e^{\alpha (t)I+\beta (t)A^T},
$$
where $I$ is identity endomorphism, $A^T$ is the traceless part
of an endomorphism $A=g^{-1}_0a$ and $\alpha (t)$, $\beta (t)$
are smooth functions on $M$, which depend on $t$ and look like:
$$
\alpha (t)=\frac 2n\ln\left((1+\frac t4\tr A)^2+\frac n{16}\tr
((A^T)^2)t^2 \right),
$$
$$
\beta(t)=\left\{
\begin{array}{l}
\frac{t}{1+\frac t4\tr A}, \quad {\rm if }\quad\tr ((A^T)^2)=0 \\
\frac{4}{\sqrt{n\tr((A^T)^2)}}\arctan\left( \frac{\sqrt{n\tr
((A^T)^2)}\ t}{4+t\ \tr A}\right), \quad {\rm if }\quad \tr
((A^T)^2)\ne 0,
\end{array}\right.,
$$
where $\arctan$ takes values in
$\left(-\frac{\pi}{2},\frac{\pi}{2}\right)$, in those points of
the manifold, where $\mbox{tr}(A)\ge 0$, and in points $x\in M$,
where $\mbox{tr}(A)<0$ we suppose:
$$
\arctan\left(\frac{\sqrt{n\tr((A^T)^2)}\ t}{4+t\ \tr
A}\right)=\left\{
\begin {array} {l}
\arctan \quad{\rm in }\quad[0,\frac{\pi}{2})\quad {\rm for }\quad t\in\left[0,-\frac{4}{\tr A}\right),\\
\frac{\pi}{2}\quad {\rm for }\quad t =-\ \frac{4}{\tr A},\\
\arctan\quad{\rm in }\quad(\frac{\pi}{2},\pi)\quad {\rm for
}\quad t\in\left(-\frac{4}{\tr A},\infty\right).
\end {array} \right.
$$
Let $N^a=\{x\in M;\ A^T(x)=0\}$, and if $N^a\ne\emptyset$, then
let $t^a=\inf\{\tr A(x))\ x\in N^a\}$. Then geodesic $g(t)$ is
defined for $t\in [0,\infty)$, if $N^a=\emptyset$ or if $t^a\ge
0$, and geodesic is defined for $t\in\left[\left.0,-\frac {4}{t ^
a} \right)\right.$ if $t^a<0$.}

\vspace {5mm}

{\bf 1.4. Decomposition of space $\mathcal{M}$.} Let
$\mathrm{Vol}(M)\subset\Gamma (\Lambda^nM)$ is the space of smooth
volume forms on $M$, i.e. space of smooth nondegenerate $n$-forms
on $M$, which give the same orientation as initial orientation on
$M$. The natural projection $vol: \mathcal{M}\longrightarrow
\mathrm{Vol}(M)$ is defined it takes each metric
$g\in\mathcal{M}$ to Riemannian volume form
$\mu_g=\left(\mbox{det}g_{ij}\right)^{1/2}\ dx^1\wedge\dots\wedge
dx^n$. Fiber of bundle $vol$ over $\mu\in Vol(M)$ is the space
$\mathcal{M}_{\mu}$ of metrics with the same Riemannian volume
form $\mu$.

The bundle $vol: \mathcal{M}\longrightarrow\mathrm{Vol}(M)$ is
trivial, the volume form $\mu\in Vol(M)$ defines decomposition of
$\mathcal{M}$ in direct product:
$$
\iota_{\mu}:\ \mathrm{Vol}(M)\times
\mathcal{M}_{\mu}\longrightarrow\mathcal{M},\qquad
(\nu,h)\rightarrow\left( \frac{\nu}{\mu}\right)^{2/n}h,
$$
$$
\varphi_{\mu}:\ \mathcal{M}\longrightarrow\mathrm{Vol}(M)\times
\mathcal{M}_{\mu},\qquad
g\rightarrow\left(\mu_g,\left(\frac{\mu}{\mu_g}\right)
^{2/n}g\right),
$$
where the positive function $\frac{\nu}{\mu}$ is defined by
equality $\nu=\frac{\nu}{\mu}\mu$.

The metric $g\in\mathcal{M}$ defines section
$$
S_{g}:\ \mathrm{Vol}(M)\longrightarrow\mathcal{M},\qquad\nu
\rightarrow\left(\frac{\nu}{\mu_g}\right)^{2/n}g.
$$

Using decomposition
$\mathcal{M}=\mathrm{Vol}(M)\times\mathcal{M}_{\mu}$, in the work
\cite {Fre-Gro} the following expressions of geodesics on space
$\mathcal{M}$ are obtained.

{\bf Theorem 1.2} \cite {Fre-Gro} {\it The geodesics in
$\mathcal{M}=\mathrm{Vol}(M)\times\mathcal{M}_{\mu}$ with initial
position $(\mu,g)$ and initial velocity $(\beta,b)\in\Gamma
(\Lambda^nM)\times S_2^T$ is
$$
g_t=\left(q(t)^2+r^2t^2\right)^{2/n}g\ \exp\left(\frac
1r\arctan\left( \frac{rt}{q}\right) B\right),
$$
where $B=g^{-1}b$, \ $\exp$  is the exponential mapping,
$q(t)=1+\frac 12\frac{\beta}{\mu}t$, $r=\frac
14\sqrt{n\mbox{tr}(B^2)}$. (If $r=0$, replace the exponential
term by 1.) The change in the volume form is given by the formula}
$$
\mu(g_t)=\left(q(t)^2+r^2t^2\right)\mu.
$$
{\bf Remark.} All facts explained above are true for the space
$\mathcal{M}^s$ of Riemannian metrics on $M$ of Sobolev class
$H^s$, $s>\frac n2+2$.

\vspace {3mm}

{\bf 1.5. Action of group of diffeomorphisms on the space of
metrics.}

Let $\mathcal{D}^{s+1}$ be the group of diffeomorphisms of the
manifold $M$ of Sobolev class $H^{s+1}$, $s>\frac n2+2$. Than
group $\mathcal{D}^{s+1}$ acts on the smooth Hilbert manifold
$\mathcal{M}^s$ in the following way:
$$
A:\ \mathcal{M}^s\times\mathcal{D}^{s+1}\longrightarrow
\mathcal{M}^s,\qquad A(g,\eta)=\eta^*g,\eqno {(1.7)}
$$
$$
\eta^*g(x)(X,Y)=g(\eta (x))\left(d\eta (X),d\eta (Y)\right),
$$
for any vector fields $X,Y$ on $M$ and any $x\in M$.

The group $\mathcal{D}^{s+1}$ is the smooth Hilbertian manifold,
but it is not a Lie group, because group operations are only
continuous. The action $\mathcal{D}^{s+1}$ on $\mathcal{M}^s$ is
also only continuous. Nevertheless D.Ebin has shown \cite {Ebi1},
that for any metric $g\in\mathcal{M}$ and its group of isometries
$I(g)$, orbit $g\mathcal{D}^{s+1}$ of action $A$ is a smooth
closed submanifold, which is diffeomorphic to quotient space
$\mathcal{D}^{s+1}/I(g)$. Differential of mapping
$$
A_g:\ \mathcal{D}^{s+1}\longrightarrow\mathcal{M}^s,\qquad
A_g(\eta)=\eta^*g,
$$
is Lie derivative:
$$
dA_g:\ \Gamma(TM)\longrightarrow S_2^s,\quad dA_g(X) =
L_Xg=2\alpha_g(X).
$$
Moreover, there is following

{\bf Theorem 1.3.} (Slice theorem, \cite {Ebi1}). {\it Let
$s>\frac n2+2$. For each metric $g\in\mathcal{M}^s$ there exists
a smooth submanifold $\mathcal{S}_g^s\subset\mathcal{M}^s$
containing $g$, such that

1) If $\eta\in I(g)$, then
$\eta^*(\mathcal{S}_g^s)=\mathcal{S}_g^s$,

2) If $\eta\in\mathcal{D}^{s+1}$ and $\eta^*(\mathcal{S}_g^s)
\cap\mathcal{S}_g^s\neq\emptyset$, then $\eta\in I(g)$,

3) There exists a local cross section $\chi:\
\mathcal{D}^{s+1}/I(g)\longrightarrow\mathcal{D}^{s+1}$ defined
in a neighbourhood $U^{s+1}$ of identity coset
$[e]\in\mathcal{D}^{s+1}/I(g)$, such that mapping
$$
F:\ \mathcal{S}_g^s\times U^{s+1}\longrightarrow\mathcal{M}^{s},
\qquad F(h,u)=(\chi (u))^*h,
$$
is a homeomorphism onto a neighbourhood $V^s$ of element
$g\in\mathcal{M}^s$.}

Notice, that the map $F:\ \mathcal{S}_g\times
U\longrightarrow\mathcal{M}$ of the slice theorem is ILH-smooth,
as for any $s\geq 2n+5$ is hold \cite{Koi}:
$$
\mathcal{S}_g^s=\mathcal{S}_g^{2n+5}\cap\mathcal{M}^s,\qquad
U^{s+1}=U^{2n+5}\cap(\mathcal{D}^{s+1}/I(g)),
$$
$$
V^s=V^{2n+5}\cap\mathcal{M}^s,\qquad\chi^{s+1}=\chi^{2n+5}|U^{s+1}
$$
and for any $k\ge 0$ mappings
$$
F^{s+k}:\ \mathcal{S}_g^{s+k}\times U^{s+1+k}\longrightarrow
V^{s},
$$
$$
p^{s+k}\times q^{s+k}:\
V^{s+k}\longrightarrow\mathcal{S}_g^{s}\times U^{s+1}
$$
are $C^k$ differentiable.

The quotient space $\mathcal{M}/\mathcal{D}$ is not a manifold.
Actually, $\mathcal{D}$ does not act freely. Elements
$g\in\mathcal{M}$ have isotrophy groups $I(g)$ depending on
$g\in\mathcal{M}$. It is known \cite {Ebi1}, that the set of
metrics $g\in\mathcal{M}$ with trivial group of isometries is an
open dense set $\mathcal{M}^*$ in $\mathcal{M}$. Group
$\mathcal{D}$ acts on the space $\mathcal{M}^*$ freely. On the
slice theorem, we obtain, that the quotient space
$\mathcal{M}^*/\mathcal{D}$ is a ILH-smooth manifold.

Since the tangent space to an orbit $g\mathcal{D}$ is identified
with space
$$
\alpha_g\left(\Gamma (TM)\right)=\{h\in S_2;\ h=L_Xg,\ X\in\Gamma
(TM)\},
$$
then (Berger-Ebin decomposition):
$$
T_g\mathcal{M}=S_2=S_2^0\oplus T_g(g\mathcal{D}).
$$
We obtain, that the tangent space
$T_{[g]}(\mathcal{M}^{*}/\mathcal{D})$ is identified with space
$$
S_2^0=\{h\in S_2;\ \delta_gh=0\}
$$
of divergenceless 2-forms $h$.

\vspace {5mm}

{\bf 1.6. ILH-manifolds.}

{\bf Definition 1.1.} {\it Topological vector space $E$ is called
ILH-space, if $E$ is an inverse limit of Hilbert spaces $\{E^s;\
s=1,2,\dots\}$, such that $E^l\subset E^s$, if $l\geq s$, and
this inclusion is a bounded linear operator.}

We will denote $E=\lim_{\leftarrow}E^s$.

{\bf Definition 1.2.} {\it Topological space $X$ is called
$C^k$-ILH-manifold, which are modeled on ILH-space $E$, if

a) $X$ is an inverse limit of $C^k$-smooth Hilbert manifolds
$\{X^s\}$, modeled on $\{E^s\}$, and if $l\geq s$, $X^l\subset
X^s$;

b) For any point $x\in X$ there are open neighbourhoods $U^s(x)$
of point $x$ in $X^s$ and homeomorphisms $\psi^s$ of
neighbourhoods $U^s(x)$ on open subsets $V^s(x)\subset E^s$,
which define $C^k$-coordinates on $X^s$ in a neighbourhood of the
point $x$, such that $U^l(x)\subset U^s(x)$ at $l\geq s$ and
$\psi^{s+l}(y)=\psi^{s}(y)$ for any point $y\in U^{s+l}(x)$.}

{\bf Definition 1.3.} {\it Let $X,Y$ are $C^k$-ILH-manifolds. The
mapping $\varphi:\ X\longrightarrow Y$ is called
$C^k$-ILH-differentiable, if $\varphi$ is an inverse limit of
$C^k$-differentiable mappings of Hilbert manifolds $X^l$ and
$Y^s$, i.e. if for any $s$ there is a number $l(s)$ and
$C^k$-differentiable mapping $\varphi^s:\ X^{l(s)}\longrightarrow
Y^{s}$, such that $\varphi^s(x)=\varphi^{s+1}(x)$, $\forall\ x\in
X^{l(s+1)}$ and $\varphi=\lim_{\leftarrow}\varphi^s$.}

\vspace {3mm}

{ \bf 1.7. Some final notation.}

$\mathrm{End}(TM)$ is the vector bundle of endomorphisms $K:TM\to
TM$ of the tangent bundle. Fiber over a point $x\in M$ consists
of all endomorphisms of the tangent space $T_xM$. Endomorphism
$K:TM\to TM$ will be also called an operator, acting on the
tangent bundle.

$A(M)$ is bundle over $M$ whose fiber $A_x(M)$ over the point
$x\in M$ consists of automorphisms $J_x$ of the tangent space
$T_xM$, such that: $J_x^2=-I_x$, where $I_x$ is identity
automorphism of the space $T_xM$. It is supposed, that dimension
$n$ of manifold $M$ is even, $n=2m$.

{\bf Definition 1.4.} {\it An almost complex structure (a.c.s.) on
$M$ is the smooth section $J$ of bundle $A(M)$.}

Thus, the almost complex structure on $M$ is a smooth automorphism
$J:\ TM\longrightarrow TM$, such that $J^2=-I$, where $I$ is
identity automorphism.

$\mathcal{A}=\Gamma (A(M))$ is space of all smooth almost complex
structures on $M$. This is the space of smooth sections of bundle
$A(M)$, therefore $\mathcal{A}$ is infinite-dimensional smooth
ILH-manifold \cite {Abr}.

\newpage

\vspace {10mm}

\centerline{\bf \S 2. Natural weak Riemannian structures}
\centerline {\bf on space of the Riemannian metrics.}

\vspace {7mm}

As above, there is a canonicaly weak Riemannian structure on the
space $\mathcal{M}$ of all Riemannian metrics on manifold $M$
(1.4). In this paragraph we shall consider a series of other
natural weak Riemannian structures on $\mathcal{M}$ and we shall
obtain their formulas for covariant derivative, curvature tensor,
sectional curvatures and geodesic.

\vspace {2mm}

{\bf 2.1. A flat structure.} Let $g_0$ is fixed Riemannian metric
on $\mathcal{M}$. The formula
$$
(a,b)_0^{\alpha}=\int_M\tr(g_0^{-1}ag_0^{-1}b)d\mu_0+
\alpha\int_M\tr(g_0^{-1}a)\tr(g_0^{-1}b)d\mu_0,\eqno{(2.1)}
$$
where $a,b\in S_2$, $\alpha\in{\bf R}$ is some number and
$d\mu_0=d\mu(g_0)$ is the Riemannian volume form, defines the
symmetric form in a vector space $S_2$. It is positive defined
for $\alpha>-\frac 1n$ and nondegenerate for $\alpha\ne -\frac
1n$. Therefore for $\alpha >-\frac 1n$ we obtain a flat weak
Riemannian structure on $S_2$.

The covariant derivative $d$ is a usual directional derivative in
a vector space $S_2$: if $b=b(g)$ is vector field on
$\mathcal{M}$, $d_ab=\left.\frac {d}{dt}\right |_{t=0}b(g+ta)$.
The curvature tensor is equal to zero. Geodesics are: $g_t=g+ta$.

\vspace {3mm}

{\bf 2.2. A conformally flat structure.} Consider the following
weak Riemannian structure on $\mathcal{M}$:
$$
\left<a,b\right>_g=\int_M\tr(g_0^{-1}ag_0^{-1}b)d\mu(g), \eqno
{(2.2)}
$$
where $a,b\in S_2$, $g_0$ is fixed metric on $\mathcal{M}$. In
contrast to the previous case, the volume form $\mu(g)$ depends
on $g\in\mathcal{M}$. It differs from $\mu(g_0)$ on smooth
positive function: $\mu(g)=\rho(g)\mu(g_0)$. Function $\rho(g)$
we shall call a denseness of the Riemannian metric $g$ with
respect to $g_0$.

Write weak Riemannian structure (2.2) as
$$
\left<a,b\right>_g=\int_M\tr(g_0^{-1}ag_0^{-1}b)\rho(g)d\mu(g_0),\eqno
{(2.3)}
$$
To find a covariant derivative $D$ of metric (2.2) on
$\mathcal{M}$ it is used usual "six-term formula"
\cite{Gr-Kl-Me}. As volume form $\mu (g_0)$ is constant in the
integral (2.3) (i.e. does not depend on $g\in\mathcal{M}$) and
application of six-term formula has no differentiation on $x\in
M$, then it is enough to calculate a covariant derivative $D_x$,
curvature tensor $K_x$ and geodesics for the Riemannian metric
$$
\left<a,b\right>_{g,x}=\tr (g_0^{-1}ag_0^{-1}b)(x)\rho(g)(x)
\eqno {(2.4)}
$$
on the space $\mathcal{M}_x$ of inner products on tangent space
$T_xM$ at each point $x\in M$.

We will use following notation:
$$
A=g^{-1}a,\ G=g^{-1}g_0,\ A_0=g_0^{-1}a,\ G_0=g_0^{-1}g.\eqno
{(2.5)}
$$

{\bf Theorem 2.1.} {\it Weak Riemannian structure (2.2) on
$\mathcal{M}$ has the following characteristics:

1) Covariant derivative
$$
D_ab=d_ab+\frac 14\left(\tr (A)b+\tr (B)a\right)-\frac 14\tr
(A_0B_0)g_0G,
$$

2) Curvature tensor
$$
K(a,b)c=
$$
$$
=-\frac {1}{16}\left(\left(4\tr (AC)+\tr (A)\tr (C)\right)b -
\left(4\tr (BC)-\tr (B)\tr (C)\right)a\right)-\frac {1}{16} \tr
(A_0C_0)g_0 \left(4B+\tr (B)\right)G+
$$
$$
+\frac {1}{16}\tr (B_0C_0) g_0\left(4A+\tr
(A)\right)G+\frac{1}{16}\tr (G^2) \left(\tr (A_0C_0)b-\tr
(B_0C_0)a\right).
$$

3) The sectional curvature $K_{\sigma}$ in a plane section
$\sigma$, given by orthonormal vectors $a,b\in T_g\mathcal{M}$ is
equal to zero, if $a,b$ are scalar tensors ($a=\alpha g$,
$b=\beta g$), and in case, when $a,b$ are traceless ($\tr A=0$,
$\tr B=0$), it is expressed by the formula
$$
K_{\sigma}=\int_M\left(-\frac {1}{2}\tr (AB)\tr (A_0B_0)+\frac
{1}{4}\left( \tr (A^2)\tr (B_0^2)+\tr (B^2) \tr
(A_0^2)\right)-\right.
$$
$$
\left.-\frac {1}{16}\tr (G^2)\left(\tr (A_0^2)\tr
(B_0^2)-(\tr(A_0B_0))^2\right) \right)d\mu(g).
$$

4) Geodesics $g_t$ on $\mathcal{M}$ are solutions of the following
differential equation of the second order,
$$
\frac {d}{dt}\left(\rho(g)\frac {dg}{dt}\right)=kg_0G,
$$
where $k=k(x)$ is positive function on $M$, which does not depend
on $t$ and is equal to one eighth length of an initial velocity
of geodesics $g_t(x)$ on $\mathcal{M}_x$, $x\in M$.}

{\bf Proof.} We will make computations on $\mathcal{M}_x$, $x\in
M$. The metric (2.4) on $\mathcal{M}_x$ is conformally equivalent
to the metric $\left<a,b\right>_x=\tr (g_0^{-1}ag_0^{-1}b)(x)$.
Therefore to find a derivative $D_x$ we can apply usual formula
(see for example \cite{Gr-Kl-Me}):
$$
D_ab=d_ab+\frac 12\left(a(\psi)+b(\psi)a -
\left<a,b\right>_xd\psi\right), \eqno {(2.6)}
$$
where $\psi=\ln (\rho (g))$, $d\psi$ is gradient of function
$\psi$. A curvature tensor is found in the same way. Directional
derivative $a(\psi)$ and gradient $d\psi$ of function $\psi$ are
simply found: if $g_t=g+ta$, then from equality $\mu (g_t)=\rho
(g_t)\mu$ is obtained
$$
a(\rho)=\left.\frac {d}{dt}\right|_{t=0}\rho(g_t)=\frac 12\tr
(g^{-1}a)\rho(g)=\frac 12\tr (A)\rho (g).
$$
For $\psi=\ln(\rho)$ is obtained,
$$
a(\psi)=\frac 12\tr (g^{-1}a)=\frac 12
\left<g_0g^{-1}g_0,a\right>_x.
$$
Therefore
$$
a(\psi)=\frac 12\tr (A),\qquad d\psi=\frac 12 g_0G,\eqno{(2.7)}
$$
Now expression for a covariant derivative is:
$$
D_ab=d_ab+\frac 14\left(\tr (A)b+\tr (B)a\right)-\frac 14\tr
(A_0B_0)g_0G,
$$
To find a curvature tensor $K_x$ of conformally equivalent metric
it is needed a Hessian tensor $H_{\psi}(a)=d_a(d\psi)$ and Hessian
$h_{\psi}(a,b)=\left<d_a(d\psi),b\right>_{0,x}$. They are easily
calculated,
$$
d_a(d\psi)=\left.\frac {d}{dt}\right|_{t=0}\left(\frac 12
g_0g_t^{-1}g_0 \right)=-\frac 12 g_0g^{-1}ag^{-1}g_0=-\frac 12
g_0AG,
$$
$$
h_{\psi}(a,b)=-\frac
12\tr\left(g_0^{-1}g_0g^{-1}ag^{-1}g_0g_0^{-1}b\right) =-\frac
12\tr\left(g^{-1}ag^{-1}b\right)=-\frac 12\tr (AB).
$$
Put the given expressions in the formula for a curvature tensor
of conformally equivalent metric \cite{Gr-Kl-Me}, we obtain a
curvature tensor $K_x$ on $\mathcal{M}_x$ and, therefore, on
$\mathcal{M}$.

Consider the equation of geodesics
$$
D_aa=a'+\frac 12\tr (A)a-\frac 12\tr (A_0^2)d\psi=0. \eqno {(2.8)}
$$
If $g_t$ is geodesic, then velocity $a=g_t'$ has a constant
length: $\left<a,a\right>_{g,x}=\tr (A_0^2)\rho(g)(x)= c(x)$. Let
multiply the equation (2.8) on $\rho(g)$,
$$
a'\rho+\frac 12\tr (A)a\rho-\frac 12\tr (A_0^2)\rho d\psi=0.
$$
Then
$$
\frac {d}{dt}(a\rho)=\frac 12\left<a,a\right>_{g,x}d\psi=k(x)g_0G.
$$
The theorem is proved.

{\bf Remark 1.} Setting formulas for a covariant derivative,
curvature tensor and geodesics on $\mathcal {M}$ give us the same
characteristic of Riemannian manifold $\mathcal{M}_x$ with the
metric
$$
\left<a,b\right>_{g,x}=\tr (g_0^{-1}ag_0^{-1}b)(x)\sqrt{\det
g_{ij}(x)}.
$$

{ \bf Remark 2.} For more general Riemannian structure on
$\mathcal{M}_x$
$$
\left<a,b\right>_{g,x}^{\alpha}= \tr
(g_0^{-1}ag_0^{-1}b)\rho(g)(x)+ \alpha\tr (g_0^{-1}a)\tr
(g_0^{-1}b)\rho(g)(x)
$$
expressions for a gradient and Hessian are the following:
$$
d\psi=\frac 12 g_0G-\frac{\alpha}{2(1+\alpha n)}\tr (G)g_0,
$$
$$
H_{\psi}(a)=-\frac 12 g_0AG+\frac {\alpha}{2(1+\alpha n)}\tr
(AG)g_0,
$$
$$
h_{\psi}(a,b)=-\frac {1}{2}\tr (AB).
$$
Using these formulas we can easily obtain a covariant derivative
and curvature tensor on $\mathcal{M}$ of a weak Riemannian
structure
$$
\left<a,b\right>_{g}^{\alpha}=\int_M\tr
(A_0B_0)d\mu(g)+\alpha\int_M \tr (A_0)\tr (B_0)d\mu (g).
$$

\vspace {3mm}

{ \bf 2.3. A homogeneous structure.} The inner product in a point
$g\in\mathcal{M}$ is defined by the formula:
$$
(a,b)_{g}=\int_M\tr (AB)d\mu (g_0).\eqno {(2.9)}
$$
where $a,b\in T_g\mathcal{M}$, $g_0$ is fixed metric on $M$,
$A=g^{-1}a$.

\vspace {3mm}

{\bf Theorem 2.2.} {\it Weak Riemannian structure (2.9) has the
following geometric characteristics:

1) Covariant derivative,
$$
\nabla_ab=d_ab-\frac 12\left(aB+bA\right),
$$

2) Curvature tensor,
$$
R(a,b)c=-\frac 14g\ [[A,B],C],
$$

3) Sectional curvature $K_{\sigma}$ in a plane section $\sigma$,
given by orthonormal pair $a,b\in T_g\mathcal{M}$,
$$
K_{\sigma}=\frac 14\int_M\tr \left([A,B]^2\right)d\mu (g_0),
$$

4) Geodesics, going out from a point $g\in\mathcal{M}$ in
direction $a\in T_g\mathcal{M}$ look like $g_t=g\ e^{tA}$.}

\vspace {3mm}

{\bf Proof.} As the volume form $\mu (g_0)$ is constant (does not
depend from $g\in\mathcal{M}$), then computations can be done
with integrand expression from (2.9), i.e. on the
finite-dimensional manifold $\mathcal{M}_x$ with the metric
$(a,b)_{g,x}=\tr(AB)(x)$

Let $a,b\in S_{2,x}(M)$. We will think, that they define parallel
vector fields on $\mathcal{M}_x\subset S_{2,x}(M)$, then $d_ab=0$
and $[a,b]=d_ab-d_ba=0$. For $g\in\mathcal{M}_x$ we have
$A=g^{-1}a$ and $B=g^{-1}b$. Let $g_t=g+tc$ is curve on
$\mathcal{M}_x$, going out in direction $c\in S_{2,x}(M)$. Then,
$$
d_cA=\left.\frac {d}{dt}\right|_{t=0}A=\left.\frac {d}{dt}
\right|_{t=0}g^{-1}a=-g^{-1}cg^{-1}a=-CA.
$$
To find a covariant derivative $\nabla_ab$ we apply six-term
formula
$$
2(\nabla_ab,c)=a(b,c)+b(c,a)-c(b,a)+(c,[a,b])+(b,[c,a])-(a,[b,c])=
$$
$$
=d_a(\tr BC)+d_b(\tr CA)-d_c(\tr BA)=-\tr (ABC)-\tr (BAC) -\tr
(BCA)-
$$
$$
-\tr (CBA)+\tr (CBA)+\tr (BCA)=-\tr ((AB+BA)C)=-(aB+bA,c)_{g,x}.
$$
Therefore for a constant vector fields $a$ and $b$ on
$\mathcal{M}_x$ we have:
$$
\nabla_ab=-\frac 12 (aB+bA)=\Gamma_g(a,b).
$$

Calculate a curvature tensor
$R(a,b)c=\nabla_a\nabla_bc-\nabla_b\nabla_ac$.
$$
\nabla_a\nabla_bc=d_a(\nabla_bc)+\Gamma_g(a,\nabla_bc)=
\left.\frac {d}{dt}\right|_{t=0}\left(-\frac 12 (bg_t^{-1}c+cg
t^{-1}b)\right)-\frac 12(ag^{-1}\nabla_bc+\nabla_bc\ A)=
$$
$$
=-\frac 12(-bg^{-1}ag^{-1}c-cg^{-1}ag^{-1}b)+\frac
14(a(BC+CB)+(bCA+cBA))=
$$
$$
=\frac 12(bAC+cAB)+\frac 14(aBC+aCB+bCA+cBA).
$$
Similarly,
$$
\nabla_b\nabla_ac=\frac 12(aBC+cBA)+\frac 14(bAC+bCA+aCB+cAB).
$$
From this we obtain,
$$R(a,b)c=\nabla_a\nabla_bc-\nabla_b\nabla_ac=
-\frac 14\left((aB-bA)C-c(AB-BA)\right)=
$$
$$
=-\frac 14g\left((AB-BA)C-C(AB-BA)\right)=-\frac 14g[[A,B],C].
$$

Find geodesics of this metric. Let $g_t$ is geodesic and $a=g_t'$.
Then $\nabla_aa=a'-aA=0$. Write the equation $a'-aA=0$ as
$g_t^{-1}a'-A^2=0$. Then the last equation is equivalent to the
following: $A'=0$, where $A=g_t^{-1}a$. Solution of the last one
is a constant operator: $A(t)=A$. Then $a(t)=g_tA$, or
$g_t'=g_tA$. Therefore $g_t=g\ e^{tA}$.

\vspace{3mm}

{\bf Corollary.} {\it Submanifold $\mathcal{M}_{\mu}$ of metrics
$g$ with the same volume form $\mu$ and submanifold
$\mathcal{P}g$ of pointwise conformally equivalent to
$g\in\mathcal{M}$ metrics are totally geodesic in $\mathcal{M}$
with respect to weak Riemannian structure (2.9).}

\vspace{3mm}

{\bf Proof.} If the initial velocity $a$ of geodesic $g_t$
touches a submanifold $\mathcal{M}_{\mu}$, $\tr\ A=0$. Therefore
for geodesic $g_t=ge^{tA}$ we have $\mu (g_t)=\mu
(g)\det^{1/2}(e^{tA})=\mu (g)$, therefore,
$g_t\in\mathcal{M}_{\mu}$. Similarly, if $a$ touches submanifold
$\mathcal{P}g$, then $a=\alpha g$, where $\alpha$ is function on
$M$. Then $A=\alpha I$ and $e^{tA}=e^{t\alpha}I$.

\vspace{3mm}

{\bf Remark 1.} Just the same results are obtained for more
general weak Riemannian structure
$$
(a,b)_{g,\alpha}=\int_M\tr (AB)d\mu (g_0)+\alpha\int_M\tr (A)\tr
(B)d\mu (g_0).
$$

{\bf Proof.} It is enough to calculate a covariant derivative. On
the six-term formula we have:
$$
2(\nabla_ab,c)=a(b,c)_{g,\alpha}+b(c,a)_{g,\alpha}-c(b,a)_{g,\alpha}=
$$
$$
=d_a(\tr (BC)+\alpha\tr (B)\tr (C))+d_b(\tr (CA)\alpha\tr (C)\tr
(A))-d_c( \tr (BA)\alpha\tr (B)\tr (A))=
$$
$$
=-\tr (ABC)-\tr (BAC)-\alpha\tr (AB)\tr (C)-\alpha\tr (B)\tr (AC)-
$$
$$
-\tr (BCA)-\tr (CBA)-\alpha\tr (BC)\tr (A)-\alpha\tr (C)\tr (BA)+
$$
$$
+\tr (CBA)+\tr (BCA)+\alpha\tr (CB)\tr (A)+\alpha\tr (B)\tr (CA)=
$$
$$
=\tr ((AB+BA)C)-\alpha(\tr(AB+BA)\tr (C))=-(aB+bA,c)_{g,\alpha}.
$$
Sectional curvature
$$
K_{\sigma,\alpha}=(R(a,b)b,a)_{g,\alpha}=-\frac 14\tr ([[A,B]B]A)-
\frac {\alpha}{4}\tr ([[A,B]B])\tr (A)=
$$
$$
=-\frac 14\tr ([[A,B]B]A)=-\frac 14\tr ([A,B][B, A])=\frac 14\tr
([A,B]^2).
$$

{\bf Remark 2.} Weak Riemannian structure (2.9) is most
convenient for study of the manifold $\mathcal{M}$. It has simple
formulas for a covariant derivative, curvature and geodesics, the
submanifolds $\mathcal{M}_{\mu}$ and $\mathcal{P}g$ are totally
geodesic in $\mathcal{M}$. But the structure (2.9) is
non-invariant with respect to action of group of diffeomorphisms
$\mathcal{D}(M)$.

\vspace {3mm}

{\bf 2.4. General canonical structure.} We shall consider a weak
Riemannian structure, which is more general, than canonical,
$$
(a,b)_{g}^{\alpha}=\int_M\tr (AB)d\mu (g)+\alpha\int_M\tr (A)\tr
(B)d\mu (g). \eqno {(2.10)}
$$

If $\alpha\ne -1$, then given structure is called as Devitt
metric, it arises at Hamilton description of a general relativity
theory \cite {Pek1}.

\vspace {3mm}

{\bf Theorem 2.3.} {\it Weak Riemannian structure (2.10) has the
following geometric characteristics:

1) Covariant derivative
$$
\nabla_a^{\alpha} b = d_a b - \frac 12 (aB+bA) + \frac 14
\left(\tr (A)b +\tr (B)a\right) - \frac {1}{4(1+\alpha
n)}\left(\tr (AB) + \alpha \tr (A)\tr (B)\right)g,
$$

2) Curvature tensor,
$$
R^{\alpha}(a,b)c=-\frac {1}{4}g[[A,B],C] - \frac {1}{16}\tr
(C)\left(\tr (A)b-\tr (B)a\right)+
$$
$$
+\frac {n}{16(1+\alpha n)}\left((a,c)_{x}^{\alpha}b^T -
(b,c)_{x}^{\alpha}a^T\right),
$$
where $(a,c)_{x}^{\alpha} =\tr (AC)+\alpha\tr (A)\tr (C)$, $b^T=b-
\frac {1}{n}\tr (B)g$ is a traceless part of a tensor $b$,

3) Sectional curvature $R_{\sigma}^{\alpha}$ in a plane section
$\sigma$, given by orthonormal pair $a,b\in T_g\mathcal{M}$:
$$
K_{\sigma}^{\alpha}=\frac {1}{4} \int_M \tr ([A,B])^2 d\mu(g) -
 \frac {n}{16(1+\alpha n)}\int_M \left((a,a)_{x}^{\alpha}(b,b)_{x}^{\alpha}-
((a,b)_{x}^{\alpha})^2 \right)d\mu(g)+
$$
$$
+\frac {1}{16} \int_M \left(\tr (A^2)\tr ^2(B)+ \tr (B^2)\tr ^2(A)
- 2\tr (A)\tr (B)\tr (AB)\right)d\mu(g);
$$

4) The geodesics in
$\mathcal{M}=\mathrm{Vol}(M)\times\mathcal{M}_{\mu}$ with initial
position $(\mu,g)$ and initial velocity $(\beta,b)\in\Gamma
(\Lambda^nM)\times S_2^T$ is
$$
g_t=\left(q(t)^2+r^2t^2\right)^{2/n}g%
\ \exp\left(\frac 1r \arctan \left(\frac {rt}{q}\right)B\right),
$$
where $B=g^{-1}b$,\ $\exp$ is exponential mapping, $q(t)=1+\frac
12\frac {\beta}{\mu}t$, $r=\frac 14\sqrt{n\tr (B^2)}$. If $r=0$,
replace the exponential term by 1.}

\vspace {3mm}

{\bf Proof.} Write weak Riemannian structure (2.10) as
$$
(a,b)_{g}^{\alpha}=\int_M\left(\tr (AB)+\alpha\tr (A)\tr
(B)\right)\rho (g) d\mu (g_0).
$$
Integrand expression defines Riemannian metric on $\mathcal{M}_x$,
which is conformally equivalent to metric $(a,b)_{x}^{\alpha}=\tr
(AB)+\alpha\tr (A)\tr (B)$, considered in item 2.3. Therefore
covariant derivative and a curvature tensor are easily obtained
from a covariant derivative $\nabla$ and a curvature tensor $R$
with using of equalities:
$$
\psi=\ln\rho (g),\quad d\psi=\frac {1}{2(1+\alpha n)}g,\quad
H_{\psi}(a)=0.
$$
Set view of geodesics. Suppouse, that $g_t$ is geodesic and
$a(t)=\frac {d}{dt}g_t$. Then from the equation
$\nabla_{a}^{\alpha}a=0$ we obtain,
$$
a'=aA+\frac 12 \tr(A)a-\frac {1}{4(1+\alpha n)}%
\left(\tr (A^2)-\alpha \tr ^2(A)\right)g=0.
$$
It is easy to see, that for $A(t)=g_t^{-1}a$ the derivative $A'$
is found from the formula $A'=g_t^{-1}(a'-aA)$. Therefore
equation of geodesics is written as:
$$
A'+\frac 12 \tr(A)A-\frac {1}{4(1+\alpha n)}%
\left(\tr(A^2)-\alpha \tr^2(A)\right) = 0.
$$
Enter the following value: $f=\tr (A)$, $E=A-\frac fn I$, $v =\tr
(E^2)$, where $I$ is identity operator. Then $A=E+\frac fn I$ and
$\tr (A^2)=v+\frac {f^2}{n}$. Substituting it in the equation of
geodesics, we obtain the pair of the equations,
$$
f'+\frac 14 f^2-\frac {1}{4(1+\alpha n)}v=0,
$$
$$
E'+\frac 12fE=0.
$$
Now use decomposition of space $\mathcal{M}$ in direct product
$\mathcal{M}=Vol(M)\times\mathcal{M}_{\mu}$. Then $g_t$ is
represented as a pair $g_t=(\mu_t,h_t)$. Let $\mu_t=\rho (t)\mu$.
We have $g_t=\rho^{2/n}(t)h_t$ and $a(t)=\frac {2}{n}\frac {\rho
'}{\rho}g_t+g_th_t^{-1}h_t'$. Therefore $A(t)=\frac{2}{n}\frac
{\rho'}{\rho}I+h_t^{-1}h_t'$. It follows that
$f=2\frac{\rho'}{\rho}$ and $E=h_t^{-1}h_t'$. Let
$u=\frac{\rho'}{\rho}$, then previous sistem of equations
rewrited as:
$$
u'+\frac 12u^2-\frac {1}{8(1+\alpha n)}v=0,\qquad E'=-uE.
$$
As $v'=(\tr (E^2))'=2\tr (EE')=2\tr (-uE^2)=-2uv$, then by
differentiation of the first equation, we obtain,
$u''+uu'=-\frac{Uv}{4(1+\alpha n)}=-2u(u'+\frac 12 u^2)$, i.e.,
$$
u''+3uu'+u^3=0.
$$
This equation, and the second one $E'=-uE$, have coincided with
the similar equations in work \cite {Fre-Gro} for a determination
of the geodesics of canonical metric. Therefore geodesics $g_t$
of our Riemannian structure $(a,b)_{g}^{\alpha}$ coincide with
geodesics of canonical Riemannian structure on $\mathcal{M}$,
found in work \cite {Fre-Gro}.

\vspace {3mm}

{\bf 2.5. Non-Riemannian connection.} In this section we consider
connection on $\mathcal{M}$, which takes an intermediate position
between Riemannian connection $\nabla$ of the section 2.3 and
canonical connection $\nabla^0$. It is defined by the following
covariant derivative:
$$
\overline{\nabla}_ab=d_ab-\frac 12(aB+bA)+\frac 18\left(\tr
(A)b+\tr (B)a\right).
$$
It is easy to see, that the bilinear form
$$
Q(a,b)=\int_M\tr (A)\tr (B)d\mu (g)
$$
is invariant with respect to $\overline{\nabla}$:
$AQ(b,c)=Q(\overline{\nabla}_ab,c)+Q(b,\overline{\nabla}_ac)$.

\vspace {3mm}

{\bf Theorem 2.4.} {\it 1) Curvature tensor of connection
$\overline{\nabla}$ expressed by the formula
$$
\overline{R}(a,b)c=-\frac {1}{4}g[[A,B],C]-\frac{1}{64}\left(\tr
(A)b- \tr (B)a\right)\tr (C).
$$

2) Geodesics on $\mathcal{M}$, going out from a point
$g\in\mathcal{M}$ in direction $A_0=\frac{4\beta}{n}I+B$, $\tr
(B)=0$ look like:
$$
g_t=\left\{
\begin{array}{l}
(\beta t+1)^{4/n}g\ \exp (\frac {1}{\beta}\ln (\beta t+1)B),\quad
\beta\ne 0\\
g\ \exp (Bt),\quad \beta=0
\end{array}
\right..
$$
}

{\bf Proof.} Curvature tensor is calculated immediately. Consider
the equation of geodesics $\overline{\nabla}_aa=a'-aA+\frac 14\tr
(A)a=0$. Applying an operator $A(t)=g_t^{-1}a(t)$, we reduce it
to a view
$$
A'+\frac 14\tr (A)A=0.
$$
Let $f=\tr (A)$ and $A=\frac {f}{n}I+E$ is decomposition $A$ on a
scalar and traceless parts. The equation for $A$ falls to two
ones:
$$
f'+\frac 14f^2=0,\qquad E'+\frac 14fE=0.
$$
Solution of the first equation, with the initial condition
$f(0)=\beta $, is found immediately: $f=4\beta (\beta t+4)^{-1}$.
Before to decide the second equation, we decompose $g_t$ on two
components according to the decomposition
$\mathcal{M}=Vol(M)\times\mathcal{M}_{\mu}$. Then
$g_t=(\mu_t,h_t)$. Let $\mu_t=\rho (t)\mu$. We have
$g_t=\rho^{2/n}(t)h_t$ and $a(t)=\frac
{2}{n}\frac{\rho'}{\rho}g_t+g_th_t^{-1}h_t'$. Therefore
$A(t)=\frac {2}{n}\frac {\rho'}{\rho}I+h_t^{-1}h_t'$.
Consequently $f=2\frac {\rho'}{\rho}$ and $E=h_t^{-1}h_t'$. From
the first equation, taking into account the initial condition
$\rho (0)=1$, we find: $\rho=\frac {1}{16}(\beta t+4)^{2}$.
Solution of the second equation $E'+\frac 14fE=0$ is easily found
now: $E(t)=\frac {E(0)}{\sqrt{\rho}}=\frac {4B}{\beta t+4}$. Let
$g_0$ is initial point of geodesic $g_t$. It is easy to see, that
$E(t)=h_t^{-1}h_t'=\left(\ln (g_0^{-1}h_t)\right)'$. From the
equation $\left(\ln(g_0^{-1}h_t)\right)'=\frac {4B}{\beta t+4}$
is obtained $\ln(g_0^{-1}h_t)=\frac{4\ln(\beta t+4)B}{\beta}+c$.
Since  when $t=0$ it is fulfilled $h_0=g_0$, $c=-\frac {4\ln
4}{\beta}B$. Therefore
$$
g_0^{-1}h_t=e^c\ \exp\left(\frac {4\ln (\beta
t+4)}{\beta}B\right)=\exp \left(\frac
{4}{\beta}\ln\left(\frac{\beta t+4)}{4}\right)
B\right),\quad\beta\ne 0,
$$
$$
g_0^{-1}h_t=\exp(Bt),\quad\beta=0.
$$
As $g_t=\rho^{2/n}(t)h_t$, from the last expression, after
replacement $\beta$ on $4\beta$, we just obtain the view of
geodesic $g_t$, which is pointed in the theorem.

\vspace {3mm}

{\bf Remark.} We have considered basic cases of natural weak
Riemannian structures on $ \mathcal {M} $. Other variants also
can be studied, it is possible to find corresponding covariant
derivatives and curvature tensor. For example, if
$$
(a,b)_{g}=\int_M\tr (A_0B_0)d\mu(g_0)+\alpha\int_M\tr (A)\tr
(B)D\mu (g),
$$
then the covariant derivative has view:
$$
\nabla_ab=d_ab-\frac 12(aB+bA)-\frac {\alpha}{1+\alpha\rho\tr
(G^2)} \left(\tr (AB)-\frac 14\tr (A)\tr (B)\right)g_0g^{-1}g_0,
$$
where the function $\rho$ is found from a condition $\mu (g)=\rho
(g)\mu (g_0)$. The remaining characteristics are found similarly.

\newpage

\vspace{10mm}

\centerline {\bf \S 3. The space of associated Riemannian metrics
on} \centerline {\bf a symplectic manifold.}

\vspace {7mm}

In this paragraph we shall suppose, that manifold $M$ is
symplectic. It means, that on $M$ the closed nondegenerate 2-form
$\omega$ of class $C^{\infty}$ is given. The manifold $M$ has an
even dimension, $\dim M=2n$.

\vspace {3mm}

{\bf 3.1. The spaces of associated metrics and almost complex
structures.} On a symplectic manifold it is natural to consider
the metrics, which are compatible with the symplectic form
$\omega$. It seams, that such metrics are well connected with
almost complex structures (further a.c.s.) on $M$.

If the a.c.s. $J$ is given on $M$, then it is always convenient
to have also the Riemannian metric $g$. It is natural to demand
from $g$ to be Hermitian with respect to $J$:\ $g(JX,JY)=g(X,Y)$.
As known \cite{Ko-No}, there exists Hermitian metric $g$, for
every a.c.s. $J$, but it is not unique.

The task is to connect only one Hermitian metric with every
a.c.s. $J$. In case of Riemannian surfaces the complex structure
defines a class of the conformally equivalent metrics, the choice
of metric in this class is defined by the demand of a constancy
of curvature. If $n>1$ then the almost complex structure $J$ does
not define a class of the conformally equivalent metrics. However
there is means of unique choise of Hermitian metric $g$ on a
symplectic manifold. Let's give the necessary definitions.

\vspace {3mm}

{\bf Definition 3.1.} {\it An almost complex structure on a
manifold $M$ is the endomorphism of a tangent bundle $J:\
TM\rightarrow TM $, such that: $J^2=-I$, where $I$ is identity
automorphism.}

\vspace {3mm}

{\bf Definition 3.2.} {\it An almost complex structure $J$ on $M$
is called \textit{positive associated} with the symplectic form
$\omega$, if for any vector fields $X,Y$ on $M$ the following
conditions are hold:

1) $\omega (JX,JY)=\omega (X,Y)$,

2) $\omega (X,JX)>0$, if $X\neq 0$.}

\vspace {3mm}

{\bf Definition 3.3.}{\it Each positive associated a.c.s. $J$
defines the Riemannian metric $g$ on $M$ by equality
$$
g(X,Y)=\omega (X,JY),\eqno {(3.1)}
$$
which is also called \textit{associated}.}

\vspace {3mm}

The associated metric $g$ has the following properties:

1) $g(JX,JY)=g(X,Y)$,

2) $g(JX,Y)=\omega (X,Y)$.

\vspace {3mm}

{\bf Remark.} Sometimes positive associated almost complex
structure $J$ is called \textit{calibrating} $2$-form $\omega$
(an exterior $2$-form $\omega$ on $X$ is called $J$-calibrated if
$\omega (JX,JY)=\omega (X,Y)$ and $\omega (X,JX)>0$ for $X\neq 0$
\cite {Gro}).

Almost complex structure $J$, satisfying to the condition of positiveness \\
$\omega (X,JX)>0$, if $X\neq 0$ is also called \textit{tamed} to
form $\omega$ (we say that an exterior $2$ -form $\omega$ on $M$
tames an almost-complex structure $J$ if $\omega (X,JX)>0$ for
$X\neq 0$ \cite{Gro}).

Our terminology is offered by D.Blair \cite{Bla1}, \cite{Bla-Ian}
and it seems more natural. It also corresponds to the terminology
used in case of contact manifolds.

\vspace {3mm}

Let $\mathcal{A}$ is the space of all smooth almost complex
structures on $M$. It is the space of smooth sections of the
bundle $A(M)$ over $M$, fiber of which is $A_x(M)$ over the point
$x\in M$ consists of automorphisms $J_x$ of tangent space $T_xM$,
such that: $J_x^2=-I_x$, where $I_x$ is identity automorphism of
space $T_xM$. As the space $\mathcal{A}=\Gamma (A(M))$ is spase
of sections it is infinite-dimensional, smooth ILH-manifold
\cite{Abr}.

In this paragraph we shall also consider the following spaces:

\vspace {2mm}

$\mathcal{A}_\omega$ is the space of all smooth positive
associated almost complex structures on a symplectic manifold
$M^{2n},\omega$;

\vspace {2mm}

$\mathcal{AM}$ is the space of all smooth associated metrics on a
symplectic manifold $M^{2n},\omega$. It is clear, that
$\mathcal{AM}$ is the space of all smooth almost K\"{a}hlerian
metrics on a symplectic manifold, which fundamental form
coinciding with $\omega$.

\vspace {2mm}

These spaces $\mathcal{A}_\omega$ and $\mathcal{AM}$ are spaces
of smooth sections of corresponding bundles over $M$. Therefore
\cite{Abr} they are infinite-dimensional, smooth ILH-manifolds.
The corresponding series of smooth Hilbert manifolds form spaces
$\mathcal{A}_\omega^s$ and $\mathcal{AM}^s$ of positive
associated almost complex structures and, respectively,
associated metrics of Sobolev class $H^s$, $s>n+1$. In this
paragraph we shall not use the spaces $\mathcal{A}_\omega^s$ and
$\mathcal{AM}^s$, we only remark, that for these spaces are
fulfilled all facts, which are obtained for spaces
$\mathcal{A}_\omega$ and $\mathcal{AM}$.

Let $J\in\mathcal{A}$ is an almost complex structure on $M$. Find
tangent space $T_J\mathcal{A}$. For this purpose there is enough
to differentiate the condition $J^2=-1$. Let $J_t$ is
differentiable set of almost complex structures and
$K=\left.\frac{d}{dt}\right|_{t=0}J_t$ tangent element to the
curve $J_t$ on $\mathcal{A}$. Differentiating a condition
$J_t^2=-1$, we obtain,
$$
JK+KJ=0.
$$
Let $\mathrm{End}_{J}(TM)$ is the space of smooth endomorphisms
$K:\ TM\to TM$, anticommutating with $J$. Thus,
$$
T_J\mathcal{A}=\mathrm{End}_{J}(TM).
$$

Now we find a tangent space $T_J\mathcal{A}_\omega$ to a manifold
of positive associated almost complex structures. Let
$J\in\mathcal{A}_\omega$ is positive associated a.c.s. And $g$ is
its associated metric on $M$. It is enough to differentiate two
conditions:
$$
J_t^2=-1,\qquad\omega (J_tX,J_tY)=\omega (X,Y).
$$
Let $P=\left.\frac{d}{dt}\right|_{t=0}J_t$ is tangent element to
a curve $J_t$ on $\mathcal{A}_\omega$. Differentiating the written
out conditions, we obtain
$$
JP+PJ=0,\qquad\omega (PX,JY)+\omega (JX,PY)=0.
$$
The second condition is written more convenient via the
associated metric:
$$
g(PX,Y)-g(X,PY)=0
$$
This is condition of symmetry of an endomorphism $P$.

Thus, the tangent space $T_J\mathcal{A}_\omega$ consists of
symmetrical endomorphisms $P$, anticommutating with $J$
$$
T_J\mathcal{A}_\omega=\{P\in\mathrm{End}(TM);\quad PJ=-JP,\ \
g(PX,Y)=g(X,PY)\}.
$$
Let $\mathrm{End}_{SJ}(TM)$ is the space of smooth symmetrical
endomorphisms $P:\ TM\to TM$, anticommutating with $J$. We
obtain, that
$$
T_J\mathcal{A}_\omega=\mathrm{End}_{SJ}(TM).
$$

Similarly one can show, that the tangent space $T_g\mathcal{AM}$
to manifold $\mathcal{AM}$ at a point $g$ consists of the
anti-Hermitian symmetric 2-forms on $M$,
$$
T_g\mathcal{AM}=\{h\in S_{2};\qquad h(JX,JY)=-h(X,Y),\ \forall
X,Y\in\Gamma (TM)\}.
$$
The condition to be anti-Hermitian is $h_{ik}J^k_j=h_{kj}J^k_i$
in local coordinates. From property of $h$ to be anti-Hermitian
we obtain, in particular, that $\tr h=0$.

Denote spaces of anti-Hermitian and Hermitian symmetric 2-forms
on $M$ as $S_{2A}$ and $S_{2H}$ respectively. The
natural(pointwise) decomposition takes place
$$
S_2=S_{2A}\oplus S_{2H},\eqno {(3.2)}
$$
$$
h(X,Y)=\frac 12\left(h(X,Y)-h(JX,JY)\right)+\frac 12
\left((h(X,Y)+h(JX,JY)\right),
$$
It is orthogonal with respect to inner product (1.4) in $S_2$.
Then,
$$
T_g\mathcal{AM}=S_{2A}.
$$

The correspondence between positive associated almost complex
structures and associated metrics is defined by the diffeomorphism
$$
\mathrm{G}:\ \mathcal{A}_\omega\longrightarrow\mathcal{AM},
$$
$$
J\longrightarrow G(J)=g,\qquad g(X,Y)=\omega (X,JY).\eqno {(3.3)}
$$
In coordinates, $g_{ij}=(\mathrm{G}(J))_{ij}=\omega_{ik}J^k_j$.

Inverse diffeomorphism:
$$
\mathrm{J}:\ \mathcal{AM}^s\longrightarrow\mathcal{A}_\omega^s,
$$
$$
g\longrightarrow J,\qquad J^i_j=\omega^{ik}g_{kj}.
$$

Differentiating the relation $g_t(X,Y)=\omega (X,J_tY)$, we find
the differential of the diffeomorphism $\mathrm{G}$:
$$
d\mathrm {G}:\ T_J\mathcal{A}_\omega\longrightarrow
T_g\mathcal{AM},\qquad P\longrightarrow h=\omega P,
$$
$$
h(X,Y)=\omega (X,PY)=g(X,PJY).\eqno{(3.4)}
$$
There is also the following relation
$$
h(JX,Y)=h(X,JY)=-g(X,PY).\eqno {(3.5)}
$$

\vspace {3mm}

{\bf 3.2. A parametrization of the spaces $\mathcal{A}_\omega$ and
$\mathcal{AM}$.} As it already was marked, the spaces
$\mathcal{A}$, $\mathcal{A}_\omega$ and $\mathcal{AM}$ are smooth
ILH-manifolds. Therefore one can enter local maps on them by the
usual way \cite{Abr}. We shall show, that these spaces allow a
more natural parametrization via Cayley transformation.

Let $J_0$ is some fixed almost complex structure. As above, the
tangent space $T_{J_0}\mathcal{A}$ consists of endomorphisms
$K:TM\to TM$, anticommutating with a.c.s. $J_0$,\ $KJ_0=-J_0K$.
Therefore, for exponential mapping $e^K$ we have
$J_0e^K=e^{-K}J_0$. It just follows from here, that
$$
J=J_0e^K
$$
is an almost complex structure. The last relation gives a
parametrization of the space $\mathcal{A}$ in a neighbourhood of
the element $J_0$ by endomorphisms $K$, anticommutating with
$J_0$:
$$
E:\ \mathrm{End}_{J_0}(TM)\longrightarrow\mathcal{A},\quad
K\mapsto J=J_0e^K.
$$
Sometimes in the theory of matrices instead of a transcendental
dependence $w=e^{iz}$ the rational one: $w=\frac{1+iz}{1-iz}$,
$z=i\frac{1-w}{1+w}$ is used. Apply this transformation to an
operator $K$, possessing a property $KJ_0=-J_0K$, we obtain:
$$
J=J_0(1+K)(1-K)^{-1}.
$$
It is easy to see, that
$$
J_0(1+K)(1-K)^{-1}=(1-K)(1+K)^{-1} J_0.
$$
Therefore, $J$ is an almost complex structure. At the definition
of such almost complex structure the nondegeneracy of an operator
$1-K$ was assumed. It is enough for this purpose, that at each
point $x\in M$ the operator $K(x)$ does not have eigenvalues
equal to unit. It is clear, that the set of such endomorphisms is
an open set in the space $\mathrm{End}_{J_0}(TM)$ of all
endomorphisms $K:TM\to TM$, anticommutating with $J_0$. Denote
this set as $\mathcal{V}(J_0)$,
$$
\mathcal{V}(J_0)=\{K\in\mathrm{End}(TM);\ \ KJ_0=-J_0K,\ 1-K\
\mbox{is converted} \}.
$$

{\bf Proposition 3.1.} {\it Relations
$$
J=J_0(1+K)(1-K)^{-1},\eqno {(3.6)}
$$
$$
K=(1-JJ_0)^{-1}(1+JJ_0),\eqno {(3.7)}
$$
state the one-to-one correspondence between the set of
endomorphisms $K:TM\to TM$, anticommutating with a.c.s. $J_0$,
such that $1-K$ is converted and the set of almost complex
structures $J$ on $M$, for which an endomorphism $1-JJ_0$ is
converted.}

{\bf Proof.} Let $1-K$ is converted. Then for a.c.s.
$J=J_0(1+K)(1-K)^{-1}$ we have $JJ_0=J_0(1+K)(1-K)^{-1}J_0$,
$$
1-JJ_0=-J_0J_0-J_0(1+K)(1-K)^{-1}J_0=
-J_0(1-K+(1+K))(1-K)^{-1}J_0=
$$
$$
=-2J_0(1-K)^{-1}J_0
$$
is converted operator.

Conversely, suppose, that $1-JJ_0$ is converted, then it follows
from $1-JJ_0=-J(J+J_0)$, that the operator $J+J_0$ is converted.
From the relation (3.6), we obtain
$$
J-J_0=(J+J_0)K,\quad J(1+JJ_0)=J(1-JJ_0)K,\quad
K=(1-JJ_0)^{-1}(1+JJ_0).
$$
Moreover
$$
K=(J+J_0)^{-1}(J-J_0).
$$
Then
$1-K=1+\left(J+J_0\right)^{-1}\left(J-J_0\right)=\left(J+J_0\right)^{-1}
\left(J+J_0+J-J_0\right)=\left(J+J_0\right)^{-1}2J$ is converted
operator.

\vspace {3mm}

{\bf Remark.} The algebraic more understandable relation follows
from the relation (3.6)
$$
J=(1-K)J_0(1-K)^{-1},\eqno {(3.8)}
$$

\vspace {3mm}

The set
$$
\mathcal{U}(J_0)=\{J\in\mathcal{A}; \ 1-JJ_0\ {\rm is isomorphism
}TM \}
$$
is an open set in the space $\mathcal{A}$. Therefore map
$$
\Phi:\ \mathcal{U}(J_0)\longrightarrow\mathcal{V}(J_0),\qquad
J\mapsto K,
$$
$$
K=(1-JJ_0)^{-1}(1+JJ_0),\eqno {(3.9)}
$$
gives local coordinates in a neighbourhood of the element $J_0$.
If $K=\Phi(J)$, it is obvious that
$$
J=J_0(1+K)(1-K)^{-1}.
$$

\vspace {3mm}

The "change-of-coordinates formulas" are easily obtained from
(3.9). If $J\in\mathcal{U}(J_0)\cap\mathcal{U}(J_1)$ and
$K=(1-JJ_0)^{-1}(1+JJ_0)$,\ $P=(1-JJ_1)^{-1}(1+JJ_1)$, then
$$
P=(1-(1-K)(1+K)^{-1}J_0J_1)^{-1}(1+(1-K)(1+K)^{-1}J_0J_1).\eqno
{(3.10)}
$$

\vspace {3mm}

Parametrize the space $\mathcal{A}_\omega$ of positive associated
almost complex structures. The element $J\in\mathcal{A}_\omega$
has two properties:

1) $\omega (JX,JY)=\omega (X,Y)$,

2) $\omega (X,JX)>0$, if $X\neq 0$.

It was shown earlier, that the first property insures symmetry of
the tangent element $K\in T_J\mathcal{A}_\omega$. The second
property of a positiveness marks out an open set in the space
$\mathcal{A}$ of all almost complex structures. Enter notation
for this set:
$$
\mathcal{U}=\{J\in\mathcal{A};\ \omega (X,JX)>0,\ \mathrm{if}\
X\neq 0 \}.
$$
Further, the simple analysis shows, that if $J$ and $J_0$ are
positive almost complex structures, they are close, in the sense
that the both belong to coordinate neighbourhood
$\mathcal{U}(J_0)=\{J\in\mathcal{A};\ 1-JJ_0\mbox{ is isomorphism
}TM\}$ entered above.

\vspace {3mm}

{\bf Proposition 3.2.} {\it Let $J_0$ is positive associated
almost complex structure. An almost complex structure $J$ is
positive if and only if the following operators are positive
$$
-J_0J,\quad J^TJ_0,\quad -JJ_0,\quad J_0J^T,
$$
where the transposition is taken with respect to the metric
$g_0$, associated with $J_0$.}

\vspace {3mm}

{\bf Proof.} For associated a.c.s. $J_0$ we have
$g_0(J_0X,J_0Y)=g_0(X,Y)$ and $g_0(X,J_0Y)=-\omega (X,Y)$. Let
$J\in\mathcal{U}$ is positive a.c.s. Then the operator $S=-J_0J$
is positive if and only if $J$ is positive:
$g_0(X,SX)=-g_0(X,J_0JX)=\omega (X,JX)>0$. For an operator
$-JJ_0$ we represent an arbitrary vector $X$ as $X=J_0Y$, then
$g_0(X,-JJ_0X)=-g_0(J_0Y,JJ_0J_0Y)=g_0(J_0Y,JY)=\omega (Y,JY)>0$.
The transposed operators $J^TJ_0,\ J_0J^T$ are also positive.

\vspace {3mm}

{\bf Corollary 1.} {\it Let $J_0$ is positive associated almost
complex structure. An almost complex structure $J$ is positive
iff it can be presented as
$$
J=J_0S,
$$
where the operator $S$ is positive with respect to the metric
$g_0$, associated with $J_0$.}

\vspace {3mm}

{\bf Corollary 2.} {\it Let $J_0$ is positive associated almost
complex structure. If $J$ is any other positive a.c.s., then the
operator $1-JJ_0$ is converted.}

\vspace {3mm}

{\bf Corollary 3.} {\it Positive almost complex structure $J$ can
be represented as $J=J_0(1+K)(1-K)^{-1}$, where an operator
$(1+K)(1-K)^{-1}$ is positive defined with respect to the metric
$g_0$, associated with $J_0$.}

\vspace {3mm}

{\bf Proof.} Operator $1-JJ_0$ is converted because it is a sum of
two positive operators: identity $1$ and $-JJ_0$. As $1-JJ_0$ is
converted, the almost complex structure $J$ can be represented as
$J=J_0(1+K)(1-K)^{-1}$, where the operator $S=(1+K)(1-K)^{-1}$ is
positively defined with respect to the metric $g_0$, associated
with $J_0$.

\vspace {3mm}

The elements of the space $\mathcal{A}_\omega$ are characterized
as follows.

\vspace {3mm}

{\bf Proposition 3.3.} {\it Let $J_0$ is positive associated
almost complex structure and $g_0$ is \\ corresponding to $J_0$
associated metric. The almost complex structure $J$ is positive
associated iff it is represented as $J=J_0(1+P)(1-P)^{-1}$, where
the endomorphism $P:\ TM\to TM$ has properties:

1) $PJ_0=-J_0P$,

2) $P$ is symmetric with respect to $g_0$,

3) $1-P^2$ is positive with respect to $g_0$.}

\vspace {3mm}

{\bf Proof.} Let $J\in\mathcal{A}_\omega$. It is easily checked,
that the operator $S=-J_0J$ is symmetric with respect to $g_0$:
$$
g_0(X,-J_0JY)=\omega (X,JY)=-\omega (JX,Y)=\omega
(JX,J_0J_0Y)=g_0(JX,J_0Y) =g_0(-J_0JX,Y).
$$
Therefore, $P=-(1+S)^{-1}(1-S)$ is also symmetric with respect to
$g_0$. Moreover, the operator $S$ is positive, we shall express
$S$ through $P$, $S=(1+P)(1-P)^{-1}=(1-P^2)((1-P)^{-1})^2$. Then
$1-P^2=S(1-P)^2>0$.

Conversely, let operator $P$ is symmetric with respect to $g_0$,
and $1-P^2$ is positive. Then $1-P$ is converted and the operator
$S=(1+P)(1-P)^{-1}$ is symmetric. As $PJ_0=-J_0P$, it is easy to
see, that $J=J_0S$ is a.c.s. and
$$
\omega (JX,Y)=\omega (J_0SX,Y)=-\omega
(SX,J_0Y)=-g_0(SX,Y)=-g_0(X,SY)=
$$
$$
=g_0(J_0J_0X,SY)=-g_0(J_0X,J_0SY)=\omega (X,-J_0SY)=-\omega
(X,JY).
$$
It follows from here, that $\omega (JX,JY)=\omega (X,Y)$. The
second property, $\omega (X,JX)>0$, follows from a positiveness
of $1-P^2$. In fact, $S=(1-P^2)((1-P)^{-1})^2>0$, therefore
$$
\omega (X,JX)=\omega (X,J_0SX)=g_0(X,SX)>0,\quad\mbox{if }\ X\ne
0.
$$

\vspace {3mm}

{\bf Remark.} If $J=J_0(1+P)(1-P)^{-1}$ is positive associated
a.c.s. and $g$ is associated metric, corresponding to it, then
$J_0=J(1-P)(1+P)^{-1}$ and it is easily checked, that the
endomorphism $P$ has properties:

1) $J\ P=-P\ J$,

2) $P$ is symmetric with respect to $g$,

3) $1-P^2$ is positive with respect to $g$.

\vspace {3mm}

Recall, that the tangent space $T_J\mathcal{A}_\omega$ at a point
$J_0$ coincides with the space $\mathrm{End}_{SJ_0}(TM)$ of
smooth symmetrical endomorphisms $P:\ TM\to TM$, anticommutating
with $J_0$. The condition of positiveness $1-P^2>0$ marks out an
open set in this space. Denote it as $\mathcal{P}_{J_0}$:
$$
\mathcal{P}_{J_0}=\{P\in\mathrm{End}_{SJ_0}(TM):\ 1-P^2>0\}.
$$

It follows from the proposition 3.3 that map
$$
\Psi:\ \mathcal{P}_{J_0}\longrightarrow\mathcal{A}_\omega,\quad
P\mapsto J=J_0(1+P)(1-P)^{-1}\eqno {(3.11)}
$$
gives a global parametrization of the space $\mathcal{A}_\omega$
of positive associated almost complex structures.

\vspace {3mm}

{\bf Proposition 3.4.} {\it Relation $J=J_0e^P$ defines the
one-to-one correspondence between the space
$\mathrm{End}_{SJ_0}(TM)$ of symmetrical endomorphisms $P:TM\to
TM$, anticommutating with $J_0$ and the space
$\mathcal{A}_\omega$ of positive associated almost complex
structures.}

\vspace {3mm}

{\bf Proof.} As $P$ anticommutates with $J_0$, $J=J_0e^P$ is
a.c.s. It follows from symmetry of $P$, that $e^P$ is symmetrical
and positive, therefore $J=J_0e^P$ is the positive associated
almost complex structure. Conversely, if
$J\in\mathcal{A}_\omega$, $J=J_0(1+P)(1-P)^{-1}$ and the operator
$(1+P)(1-P)^{-1}$ is positive and symmetrical. Therefore the
logarithm is uniqe defined to it: a symmetrical operator $Q$,
such that $e^Q=J_0(1+P)(1-P)^{-1}$.

\vspace {3mm}

It follows from a proposition 3.4, that one more parametrization
of space $\mathcal{A}_\omega$ is given by map
$$
E_S:\ \mathrm{End}_{SJ_0}(TM)\longrightarrow\mathcal{A}_\omega,
\quad P\mapsto J=J_0e^P.\eqno {(3.12)}
$$

Consider the matter of difference of the positive associated
almost complex structures.

Let $J_0$ is some positive associated almost complex structure,
and $g_0$ is corresponding metric. The space $\mathcal{A}$ of all
almost complex structures is parametrized by endomorphisms
$K:TM\to TM$, anticommutating with $J_0$. As $J_0^T=-J_0$, it
follows from equality $KJ_0=-J_0K$, that
$$
J_0K^T=-K^TJ_0.
$$
Therefore operator $K$ is decomposed in a sum $K=P+L$ of
symmetrical $K$ and skew-symmetric $L$ endomorphisms, each of
which also anticommutates with $J_0$,
$$
\mathrm{End}_{J_0}(TM)=\mathrm{End}_{SJ_0}(TM)\oplus\mathrm{End}_{AJ_0}(TM).
$$

At an exponential parametrization of space $\mathcal{A}$
$$
E:\ \mathrm{End}_{J_0}(TM)\longrightarrow\mathcal{A},\quad
K\mapsto J=J_0e^K,
$$
the subspace $\mathrm{End}_{SJ_0}(TM)$ of symmetrical
endomorphisms parametrizes associated almost complex structures,
and subspace $\mathrm{End}_{AJ_0}(TM)$ of antisymmetric
endomorphisms is used for a parametrization of the other, which
are not associated almost complex structures.

Thus, the submanifold, which is transverse to
$\mathcal{A}_\omega$ is parametrized by map
$$
E_A:\ \mathrm{End}_{AJ_0}(TM)\longrightarrow\mathcal{A},\quad
L\mapsto J=J_0e^L.
$$
As the endomorphism $L$ is skew-symmetric, then $e^L$ is
orthogonal transformation, \\ anticommutating with $J_0$.

Therefore, the submanifold, which is transverse to
$\mathcal{A}_\omega$ in a neighbourhood of the element $J_0$
forms an orthogonal almost complex structures $J$ of the view
$J=J_0O$, where $O$ is orthogonal transformation, anticommutating
with $J_0$.

Finaly, give a parametrization of the space of the associated
metrics. As we know, there is a natural diffeomorphism between
positive associated almost complex structures and associated
metrics:
$$
\mathrm{G}:\ \mathcal{A}_\omega\longrightarrow\mathcal{AM},
$$
$$
J\longrightarrow G(J)=g,\qquad g(X,Y)=\omega (X,JY).
$$
As $J=J_0(1+P)(1-P)^{-1}$,
$$
g(X,Y)=\omega (X,JY)=\omega (X,J_0(1+P)(1-P)^{-1}Y)=
$$
$$
=g_0(X,(1+P)(1-P)^{-1}Y).
$$
We obtain a global parametrization of the space $\mathcal{AM}$ of
the associated metrics:
$$
\Psi_{AM}:\ \mathcal{P}_{J_0}\longrightarrow\mathcal{AM},\quad
P\to g=g_0(1+P)(1-P)^{-1},\eqno {(3.13)}
$$
$$
g(X,Y)=g_0(X,(1+P)(1-P)^{-1}Y).
$$
Other parametrization of the space $\mathcal{AM}$ is given by map
$$
E_{AM}:\ \mathrm{End}_{SJ_0}(TM)\longrightarrow\mathcal{AM},\quad
P\mapsto g=g_0e^P,\eqno {(3.14)}
$$
$$
g(X,Y)=g_0(X,e^PY).
$$

\vspace {3mm}

Find an expression of differential of the mapping $\Psi_{AM}$. As
$\mathcal{P}_{J_0}$ is domain in the space
$\mathrm{End}_{SJ_0}(TM)$, then
$T_P\mathcal{P}_{J_0}=\mathrm{End}_{SJ_0}(TM)$. Therefore,
differential $d\ \Psi_{AM}$ at a point $P$ is mapping of the
following spaces:
$$
d \Psi_{AM}: \mathrm{End}_{SJ_0}(TM)\longrightarrow
T_g\mathcal{AM}.
$$

For the fixed element $P\in\mathcal{P}_{J_0}$ and any symmetric
operator $A$, anticommutating with $J_0$, we consider the line
$P_t=P+tA$ on domain $\mathcal{P}_{J_0}$. Then
$$
g_t=g_0(1+P_t)(1-P_t)^{-1},\quad J_t=J_0(1+P_t)(1-P_t)^{-1}.
$$
Let $h_A=\left.\frac{d}{dt}\right|_{t=0}g_t$ and
$K_A=\left.\frac{d}{dt}\right|_{t=0}J_t$. To find these values
the following obvious equality is used:
$$
\left.\frac{d}{dt}\right|_{t=0}(1-P_t)^{-1}=(1-P)^{-1}A(1-P)^{-1}.
$$
Then
$$
d \Psi_{AM}(P;\
A)=h_A=g_0\left(A(1-P)^{-1}+(1+P)(1-P)^{-1}A(1-P)^{-1}\right).
\eqno {(3.15)}
$$
The expression obtained above can be transformed by three means.

\vspace {3mm}

1) We take into account, that $g=g_0(1+P)(1-P)^{-1}$, then
$$
h_A=g_0\ A(1-P)^{-1}+g\ A(1-P)^{-1}.\eqno {(3.16)}
$$

2) Transform (3.15) as follows:
$$
g_0\left(A(1-P)^{-1}+(1+P)(1-P)^{-1}A(1-P)^{-1}\right)=
g_0\left(1+(1+P)(1-P)^{-1}\right)A(1-P)^{-1}=
$$
$$
=g_0\left(1-P+1+P\right)(1-P)^{-1}A(1-P)^{-1}=2g_0(1-P)^{-1}A(1-P)^{-1}.
$$
Then
$$
h_A=2g_0(1-P)^{-1}A(1-P)^{-1}.\eqno {(3.17)}
$$

3) Instead of $g_0$ it is more convenient to have $g$,
$$
g_0\left(A(1-P)^{-1}+(1+P)(1-P)^{-1}A(1-P)^{-1}\right)=
g_0\left(1+(1+P)(1-P)^{-1}\right)A(1-P)^{-1}=
$$
$$
=g_0(1+P)(1-P)^{-1}\left((1-P)(1+P)^{-1}+1\right)A(1-P)^{-1}=
$$
$$
=g_0(1+P)(1-P)^{-1}\left(1-P+1+P\right)(1+P)^{-1}A(1-P)^{-1}=
$$
$$
=2g(1+P)^{-1}A(1-P)^{-1}.
$$
Thus,
$$
h_A=2g(1+P)^{-1}A(1-P)^{-1}.\eqno {(3.18)}
$$
Write an operator $(1+P)^{-1}$ as
$(1-P)(1-P)^{-1}(1+P)^{-1}=(1-P)(1-P^2)^{-1}$, then
$$
h_A=2g(1-P)(1-P^2)^{-1}A(1-P)^{-1}.\eqno {(3.19)}
$$
We shall accept the last formula as the basic expression of
differential
$$
d \Psi_{AM}: \mathrm{End}_{SJ_0}(TM)\longrightarrow
T_g\mathcal{AM},
$$
$$
d\ \Psi_{AM}(A)=h_A=2g(1-P)(1-P^2)^{-1}A(1-P)^{-1}.
$$
It is easily found the inverse mapping:
$$
d \Psi_{AM}^{-1}(h)=A=\frac 12(1-P)^{-1}(1-P^2)g^{-1}h(1-P).\eqno
{(3.20)}
$$

In a case of an almost complex structure
$J_t=J_0(1+P_t)(1-P_t)^{-1}$ is similarly obtained for an
operator $K_A=\left.\frac{d}{dt}\right|_{t=0}J_t$:
$$
K_A=J_0\ A(1-P)^{-1}+J\ A(1-P)^{-1},\eqno {(3.21)}
$$
$$
K_A=2J_0(1-P)^{-1}A(1-P)^{-1},\eqno {(3.22)}
$$
$$
K_A=2J(1+P)^{-1}A(1-P)^{-1},\eqno{(3.23)}
$$
$$
K_A=2J(1-P)(1-P^2)^{-1}A(1-P)^{-1}.\eqno{(3.24)}
$$

\vspace {3mm}

{\bf 3.3. A complex structure of the space $\mathcal{AM}$.} The
space $\mathcal{AM}$ has a natural almost complex structure,
which is constructed as follows.

The tangent space $T_g\mathcal{AM}$ at $g\in\mathcal{AM}$
consists of all symmetric $J$-anti-Hermitian 2-forms $h$ on $M$,
where $J$ is almost complex structure corresponding to the metric
$g$. As the form $h$ is anti-Hermitian, i.e. $h(JX,JY)=-h(X,Y)$,
the 2-form $hJ$, defined by equality $(hJ)(X,Y)=h(X,JY)$, is also
symmetric and anti-Hermitian. Therefore on each tangent space
$T_g\mathcal{AM}$ the operator acts
$$
{\bf J}_g:\ T_g\mathcal{AM}\longrightarrow T_g\mathcal{AM},\quad
{\bf J}_g(h)=hJ.\eqno{(3.25)}
$$
It is obvious, that ${\bf J}_g^2=-1$. Therefore, on the manifold
$\mathcal{AM}$ the almost complex structure ${\bf J}$ is defined.

On the other hand, the model space $\mathrm{End}_{SJ_0}(TM)$ of a
global parametrization
$$
\Psi_{AM}:\ \mathcal{P}_{J_0}\longrightarrow\mathcal{AM},\quad
P\to g=g_0(1+P)(1-P)^{-1},
$$
has a complex structure:
$$
\mathrm{End}_{SJ_0}(TM)\longrightarrow\mathrm{End}_{SJ_0}(TM),\quad
A\mapsto AJ_0.\eqno {(3.26)}
$$
Therefore one can think, that the space $\mathcal{AM}$ is
infinite-dimensional complex manifold.

\vspace {3mm}

{\bf Theorem 3.1.} {\it The almost complex structure ${\bf J}$ on
the manifold $\mathcal{AM}$ is integrable. The corresponding
complex structure also coincides with a complex structure on
$\mathcal{AM}$, obtained by a parametrization $\Psi_{AM}:\
\mathcal{P}_{J_0}\longrightarrow\mathcal {AM} $.}

\vspace {3mm}

{\bf Proof.} $d \Psi_{AM}(A)=h_A=2g(1-P)(1-P^2)^{-1}A(1-P)^{-1}$,
$$
A\mapsto h_A\mapsto
h_AJ=2g(1-P)(1-P^2)^{-1}A(1-P)^{-1}(1-P)J_0(1-P)^{-1}=
$$
$$
=2g(1-P)(1-P^2)^{-1}AJ_0(1-P)^{-1}\mapsto AJ_0.
$$

\vspace {3mm}

Remark, that the weak Riemannian structure on $\mathcal{AM}$ is
Hermitian with respect to the complex structure ${\bf J}$.
Indeed, if $a,b\in T_g\mathcal{AM}$ are arbitrary tangent
elements, then operators corresponding to them, $A=g^{-1}a,
B=g^{-1}b$ anticommutate with $J$ and we obtain
$$
\left({\bf J}(a),{\bf
J}(b)\right)_g=\left(aJ,bJ\right)_g=\int_M\tr(AJBJ)d\mu=
\int_M\tr(AB)d\mu=\left(a,b\right)_g.
$$

Define the fundamental form of the Hermitian weak Riemannian
structure $\left(a,b\right)_g$ on $\mathcal{AM}$,
$$
\Omega_g\left(a,b\right)=\left(aJ,b\right)_g=\int_M\tr(AJB)d\mu.\eqno
{(3.27)}
$$
It is obvious, that it is the nondegenerate skew-symmetric 2-form
on $\mathcal{AM}$.

\vspace {3mm}

{\bf Theorem 3.2.} {\it The fundamental form $\Omega_g$ of the
Hermitian weak Riemannian structure $\left(a,b\right)_g$ on
$\mathcal{AM}$ is closed.}

\vspace {3mm}

{\bf Proof.} Show, that the exterior differential of the form
$\Omega_g$ is equal to zero at an arbitrary point
$g_0\in\mathcal{AM}$, $d\Omega_{g_0}=0$. For this purpose we use
coordinates on $\mathcal{AM}$ with origin at a point $g_0$:
$g=g_0(1+P)(1-P)^{-1}$. The value of an operator $P=0$
corresponds to the point $g_0$. Therefore it is enough to show,
that $d\Omega_{P}=0 $ at $P=0$. We use the standard formula for
an external product:
$$
d\Omega(A_0,A_1,A_2) =
A_0\Omega(A_1,A_2)-A_1\Omega(A_0,A_2)+A_2\Omega(A_0,A_1)-
$$
$$
-\Omega([A_0,A_1],A_2)+ \Omega([A_0,A_2],A_1) +
\Omega([A_1,A_2],A_0).
$$
Let $A_0,A_1,A_2$ are constant vector fields of operators on the
space $\mathrm{End}_{SJ_0}(TM)$. Then all Lie brackets are equal
to zero. Let's show, that the other addends $A_i\Omega(A_j,A_k)$
are equal to zero too.

The field $h_A$ on $\mathcal{AM}$ corresponds to operator
$A\in\mathrm{End}_{SJ_0}(TM)$ under the following formula:
$$
A\mapsto h_A=2g(1+P)^{-1}A(1-P)^{-1}.
$$
Then ${\bf J}(h_A)=h_{AJ_0}$. Therefore
$$
g^{-1}{\bf
J}(h_A)=g^{-1}h_{AJ_0}=2(1+P)^{-1}AJ_0(1-P)^{-1}=H_{AJ_0}.
$$
We obtain an expression of $\Omega$ in the coordinate map:
$$
\Omega_P(A,B)=\Omega_g(h_A,h_B)=\left(h_AJ,h_B\right)_g=\left(h_{AJ_0},h_B
\right)_g=
$$
$$
=4\int_M\tr(H_{AJ_0}H_B)d\mu=4\int_M\tr\left((1+P)^{-1}AJ_0(1-P)^{-1}(1+P)^{-1}
B(1-P)^{-1}\right)d\mu=
$$
$$
=4\int_M\tr\left((1-P^2)^{-1}AJ_0(1-P^2)^{-1}B\right)d\mu.
$$
Let in this formula $A=A_1,\ B=A_2$ are constant operators (i.e.
do not depend on $P$). On a linear property of an integral and
trace, it is enough to differentiate the expression
$(1-P^2)^{-1}A_1J_0(1-P^2)^{-1}A_2$ on $P$ to find derivative
$A_0\Omega (A_1,A_2)$. One can think, that $P_t=tA_0$. As
$\left.\frac{d}{dt}\right|_{t=0}(1-P_t^2)^{-1}=\left.\frac{d}{dt}
\right|_{t=0}(1-t^2A_0^2)^{-1}=0$, $A_0\Omega
(A_1,A_2)=\left.\frac{d}{dt}
\right|_{t=0}\Omega_{P_t}(A_1,A_2)=0$. The theorem is proved.

\vspace {3mm}

{\bf Corollary.}{\it Manifold $\mathcal{AM}$ is K\"{a}hler.}

\vspace {3mm}

{\bf 3.4. Local expressions. The Beltrami equation.} Let $J_0$ is
positive associated almost complex structure and $g_0$ is
corresponding associated metric. The almost complex structure
$J_0$ defines decomposition of the complexification $TM^C$ of the
tangent bundle $TM$,
$$
TM^C=T^{10}(J_0)\oplus T^{01}(J_0),
$$
on subbundles $T^{10}(J_0)$ and $T^{01}(J_0)$, on which the
complexified operator $J_0$ acts as multiplication on $i$ and
$-i$ respectively.

Let $\partial_1,\dots,\partial_n$  is local basis of sections of
the bundle $T^{10}(J_0)$,\ $dz^1,\dots,dz^n$ is dual basis of the
bundle $T^{*10}(J_0)$ and
$\overline{\partial}_1,\dots,\overline{\partial}_n$ is complex
conjugate basis of sections of the bundle $T^{01}(J_0)$,\
$d\overline{z}^1,\dots,d\overline{z}^n$ is dual basis of
$T^{*01}(J_0)$.

As $g_0$ is $J_0$-Hermitian metric, it follows, that
$$
g_{\alpha\beta}=g_0(\partial_\alpha,\partial_\beta)=0, \quad
g_{\overline{\alpha}\overline{\beta}}=g_0(\overline{\partial}_\alpha,
\overline{\partial}_\beta)=0,\quad\alpha,\beta =1,\dots, n.
$$
Let
$$
g_{\alpha\overline{\beta}}=g_0(\partial_\alpha,\overline{\partial}_\beta),
\quad
g_{\overline{\alpha}\beta}=g_0(\overline{\partial}_\alpha,\partial_\beta),
\quad\alpha,\beta=1,\dots, n.
$$
These coefficients have properties
$$
g_{\alpha\overline{\beta}}=g_{\overline{\beta}\alpha}, \quad\
g_{\alpha\overline{\beta}}=\overline{g_{\beta\overline{\alpha}}},
$$
which are implied from symmetry and hermiticity of the metric
$g_0$.

The metric $g_0$ is expressed as follows:
$$
g_0=2g_{\alpha\overline{\beta}}\ dz^\alpha d\overline{z}^{\beta}.
$$

Let now $J$ is another positive associated almost complex
structure. Then $J=J_0(1+P)(1-P)^{-1}$, where $P:\ TM\to TM$ is a
symmetric endomorphism, anticommutating with $J_0 $, satisfying
to the condition of positiveness $1-P^2>0$.

As
$$
J_0(P(\partial_\alpha))=-PJ_0(\partial_\alpha)=-P(i\partial_\alpha)=
-iP(\partial_\alpha),
$$
that $P(\partial_\alpha)$ is a local section of the bundle
$T^{01}(J_0)$, therefore
$$
P(\partial_\alpha)=P_{\alpha}^{\overline{\beta}}\
\overline{\partial}_\beta, \quad\alpha,\beta =1,\dots,n\eqno
{(3.28)}
$$
where $P_{\alpha}^{\overline{\beta}}$ is matrix of complex-valued
functions. Thus, the operator $P$ in the complex basis is locally
given by a matrix of view:
$$
P=\left(\begin{array}{cc}
0 & P_{\overline{\alpha}}^{\beta}\\
P_{\alpha}^{\overline{\beta}} & 0
\end{array}
\right),\quad P_{\overline{\alpha}}^{\beta}=
\overline{P_{\alpha}^{\overline{\beta}}}.\eqno {(3.29)}
$$
The condition of symmetry of an operator $P$ is expressed by the
relation:
$$
P_{\alpha\beta}=P_{\beta\alpha},\quad\mbox{where}\quad
P_{\alpha\beta} =g_{\alpha\overline{\gamma}}\
P_{\beta}^{\overline{\gamma}}.\eqno {(3.30)}
$$

Invariant form of an operator $P$:
$$
P=P_{\alpha}^{\overline{\beta}}\ \overline{\partial}_\beta\otimes
dz^\alpha +P_{\overline{\alpha}}^{\beta}\ \partial_\beta\otimes
d\overline{z}^\alpha \eqno {(3.31)}.
$$

Positive associated almost complex structure
$J=J_0(1+P)(1-P)^{-1}$ defines another decomposition of the
complexification $TM^C$,
$$
TM^C=T^{10}(J)\oplus T^{01}(J).
$$
Dependence of this expansion upon $J$ becomes clear if we use an
operator $P$.

\vspace {3mm}

{\bf Proposition 3.5.} {\it If $J=J_0(1+P)(1-P)^{-1}$ then
$$
T^{10}(J)=(1-P)(T^{10}(J_0)),\quad
T^{01}(J)=(1-P)(T^{01}(J_0)).\eqno {(3.32)}
$$
Vector fields
$$
\partial_\alpha(J)=\partial_\alpha-P_{\alpha}^{\overline{\beta}}\
\overline{\partial}_\beta,\quad\overline{\partial}_\alpha(J)=
\overline{\partial}_\alpha-P_{\overline{\alpha}}^{\beta}\
\partial_\beta \eqno {(3.33)}
$$
form local basises of sections of bundles $T^{10}(J)$ and
$T^{01}(J)$ respectively.}

\vspace {3mm}

{\bf Proof.} It is enough to prove, that the vector fields
$\partial_\alpha(J)=\partial_\alpha-P_{\alpha}^{\overline{\beta}}
\overline{\partial}_\beta$ form the basis of the bundle
$T^{10}(J)$. From a nondegeneracy of an endomorphism $1-JJ_0$ the
nondegeneracy of an operator $1-P$, which is equal to
$-2(1-JJ_0)^{-1}JJ_0$, follows. Therefore, vector fields
$\partial_\alpha(J)=(1-P)(\partial_\alpha)$ are linearly
independent. Let's show, that $\partial_\alpha(J)$ are sections of
the bundle $T^{10}(J)$, i.e. that
$J(\partial_\alpha(J))=i\partial_\alpha(J)$,
$$
J(\partial_\alpha(J))=J(1-P)\partial_\alpha=J_0(1+P)(1-P)^{-1}(1-P)
\partial_\alpha=
$$
$$
=J_0(1+P)\partial_\alpha=(1-P)J_0\partial_\alpha=i(1-P)\partial_\alpha
=i\partial_\alpha(J).
$$

\vspace {3mm}

{\bf Corollary.} {\it The dual basises of the forms
$dz_1(J),\dots,dz_n(J)$ and
$d\overline{z}_1(J),\dots,d\overline{z}_n(J)$ for basises of
vector fields $\partial_1(J),\dots,\partial_n(J)$ and
$\overline{\partial}_1(J),\dots,\overline{\partial}_n(J)$ have the
following expressions:
$$
dz^{\alpha}(J)=D_{\gamma}^{\alpha}\left(dz^{\gamma}+
P_{\overline{\beta}}^{\alpha}\ d\overline{z}^{\beta}\right),\eqno
{(3.34)}
$$
$$
d\overline{z}^{\alpha}(J)=D_{\overline{\gamma}}^{\overline{\alpha}}
\left(d\overline{z}^{\gamma}+P_{\beta}^{\overline{\alpha}}\
dz^{\beta}\right), \eqno {(3.35)}
$$
where $D_{\gamma}^{\alpha}$ is matrix of an operator
$D=(1-P^2)^{-1}:\ T^{10}(J_0)\ \longrightarrow T^{01}(J_0)$.}

\vspace {3mm}

It follows from the proposition 3.5 and formula (3.31), that the
operator $P$, giving an almost complex structure $J$ is
many-dimensional generalization of Beltrami coefficient in the
Beltrami equation \cite{Ear-Eel}
$$
\frac{\partial f}{\partial\overline{z}}-\mu\frac{\partial
F}{\partial z}=0.
$$
Indeed, geometrical sense (see for example \cite{Ear-Eel} ) of
Beltrami coefficient is that $\mu$ is tensor field on a Riemann
surface $M$ of view: $\mu=\mu\frac{\partial}{\partial z}\otimes
d\overline{z}$. The invariant form (3.31) of the operator $P$ has
the same sense. Therefore many-dimensional generalization of the
Beltrami equation has the following view:
$$
\overline{\partial}_\alpha(J)(f)=\frac{\partial f}
{\partial\overline{z}^{\alpha}}-P_{\overline{\alpha}}^{\beta}\
\frac{\partial f}{\partial z^{\beta}}=0,\quad\alpha
=1,\dots,n,\eqno {(3.36)}
$$
where $f$ is function on $M$.

\vspace {3mm}

Each positive associated a.c.s. $J$ defines associated Hermitian
metric $g$ by equality $g(X,Y)=\omega (X,JY)$. Recall the
expression of the metric $g$ via the Riemannian metric $g_0$ and
operator $P$:
$$
g(X,Y)=g_0((1+P)X,(1-P)^{-1}Y)=g_0((1+P)X,(1+P)DY),
$$
where $D=(1-P^2)^{-1}$.

Assume, for a simplicity that a.c.s. $J_0$ is integrable and
$z^1,\dots,z^n$ are corresponding complex local coordinates on
$M$. Let in these coordinates the Hermitian form corresponding to
a Riemannian structure $g_0$ has an aspect
$g_0=2g_{\alpha\overline{\beta}}dz^\alpha d\overline{z}^{\beta}$.
Then from (3.34) and (3.35) is obtained the following local
expression for the Hermitian form of the associated metric $g$:
$$
g=2g_{\alpha\overline{\beta}}\
D_{\overline{\gamma}}^{\overline{\beta}}
\left(dz^\alpha+P_{\overline{\sigma}}^{\overline{\alpha}}\
d\overline{z}^{\sigma}\right)\left(d\overline{z}^{\gamma}+
P_{\nu}^{\overline{\gamma}}\ dz^{\nu}\right).\eqno {(3.37)}
$$

\newpage

\vspace {10mm}

\centerline {\bf \S 4. Decomposition of the space of Riemannian
metrics} \centerline {\bf on a symplectic manifold.}

\vspace {7mm}

Let $M^{m}$ be a smooth oriented manifold and
$\mathrm{Vol}(M)\subset\Gamma (\Lambda^{m}M)$ is the space of
smooth volume forms on $M$, i.e. the space of smooth
nondegenerate $m$-forms on $M$, defined orientation which
coincides with the original one on $M$. The natural projection is
defined
$$
vol:{\cal M}\longrightarrow\mathrm{Vol}(M),\qquad
g\longmapsto\mu_g= \sqrt{\det g}dx^1\wedge\ldots\wedge dx^{m}.
$$
A fiber of the bundle $vol$ over $\mu\in Vol(M)$ is the space
$\mathcal{M}_{\mu}$ of metrics with the same Riemannian volume
form $\mu$.

The fixing of a volume form $\mu$ defines decomposition of the
space ${\cal M}$ in the direct product:
$$
\varphi_{\mu}:{\cal M\longrightarrow}\mathrm{Vol}(M)\times{\cal
M}_\mu, \quad g\longmapsto (\mu_g,\rho^{-\frac 2m} g),
$$
where $\rho$ is density of the volume form $\mu_g$ with respect
to $\mu$, i.e. function on $M$ is defined from equality
$\mu_g=\rho\mu$.

The inverse mapping:
$$
\iota_{\mu}:\mathrm{Vol}(M)\times{\cal M}_\mu\longrightarrow{\cal
M}, \quad (\nu,h)\longmapsto\rho^{\frac 2m}h,\quad\nu=\rho\mu.
$$

At such decomposition, the space $\mathrm{Vol}(M)$ of volume
forms corresponds to the space of the metrics $C_g$, which are
conformally equivalent to the fixed metric $g\in{\cal M}_\mu$:
$$
\mathrm{Vol}(M)\times\{g\}\longrightarrow
C_g,\quad\nu=\rho\mu_g,\quad\nu \longmapsto\rho^{\frac 2m}g.
$$

Thus, the space of all Riemannian metrics on a manifold $M$ is
decomposed in a direct product of the space of metrics with a
fixed Riemannian volume form and the space of pointwise
conformally equivalent metrics. The similar construction is
possible in case of symplectic and contact manifolds.

Let $M^{2n},\omega$ is symplectic manifold.

We remind, that the almost complex structure $J$ on $M$  is
called {\it positive associated} to the symplectic form $\omega$,
if for any vector fields $X,Y$ on $M$,

1) $\omega (JX,JY)=\omega (X,Y)$,

2) $\omega (X,JX)>0$, if $X\neq 0$.

Every such a.c.s. $J$ defines the Riemannian metric $g$ on $M$ by
equality:
$$
g(X,Y)=\omega (X,JY),
$$
which is also called {\it associated}.

Let ${\cal AM}_\omega$ is the space of all smooth associated
metrics on $M$.

The symplectic form $\omega$ on a manifold $M$ defines a well
known projection of the space ${\cal M}$ on ${\cal AM}_\omega$ as
follows:

Let $g'\in{\cal M}$ is any metric. There is a unique
skew-symmetric automorphism $A$ of a tangent bundle $TM$, such
that:
$$
\omega (X,Y)=g'(AX,Y),\quad A^T=-A.
$$
Apply the polar decomposition to the endomorphism $A$
$$
A=JH,
$$
where $J$ is orthogonal operator and $H$ is positive symmetrical
one. The endomorphism $-A^2=A^TA=AA^T$ is positive defined and
symmetrical with respect to $g'$. Then $H=(-A^2)^{\frac 12}$ is
positive square root from $-A^2$. As it is known, it commutes
with operators $A$ and $J$. Suppose $J=A(-A^2)^{-\frac 12}$. It
is easily checked, that $J$ is an almost complex structure:
$$
J^2=A(-A^2)^{-\frac 12}A(-A^2)^{-\frac 12}=A^2(-A^2)^{-1}=-I.
$$
The formula
$$
g(X,Y)=\omega (X,JY)
$$
defines the Riemannian metric on $M$. Indeed,
$$
g(X,Y)=\omega (X,JY)=g'(AX,JY)=g'(AX,A(-A^2)^{-\frac 12}Y)=
$$
$$
=g'(X,(-A^2)(-A^2)^{-\frac 12}Y)=g'(X,(-A^2)^{\frac 12}Y).
$$
As the operator $(-A^2)^{\frac 12}$ is positive and symmetrical,
then $g$ is the Riemannian metric. Positiveness of an almost
complex structure $J$ also follows from here:\ $\omega
(X,JX)=g(X,X)>0$.

The metric $g$ is $J$-Hermitian:
$$
g(JX,JY)=\omega (JX,-Y)=-g'(AJX,Y)=-g'(AA(-A^2)^{-\frac 12}X,Y)=
$$
$$
=g'((-A^2)^{\frac 12}X,Y)=g'(X,(-A^2)^{\frac 12}Y)=g(X,Y).
$$
The symplectic form $\omega$ is also $J$-invariant:
$$
\omega (JX,JY)=g(JX,Y)=-g(X,JY)=-\omega (X,J^2Y)=\omega (X,Y).
$$

Therefore, $J$ is a positive associated almost complex structure,
and metric $g(X,Y)=\omega (X,JY)$ is associated metric, which
corresponds to the structure $J$.

We have got a required projection
$$
p_\omega:\ {\cal M}\longrightarrow{\cal AM}_\omega,\quad
g'\mapsto g, \qquad g(X,Y)=g'(X,(-A^2)^{\frac 12}Y).\eqno {(4.1)}
$$

As the operator $J$ on construction is orthogonal with respect to
$g'$, the metric $g'$ is also Hermitian with respect to the
a.c.s. $J$. One can show, that the fiber $p_\omega^{-1}(g)$
consists of all $J$-Hermitian metrics. Let ${\cal M}_J$ denotes
the set of all $J$-Hermitian Riemannian metrics on $M$.

\vspace{3mm}

{\bf Lemma 4.1.} {\it For any associated metric $g\in{\cal
AM}_\omega$ and its corresponding a.c.s. $J$, the inverse image
of the element $g$ at a projection $p_\omega$ coincides with a
set ${\cal M}_J$ of all $J$-Hermitian Riemannian metrics on $M$:}
$$
p_\omega^{-1}(g)={\cal M}_J.
$$

\vspace{3mm}

{\bf Proof.} We already remarked, that any metric $g'$, which is
projected in $g$ is $J$-Hermitian, therefore
$p_\omega^{-1}(g)\subset{\cal M}_J$. Show the converse. Let
$g'\in{\cal M}_J$ is $J$-Hermitian metric. As $g$ and $g'$ are
$J$-Hermitian metrics, there is a symmetric positive operator $B$
commuting with $J$ and such that $g'(X,Y)=g(X,BY)=g(BX,Y)$.
$$
\omega (X,Y)=g(JX,Y)=g'(JX,B^{-1}Y)=g'(B^{-1}JX,Y).
$$
Therefore $A=B^{-1}J$, $-A^2=-B^{-1}JB^{-1}J=-B^{-2}J^2=B^{-2}$.
As $B^{-1}$ is symmetric, $(-A^2)^{-\frac 12}=(B^{-2})^{-\frac
12}=B$. Then the almost complex structure corresponding to the
metric $g'$ coincides with original: $(-A^2)^{-\frac
12}A=BB^{-1}J=J$. Therefore $p_\omega (g')=g$.

\vspace {3mm}

Let $g_0\in{\cal AM}_\omega$ is some fixed associated metric and
$J_0$ is almost complex structure, corresponding to it. Any other
associated metric $g,\ J$ can be represented as
$$
g(X,Y)=g_0(X,(1+P)(1-P)^{-1}Y),\quad J=J_0(1+P)(1-P)^{-1},
$$
where the operator $P$, anticommutating with an a.c.s. $J_0$, is
symmetric with respect to $g_0$ and $1-P^2$ is positive with
respect to $g_0$.

There is a natural question: Are other metrics of fibers
$\mathcal{M}_{J_0}$ and $\mathcal{M}_J$ connected by similar
relations? The answer is given by the following

\vspace{3mm}

{\bf Lemma 4.2.} {\it Let $g_0^{\prime}\in\mathcal{M}_{J_0}$ is an
arbitrary $J_0$-Hermitian metric and $J=J_0(1+P)(1-P)^{-1}$. Then
$g'$, defined by equality:
$$
g^{\prime}(X,Y)=\frac
12\left(g_0^{\prime}\left((1+P)(1-P)^{-1}X,Y\right)+
g_0^{\prime}\left(X,(1+P)(1-P)^{-1}Y\right)\right)\eqno{(4.2)}
$$
is $J$-Hermitian metric. The fundamental form
$\omega^{\prime}(X,Y)=g^{\prime}(JX,Y)$ of the metric
$g^{\prime}$ is expressed via $\omega$ as follows:}
$$
\omega^{\prime}(X,Y)=\frac 12\left(\omega (X,Y)+\omega\left(
(1+P)(1-P)^{-1}X,(1+P)(1-P)^{-1}Y\right)\right).\eqno {(4.3)}
$$

{\bf Proof.} Symmetry of $g^{\prime}(X,Y)$ follows from the
definition. The positive definiteness of $g^{\prime}(X,Y)$ just
follows from the positive definiteness of operator
$(1+P)(1-P)^{-1}$ with respect to $g_0$. $J$-hermiticity of
metric $g^{\prime}(X,Y)$ at once follows from $J_0$-hermiticity
of $g_0'$ and the following property
$$
J=J_0(1+P)(1-P)^{-1}=(1-P)(1+P)^{-1}J_0.
$$
The fundamental form $\omega^{\prime}(X,Y)=g^{\prime}(JX,Y)$ of
the metric $g^{\prime}$ is found from elementary calculation.

\vspace {3mm}

It follows from this lemma, that for the associated metrics
$g_0,\ J_0$ and $g,\ J$ the fibers $\mathcal{M}_{J_0}$ and
$\mathcal{M}_J$ are naturally diffeomorphic. Indeed, formulas
$$
g_0(X,Y)=g(X,(1-P)(1+P)^{-1}Y),\quad J_0=J(1-P)(1+P)^{-1},
$$
allow to define the inverse mapping of fibers
$\mathcal{M}_J\to\mathcal{M}_{J_0}$:
$$
g_0^{\prime}(X,Y)=\frac
12\left(g^{\prime}\left((1-P)(1+P)^{-1}X,Y\right)+
g^{\prime}\left(X,(1-P)(1+P)^{-1}Y\right)\right).\eqno {(4.4)}
$$

{\bf Theorem 4.1.} {\it The space $\mathcal{M}$ is a smooth
trivial bundle over $\mathcal{AM}_\omega$. A fiber over the
element $(g,J)\in\mathcal{AM}_\omega$ is the space
$\mathcal{M}_J$ of all $J$-Hermitian Riemannian metrics on $M$.}

\vspace{3mm}

{\bf Proof.} Let $(g_0,J_0)\in\mathcal{AM}_\omega$ is associated
structure. The space $\mathcal{AM}_\omega$ is parametrized by
domain
$$
\mathcal{P}_{J_0}=\{P\in\mathrm{End}_{SJ_0}(TM):\ 1-P^2>0\}
$$
of the space $\mathrm{End}_{SJ_0}(TM)$ of smooth symmetrical
endomorphisms $P:\ TM\to TM$ \\ anticommutating with $J_0$.

The fiber $\mathcal{M}_{J_0}$ of fiber bundle $p_{\omega}$ over
the point $g_0$ consists of all $J_0$-Hermitian metrics, it is an
open set in the linear Frechet space of smooth symmetric
$J_0$-Hermitian forms on $M$.

One can take the following ILH-smooth mapping as a required map :
$$
\Psi_{M}:\
\mathcal{P}_{J_0}\times\mathcal{M}_{J_0}\longrightarrow\mathcal{M},
\quad (P,g_0^{\prime})\longrightarrow g^{\prime},
$$
where $P\in\mathcal{P}_{J_0}$, $g_0^{\prime}\in\mathcal{M}_{J_0}$
and the metric $g^{\prime}$ is defined by equality (4.2). It
follows from the lemma 2, that $\mathcal{M}_{J_0}$ is
diffeomorphic mapped on fiber $\mathcal{M}_J$ of fiber bundle
$p_{\omega}$ over the point $g(X,Y)=g_0((1+P)X,(1-P)^{-1}Y)$.

The inverse mapping for $\Psi_{M}$ is defined by the
correspondence (4.4):\ $\mathcal{M}_J\to\mathcal{M}_{J_0}$,
$g^{\prime}\longrightarrow g_0^{\prime}$ together with the
projection
$p_{\omega}:\mathcal{M}\longrightarrow\mathcal{AM}_\omega$. The
theorem is proved.

\vspace{3mm}

{\bf Conclusion.} {\it The space of all Riemannian metrics on a
symplectic manifold $M$ is decomposed in a direct product of the
space of associated metrics and the space of $J$-conformally
equivalent metrics.}

\vspace {3mm}

Finaly, we remark, that one more projection
$q:\mathcal{M}\longrightarrow\Lambda_0^2(M)$ of the space of
metrics $\mathcal{M}$ on the space $\Lambda_0^2(M)$ of all smooth
nondegenerate skew-symmetric 2-forms can be defined. If
$g^{\prime}\in\mathcal{M}$ and $J=p_\omega (g^{\prime})$ is
corresponding a.c.s., then
$$
q(g^{\prime})=\omega^{\prime},\quad\omega^{\prime}(X,Y)=g^{\prime}(JX,Y).
$$

The inverse image $q^{-1}(\omega^{\prime})$ of the element
$\omega^{\prime}\in\Lambda_0^2(M)$ is a set
$\mathcal{AM}_{\omega^{\prime}}$ of Riemannian metrics,
associated with the nondegenerate 2-form $\omega^{\prime}$ on $M$.

At mapping $q$ the fiber $\mathcal{M}_J$ of projection $p_\omega$
passes in the space of the nondegenerate exterior 2-forms, which
are invariant with respect to a.c.s. $J$, i.e. such forms
$\omega^{\prime}$, that
$\omega^{\prime}(JX,JY)=\omega^{\prime}(X,Y)$.

\newpage

\vspace {10mm}

\centerline{\bf \S 5. A curvature of the space of associated
metrics.}

\vspace {7mm}

In this paragraph we shall find geodesics and sectional
curvatures of the space $\mathcal{AM}$ of almost K\"{a}hler
metrics on a symplectic manifold $M^{2n},\omega$.

It is well known \cite{Ko-No}, that the Riemannian volume form
$\mu_g$ of a metric $g\in\mathcal{AM}$ is expressed via the form
$\omega$: $\mu_g=\frac{1}{n!}\omega^n$. Therefore, the space
$\mathcal{AM}$ is into manifold $\mathcal{M}_{\mu}$ of metrics
with the same volume form $\mu=\frac{1}{n!}\omega^n$.

The manifold $\mathcal{M}_{\mu}$, is a smooth ILH-submanifold in
$\mathcal{M}$ \cite{Ebi1} and inherits a weak Riemannian
structure, which is expressed as follows. If $a,b\in
T_g\mathcal{M}_{\mu}$ are two smooth symmetric 2-forms on $M$,
representing elements of tangent space $T_g\mathcal{M}_{\mu}$,
then their inner product is defined by the formula:
$$
(a,b)_g=\int_M\tr AB\ d\mu,\eqno {(5.1)}
$$
where $A=g^{-1}a=g^{ik}a_{kj}$. Recall, that the tangent space
$T_g\mathcal{M}_\mu$ consists of traceless symmetric 2-forms,
$$
T_g\mathcal{M}_\mu=S_2^T=\{h\in S_2:\ \tr_g\ h=0\}.
$$
As in (5.1) the volume form does not depend on $g$, $\mu_g=\mu$,
it will be more convenient to use another weak Riemannian
structure, defined by the formula (2.9) from \S 2, on all the
space $\mathcal{M}$. In this case submanifold
$\mathcal{M}_{\mu}\subset\mathcal{M}$ is (see \S 2) totally
geodesic in $\mathcal{M}$. Therefore, from the theorem 2.2 we
obtain the following characteristics of the space
$\mathcal{M}_{\mu}$:

1) A covariant derivative:
$$
\nabla_ab=d_ab-\frac 12\left(aB+bA\right),
$$

2) A tensor of curvature:
$$
R(a,b)c=-\frac 14g\ [[A,B],C],
$$

3) A sectional curvature $K_{\sigma}$ in plane section $\sigma$,
given by orthonormal pair $a,b\in T_g\mathcal{M}_{\mu}$:
$$
K_{\sigma}=\frac 14\int_M\tr\left([A,B]^2\right)\ d\mu,
$$

4) Geodesics, going out from a point $g\in\mathcal{M}_{\mu}$ in
direction $a\in T_g\mathcal{M}_{\mu}$ look like $g_t=g\ e^{tA}$.

\vspace{3mm}

{\bf Proposition 5.1.} {\it The manifold $\mathcal{AM}$ is totally
geodesic submanifold in the space $\mathcal{M}_\mu$ of all
Riemannian structures on $M$ with the same volume form
$\mu=\frac{1}{n!}\omega^n$.}

\vspace{3mm}

{\bf Proof.} Geodesics on $\mathcal{M}_\mu$ look like
$g_t=ge^{tA}$, where $A=g^{-1}a$, $a\in T_g\mathcal{M}_\mu$.
Show, that if $a\in T_g\mathcal{AM}$, then geodesic $g_t$ is on
$\mathcal{AM}$. If $a\in T_g\mathcal{AM}$, then the operator
$A=g^{-1}a$ is symmetric and anticommutates with $J$, $JA=-AJ$.
Therefore, $Je^{tA}=e^{-tA}J$. Show, that $g_t=ge^{tA}$ is a
family of associated metrics, corresponding to a family of
positive associated almost complex structures $J_t=Je^{tA}$. For
this purpose it is enough to show, that $g_t(X,J_tY)=-\omega
(X,Y)$, but it is obvious:
$$
g_tJ_t=ge^{tA}Je^{tA}=gJe^{-tA}e^{tA}=gJ=-\omega.
$$
The proposition is proved.

\vspace {4mm}

Recall, that
$$
\mathrm{End}_{SJ}(TM)=\{P\in\mathrm{End}(TM);\quad PJ=-JP,\quad
g(PX,Y)=g(X,PY)\}
$$
is the space of smooth symmetrical endomorphisms $P: TM\to TM$,
anticommutating with $J$.

{\bf Corollary.} {\it Map
$$
E_{AM}: \mathrm{End}_{SJ_0}(TM)\longrightarrow\mathcal{AM},\quad
P\mapsto g= g_0e^P,\ g(X,Y)=g_0(X,e^PY)
$$
gives normal coordinates on the space $\mathcal{AM}$ in a
neighbourhood of the element $g_0$.}

\vspace {4mm}

The tangent space $T_g\mathcal{AM}$ at a point $g\in\mathcal{AM}$
consists of all symmetric anti-Hermitian 2-forms $h$ on $M$. As
the form $h$ is anti-Hermitian, i.e. $h(JX,JY)=-h(X,Y)$, the
2-form $hJ$, defined by equality $(hJ)(X,Y)=h(X,JY)$, is also
symmetric and anti-Hermitian. Therefore, on the space
$T_g\mathcal{AM}$ the operator
$$
{\bf J}:\ T_g\mathcal{AM}\longrightarrow
T_g\mathcal{AM},\quad{\bf J}(h)=hJ.
$$
is defined. Obviously, that ${\bf J}^2=-1$. Therefore, on the
space $\mathcal{AM}$ the almost complex structure ${\bf J}$ is
defined. It is immediately checked, that the inner product (5.1)
is Hermitian with respect to the a.c.s. ${\bf J}$ on
$\mathcal{AM}$ (see also \S 3).

It follows from the property of the submanifold $\mathcal{AM}$ to
be totally geodesic in $\mathcal{M}_\mu$, that all properties of
a curvature pointed above for the space $\mathcal{M}_\mu$ hold
for the space $\mathcal{AM}$.

\vspace {3mm}

{\bf Theorem 5.2.} {\it The space $\mathcal{AM}$ has the following
geometric characteristics.

1) A tensor of curvature is following
$$
R(a,b)c=-\frac 14g\ [[A,B],C],\eqno {(5.2)}
$$
where $a,b,c\in T_g\mathcal{AM}$, $A=g^{-1}a$.

2) A sectional curvature $K_{\sigma}$ in plane section $\sigma$,
given by orthonormal pair $a,b\in T_g\mathcal{AM}$ is expressed
by the formula
$$
K_{\sigma}=\frac 14\int_M\tr\left([A,B]^2\right)\ d\mu. \eqno
{(5.3)}
$$
In particular, a holomorphic sectional curvature has a view:
$$
K(a),{\bf J}a)=-\int_M\mbox{tr}(A^4)d\mu.\eqno {(5.4)}
$$

3) Geodesics, going out from a point $g\in\mathcal{AM}$ in
direction $a\in T_g\mathcal{AM}$ are given by the following way:\
$g_t=g\ e^{tA}$.}

\vspace {3mm}

Fix a positive associated a.c.s. $J_0$ and corresponding
associated metric $g_0$. Let $\mathcal{P}_{J_0}$ is domain in the
space $\mathrm{End}_{SJ_0}(TM)$, consisting of endomorphisms
$P:TM\to TM$, for which the operator $1-P^2$ is positively
defined with respect to $g_0$.

Consider a global parametrization of the space $\mathcal{AM}$ of
associated metrics:
$$
\Psi_{AM}: \mathcal{P}_{J_0}\longrightarrow\mathcal{AM},\quad P\to
g=g_0(1+P)(1-P)^{-1},
$$
$$
g(X,Y)=g_0(X,(1+P)(1-P)^{-1}Y).
$$

Find expression of metric and curvature of the space
$\mathcal{AM}$ in these coordinates. The differential of mapping
$\Psi_{AM}$ acts as follows:
$$
d\ \Psi_{AM}:\ \mathrm{End}_{SJ_0}(TM)\longrightarrow
T_g\mathcal{AM},
$$
$$
d\ \Psi_{AM}(A)=h_A=2g(1-P)(1-P^2)^{-1}A(1-P)^{-1}.
$$
The inverse mapping:
$$
d\ \Psi_{AM}^{-1}(h)=A=\frac 12(1-P)^{-1}(1-P^2)g^{-1}h(1-P).
$$
It follows from these formulas and the theorem 5.2

\vspace{3mm}

{\bf Theorem 5.3.} {\it The space $\mathcal{AM}$ has the following
geometric characteristics.

1) An inner product is given by the formula:
$$
(A,B)_P=4\int_M\tr\left((1-P^2)^{-1}A(1-P^2)^{-1}B\right)d\mu,\eqno
{(5.5)}
$$
where $A,B\in T_P\mathcal{P}_{J_0}=\mathrm{End}_{SJ_0}(TM)$.

2) A covariant derivative of vector fields given by (constant)
operators $A$ and $B$:
$$
\nabla_AB=AP(1-P^2)^{-1}B+B(1-P^2)^{-1}A.\eqno {(5.6)}
$$

3) A curvature tensor has a view:
$$
R(A,B)C=-(1-P^2)\left[\left[(1-P^2)^{-1}A,(1-P^2)^{-1}B\right],(1-P^2)^{-1}C
\right],\eqno {(5.7)}
$$
where $A,B,C\in T_P\mathcal{P}_{J_0}=\mathrm{End}_{SJ_0}(TM)$.

4) A sectional curvature $K_{\sigma}$ in a plane section $\sigma$,
given by orthonormal pair $A,B$ is found by the formula:
$$
K_{\sigma}=\int_M\tr\left(\left[(1-P^2)^{-1}A,(1-P^2)^{-1}A\right]^2\right)\
d\mu.\eqno {(5.8)}
$$
In particular, a holomorphic sectional curvature has a view:
$$
K(A,AJ_0)=-\int_M\tr\
\left(\left((1-P^2)^{-1}A\right)^4\right)d\mu.\eqno{(5.4)}
$$

5) Geodesics, going out from a point $g_0$ in directions
$A\in\mathrm{End}_{SJ_ 0}(TM)$ are curves $P(t)$ on domain
$\mathcal{P}_{J_0}$ of a view:
$$
P(t)=\tanh(tA)=\left(e^{tA}+e^{-tA}\right)^{-1}\left(e^{tA}-e^{-tA}\right).
\eqno {(5.10)}
$$
}

{\bf Proof.} 1) The inner product.
$$
(A,B)_P=(h_A,h_B)_g=\int_M\tr\left(g^{-1}h_Ag^{-1}h_B\right)d\mu=
$$
$$
=4\int_M\tr\left((1-P)(1-P^2)^{-1}A(1-P^2)^{-1}B(1-P)^{-1}\right)d\mu=
$$
$$
=4\int_M \tr\left((1-P^2)^{-1}A(1-P^2)^{-1}B\right)d\mu.
$$

2) The covariant derivative is immediately calculated under the
six-term formula. In view of a constancy of operators $A,B,C$ it
is:
$$
(\nabla_AB,C)_P=\frac 12\left(A(B,C)_P+B(A,C)_P-C(A,B)_P\right).
$$
In calculations we use a linearity of an integral and trace, and
also the formula
$$
\left((1-P_t^2)^{-1}\right)'=(1-P^2)^{-1}(AP+PA)(1-P^2)^{-1},
$$
where $P_t=P+tA$ is variation of an operator $P$ in direction $A$.

3) The curvature tensor. Let $h_A,h_B,h_C\in T_g\mathcal{AM}$,
then
$$
R(h_A,h_B)h_C=-\frac 14g\left[\left[H_A,H_B\right],H_C\right],
$$
where $H_A=g^{-1}h_A=2(1-P)(1-P^2)^{-1}A(1-P)^{-1}$. Substituting
this expression, we obtain:
$$
R(h_A,h_B)h_C=-2g(1-P)\left[\left[(1-P^2)^{-1}A,1-P^2)^{-1}B\right],
(1-P^2)^{-1}C\right](1-P)^{-1}.
$$
Therefore
$$
R(A,B)C=d\
\Psi_{AM}^{-1}R(h_A,h_B)h_C=-(1-P^2)\left[\left[(1-P^2)^{-1}A,
(1-P^2)^{-1}B\right],(1-P^2)^{-1}C\right].
$$

4) The sectional curvature is calculated similarly via using of
expression $d \Psi_{AM}(A)=h_A=2g(1-P)(1-P^2)^{-1}A(1-P)^{-1}$.

5) Let $A\in\mathrm{End}_{SJ_0}(TM)$ and $h=d\
\Psi_{AM}(A)=2g_0A$. Then geodesic going out from a point $g_0$
in direction $h\in\mathcal{AM}_{g_0}$ is a curve
$g_t=g_0e^{t2A}$. Let $P(t)$ is corresponding curve on the domain
$\mathcal{P}_{J_ 0}$. Then
$g_t=g_0e^{t2A}=g_0(1+P_t)(1-P_t)^{-1}$. From the last formula
$P(t)$ we express:
$$
P(t)=\left(e^{tA}+e^{-tA}\right)^{-1}\left(e^{tA}-e^{-tA}\right)=\tanh(tA).
$$

\newpage

\vspace {10mm}

\centerline {\bf \S 6. Orthogonal decompositions of the space}
\centerline {\bf of symmetric tensors on an almost K\"{a}hler
manifold.}

\vspace {5mm}

We consider a manifold $M$ with the closed nondegenerate 2-form
$\omega$ of class $C^{\infty}$, $\dim\ M=m=2n$.

Let $\mathcal {M}$ is the space of all Riemannian metrics on $M$.
The group $\mathcal{D}$ of smooth diffeomorphisms of manifold $M$
acts naturaly on the space $\mathcal{M}$:
$$
A: \mathcal{M}\times\mathcal{D}\longrightarrow\mathcal{M},\qquad
A(g,\eta)= \eta^*g,
$$
where the metric $\eta^*g$ is defined by equality:
$$
\eta^*g(x)(X,Y)=g(\eta(x))\left(d\eta(X),d\eta(Y) \right),
$$
for any vector fields $X,Y$ on $M$ and any $x\in M$. It is known
\cite{Ebi1}, that for any metric $g\in\mathcal{M}$, orbit
$g\mathcal{D}$ of action $A$ is a smooth closed submanifold. The
tangent space $T_g(g\mathcal{D})$ to an orbit consists of
symmetric 2-forms of the view $h=L_Xg=\nabla_iX_j+\nabla_jX_i$
and, therefore, coincides with an image of a differential operator
$$
\alpha_g:\Gamma(TM)\longrightarrow S_2,\qquad\alpha_g (X)=L_Xg,
$$
where $L_Xg$ is Lie derivative along a vector field $X$ on $M$,
$$
T_g (g\mathcal{D})=\alpha_g\left(\Gamma (TM)\right).
$$

Adjoint operator for $\alpha_g$ is the covariant divergence,
$\alpha_g^*=2\delta_g$, where
$$
\delta_g:S_2\longrightarrow\Gamma (TM),\qquad (\delta_g a)^i
=-\nabla_ja^{ij}.
$$
The following orthogonal Berger-Ebin decomposition  \cite{Ber-Eb}
of the space $S_2=T_g\mathcal{M}$ hold:
$$
S_2=S_2^0\oplus\alpha_g(\Gamma (TM)),\eqno {(6.1)}
$$
where $S_2^0=\mbox{ker}\delta_g=\{a\in S_2;\ \delta_ga=0\}$ is
the space of divergence-free symmetric 2-forms. Accordingly, each
2-form $a\in S_2$ is represented in an unique way as:
$$
a=a^0 + L_Xg,
$$
where $\delta_ga^0=0$. The components $a^0$ and $L_Xg$ are
orthogonal and defined by an unique way.

Let $\mathcal{AM}$ is the space of associated metrics on a
symplectic manifold $(M,\omega)$. Under the action of all group
$\mathcal{D}$ of smooth diffeomorphisms, the space $\mathcal{AM}$
is not invariant. It is easy to see, that the group
$\mathcal{D}_\omega$ of smooth symplectic diffeomorphisms of a
manifold $M$, i.e. such diffeomorphisms $\eta :M\to M $, that
save the symplectic form $\eta^*\omega=\omega$, acts on the space
$\mathcal{AM}$.

The group $\mathcal {D}_{\omega}$ acts also on the space
$\mathcal{A}_{\omega}$ of positive associated almost complex
structures: if $\eta\in\mathcal{D}_{\omega}$ and
$J\in\mathcal{A}_{\omega}$, then
$\eta^{*}J=(d\eta)^{-1}\circ(J\circ\eta)\circ d\eta$. The
equivariance of a diffeomorphism $G:\ \mathcal {A}_\omega
\longrightarrow \mathcal {AM}$ (see \S 3) is easily checked.

In this paragraph we shall state orthogonal decompositions of the
spaces $S_2$ and $S_{2A}$, with respect to action of group
$\mathcal{D}_\omega$.

Recall, that the space $S_{2A}$ of anti-Hermitian symmetric
2-forms is tangent to a manifold $\mathcal{AM}$, on $M$,
$$
T_g\mathcal{AM}=S_{2A}=\{h\in S_{2};\ h(JX,JY)=-h(X,Y),\ \forall
X,Y\in\Gamma (TM)\},
$$
where $J$ is an almost complex structure, corresponding to the
metric $g$.

As all considered spaces are ILH-manifolds, then decompositions
have to be obtained for finite class of smoothness too.

Let $\mathcal{M}^s$ is the space of all metrics on $M$ of Sobolev
class $H^s$, $s\geq\frac 12\dim M+2$ and $\mathcal{AM}^s$ is the
space of associated metrics of class $H^s$. For
$g\in\mathcal{M}^s$ the tangent space $T_g\mathcal{M}^s$ to the
manifold $\mathcal{M}^s$ is the space $S_2^s$ of symmetric
2-forms of class $H^s$ on $M$, and the tangent space
$T_g\mathcal{AM}^s$ coincides with the space $S_{2A}^s$ of
anti-Hermitian symmetric 2-forms on $M$.

Let $\mathcal{D}_{\omega}$ and $\mathcal{D}_{\omega}^s$ are
groups of smooth and, of class $H^s$ symplectic diffeomorphisms of
manifold $M$ respectively. The group $\mathcal{D}_{\omega}$ is
ILH-Lie group with a Lie algebra $T_e\mathcal{D}_{\omega}$,
consisting of smooth locally Hamilton vector fields $X$ on $M$.
For each $s$ the group $\mathcal{D}_{\omega}^s$ is a continuous
group and smooth Hilbert manifold. The tangent space at unit
$T_e\mathcal{D}_{\omega}^s$ consists of locally Hamilton vector
fields $X$ on $M$ of class $H^s$.

At first we consider an ILH-Lie group $\mathcal{G}$ of exact
symplectic diffeomorphisms. Its Lie algebra is the algebra
$\mathcal{H}$ of Hamilton vector fields $X_F$ on $M$ \cite{Omo5}.
The arbitrary Hamilton vector field $X_F$ can be presented as
$X_F=J\grad F$, where $F$ is some function on $M$, called a
Hamiltonian of the field $X_F$ and $J$ is almost complex
structure corresponding to the metric $g$. Therefore orthogonal
complement $\mathcal{H}^{\perp}$ to $\mathcal{H}$ in the space
$\Gamma (TM)$ of all vector fields consists of vector fields $V$
on $M$, satisfying to the condition $\div JV=0$.

Fix the Riemannian metric $g\in\mathcal{AM}^s$ and consider its
orbit $g\mathcal{G}$. The tangent space $T_g(g\mathcal{G})$
consists of 2-forms of the view $h=L_{X_F}g$. In this connexion
we consider a differential operator, acting on functions:
$$
D_g:\ H^{s+2}(M,{\bf R})\longrightarrow S_2^s,\qquad
D_g(F)=L_{X_F}g.
$$

Let $\Im\ D_g$ is image of an operator $D$ and $\Ker\ D_g^*$ is
kernel of a adjoint operator $D_g^*:\ S_2^s\longrightarrow
H^{s-2}(M,{\bf R})$.

Recall, that the Riemannian metric $g$ defines an inner product
(1.4) in the space $S_2^s$, with respect to which the adjoint
operator is defined.

\vspace{3mm}

{\bf Theorem 6.1.} {\it The space $S_2^s$ is decomposed in a
direct sum of orthogonal subspaces
$$
S_2^s=\Ker\ D_g^*\oplus\Im\ D_g.\eqno {(6.2)}
$$
Accordingly with it each symmetric 2-form $h$ of class $H^s$ is
represented as
$$
h= h^*+L_Xg,\eqno {(6.3)}
$$
where $X=X_F$ is the Hamilton vector field of class $H^{s+1}$, and
$h^*$ satisfies to a condition $\div (J\delta h^*)=0$.}

{\bf Proof.} For function $F\in H^{s+1}(M,{\bf R})$ we have
$J(\grad F)=X_F$. Therefore operator $D_g$ is a composition of
three operators, $D_g=\alpha_g\circ J\circ\grad$ (recall, that
$\alpha_g(X)=L_Xg$). The operators $\grad$ and $\alpha_g$ have
injective symbols \cite{Ber-Eb}, consequently, $D_g$ has an
injective symbol too. From \cite{Ber-Eb} (theorem 6.1), we obtain
decomposition (6.2). Find a adjoint operator,
$D_g^*=(\alpha_g\circ J\circ\grad)^*=\grad^*\circ
J^*\circ(\alpha_g)^*= -\div\ \circ (-J)\circ(2\delta_g)=2\div\
\circ J\circ\delta_g$. Therefore $h^*\in\Ker\ D_g^*$ satisfies to
the condition $\div\ (J\delta_gh^*)=0$.

\vspace{3mm}

{\bf Remark.} The space of 2-forms $h^*$ of class $H^s$,
satisfying to the condition,
$$
\div (J\delta_g h^*)=0,
$$
we shall designate by the symbol $S_2^{*s}$. Sometimes we shall
write the decomposition (6.2) as
$$
S_2^s=S_2^{*s}\oplus\alpha_g(\mathcal{H}^{s+1}),\eqno{(6.2')}
$$
where $\mathcal{H}^{s+1}$ is the space of Hamilton vector fields
on $M$ of class $H^{s+1}$.

\vspace{3mm}

Consider the space $\mathcal{AM}^s$ of associated metrics. As the
group $\mathcal{G}^{s+1}$ acts on $\mathcal{AM}^s$, the tangent
space $T_g\mathcal{AM}^s=S_{2A}^s$ contains the tangent space to
an orbit $T_g(g\mathcal{G}^{s+1})$. As
$T_g(g\mathcal{G}^{s+1})=\alpha_g (T_e\mathcal {G}^{s+1})
=\alpha_g(\mathcal{H}^{s+1})$,
$$
\alpha_g(\mathcal{H}^{s+1})\subset S_{2A}^s.
$$

Therefore we obtain decomposition of the space $S_{2A}^s$ of
symmetric anti-Hermitian 2-forms.

\vspace{3mm}

{\bf Corollary 1.} {\it The space $S_{2A}^s$ is decomposed in an
orthogonal direct sum
$$
S_{2A}^s=S_{2A}^{*s}\oplus\alpha_g(\mathcal{H}^{s+1}),\quad
h=h^*+L_Xg,
$$
where $X=X_F$ is the Hamilton vector field of class $H^{s+1}$, and
$h^*$ is anti-Hermitian and satisfies to the condition $\div
(J\delta h^*)=0$.}

\vspace{3mm}

The space $S_2^{*s}$ can be decomposed on orthogonal subspaces
according to Berger-Ebin decomposition
$$
S_2^{s}=S_2^{0s}\oplus\alpha_g (\Gamma^{s+1}(TM)).
$$
As in our case $\Im\ D_g=\alpha_g(\mathcal{H})\subset\Gamma (TM)$,
the space $S_2^{0s}$ is in $\Ker\ D_g^*$. Therefore the space
$\Ker\ D_g^*$ can be decomposed on orthogonal subspaces according
to Berger-Ebin decomposition:
$$
\Ker\ D_g^*=S_2^{0s}\oplus B^s,
$$
where
$$
B^s=\Ker\ D_g^*\cap\alpha_g(\Gamma^{s+1}(TM))=\{h\in S_2^s;\
h=L_Yg,\ \div J\delta_g(h)=0\}.
$$

\vspace{3mm}

{\bf Theorem 6.2.} {\it The space $S_2^s$ is decomposed in a
direct sum of orthogonal subspaces
$$
S_2^{s}=S_2^{0s}\oplus B^s\oplus\alpha_g(\mathcal{H}^{s+1}).\eqno
{(6.4)}
$$
According to this, each symmetric 2-form $h$ of class $H^s$ is
represented in an unique way as
$$
h=h^0+L_Yg+L_Xg,
$$
where $X=X_F$ is Hamilton vector field of class $H^{s+1}$, $h^0$
has the property $\delta_g(h^0)=0$, and vector field $Y$ of class
$H^s$ is such that $\div J\delta_gL_Yg = 0$.}

\vspace{3mm}

{\bf Remark.} For the vector field $\delta_gL_Yg$ there is a
following expression from \cite{Fis-Mar1}:
$$
\delta_gL_Yg=\Delta Y-\grad(\div\ Y)-2Ric(g)Y,
$$
where $Ric(g)$ is Ricci tensor and $(\Delta Y)_\sharp=(dd^*+
d^*d)(Y_\sharp)$, $Y_\sharp$ is 1-form obtained from a vector
field $Y$ by omitting of an index.

\vspace{3mm}

State variant of the last decomposition, when instead of the space
$\mathcal {H}$ of Hamilton vector fields on $M$ the space
$T_e\mathcal{D}_\omega$ of locally Hamilton vector fields is
taken. It is known, that $T_e\mathcal {D}^{s + 1}_\omega
/\mathcal{H}^{s+1}$ is isomorphic to the first cohomology group
$H^1(M,{\bf R})$ of manifold $M$. One can think, that $H^1(M,{\bf
R})$ is represented by vector fields of class $C^{\infty}$ and
$T_e\mathcal{D}^{s+1}_\omega=\mathcal{H}^{s+1}\oplus H^1(M,{\bf
R})$ is direct sum.

The operator $\alpha_g:\ \Gamma^{s+1}(TM)\longrightarrow S_2^s$,
$\alpha_g (X) = L_Xg$, maps the finite-dimensional space
$H^1(M,{\bf R})$ on the closed finite-dimensional subspace
$\alpha_g(H^1(M,{\bf R}))\subset S_2^s$. We form a direct sum
$$
\Im\ D_g\oplus\alpha_g(H^1(M,{\bf
R}))=\alpha_g(T_e\mathcal{D}^{s+1}_\omega).
$$
Consider the orthogonal complement
$$
S_2^{\omega s}=\alpha_g(T_e\mathcal{D}^{s+1}_\omega)^{\perp}
$$
in the space $S_2^{s}$, then
$$
S_2^{s}=S_2^{\omega s}\oplus\alpha_g(T_e\mathcal{D}^{s+1}_\omega).
$$
Let $(J\delta_g h)_\sharp$ is 1-form obtained from a vector field
$J\delta_g h$ by omitting of an index via the metric tensor $g$.
Then the belonging of $h$ to the space $\Ker\ D_g$ is expressed
by condition of co-closed: $d^*(J\delta_gh)_\sharp=0$, and the
space $S_2^{\omega s}$ consists of such $h$, that 1-form
$(J\delta_gh)_\sharp$ is $d^*$-exact, where $d^*$ is
codifferential.

\vspace{3mm}

{\bf Theorem 6.3.} {\it The space $S_2^s$ is decomposed in a
direct sum of orthogonal subspaces
$$
S_2^{s} = S_2^{\omega
s}\oplus\alpha_g(T_e\mathcal{D}^{s+1}_\omega). \eqno {(6.5)}
$$
Each form $h\in S_2^{s}$ is represented in an unique way as
$h=h^\omega+L_Xg$, where $X$ is locally Hamilton vector field of
class $H^s$, and $h^\omega$ is such that the 1-form
$(J\delta_gh)_\sharp$ is $d^*$-exact.}

\vspace{3mm}

Consider the space $\mathcal{AM}^s$ of associated metrics. As the
group $\mathcal{D}_\omega^{s+1}$ acts on $\mathcal{AM}^s$, the
tangent space $T_g\mathcal{AM}^s=S_{2A}^s$ contains the space
$\alpha_g(T_e\mathcal{D}_\omega^{s+1})$. It is tangent to an
orbit $g\mathcal{D}_\omega^{s+1})$ of action of group
$g\mathcal{D}_\omega^{s+1})$ on $\mathcal{AM}^s$.

From the theorem 6.3 is obtained

{\bf Corollary 2.} {\it The space $S_{2A}^{s}$ is decomposed in
an orthogonal direct sum
$$
S_{2A}^{s}=S_{2A}^{\omega
s}\oplus\alpha_g(T_e\mathcal{D}_\omega^{s+1}), \eqno {(6.6)}
$$
where $S_{2A}^{\omega s}=S_{2}^{\omega s}\cap S_{2A}^{s}$
consists of anti-Hermitian symmetric 2-forms $h$ of class $H^s$
such, that the form $(J\delta_g h)_\sharp$ is $d^*$-exact.}

\vspace{3mm}

As $S_2^{\omega s}$ contains the space $S_2^{0s}$ of
divergence-free forms, $S_2^{\omega s}$ is decomposed further in
correspondence with Berger-Ebin decomposition:
$$
S_2^{\omega s}= S_2^{0s}\oplus\widetilde{B}^s.
$$
Thus,
$$
\widetilde {B}^s=S_2^{\omega
s}\cap\alpha_g(\Gamma^{s+1}(TM))\qquad\mbox{and} \qquad
B^s\cong\widetilde{B}^s\oplus\alpha_g(H^1(M,{\bf R})).
$$

\vspace{3mm}

{\bf Theorem 6.4.} {\it The space $S_2^s$ is decomposed in the
orthogonal direct sum
$$
S_2^{s}=S_2^{0s}\oplus\widetilde{B}^s\oplus\alpha_g
(T_e\mathcal{D}_\omega^{s+1}).\eqno {(6.7)}
$$
According to this each tensor field $h\in S_2^s$ is represented
in an unique way as
$$
h=h^0+L_Yg+L_Xg,\eqno {(6.8)}
$$
where $X$ is locally Hamilton vector field of class $H^s$, $h^0$
has the property $\delta_g(h^0)=0$, and the vector field $Y$ of
class $H^s$ is such, that the 1-form $\gamma (Y)=\left(J(\Delta
Y-\grad\ \div Y-2Ric(g)Y\right)_\sharp$ is $d^*$-exact.}

\vspace{3mm}

{\bf Proof.} It is only required to prove the property of the
field $Y$. Let $Y$ is such vector field, that $L_Yg\in B^s$. It
follows from an orthogonality of the spaces $B^s$ and
$\alpha_g(T_e\mathcal{D}_\omega^{s+1})$, that for any $X\in
T_e\mathcal{D}_\omega^{s+1}$ we have $(L_Yg,L_X g)_g = 0$ (here
$(..)_g$ is inner product in the space $S_2$). On the other hand,
$(L_Yg,L_Xg)=(L_Yg,\alpha_g(X))=(\alpha_g^*(L_Yg),X)=2(\delta_g(L_Yg),X)$.
Therefore vector field $\delta_g(L_Yg)$ is orthogonal to
$T_e\mathcal{D}_\omega^{s+1}$. It is equivalent to
$d^*$-exactness of the 1-form $(J\delta_gL_Yg)_\sharp$. Now we
use the formula from \cite{Fis-Mar1}:
$$
\delta_gL_Yg=\Delta Y-\grad(\div\ Y)-2Ric(g)Y.
$$

\vspace{3mm}

{\bf Remark.} In decomposition (6.8) vector fields $X$ and $Y$
are defined on $h$ up to a Killing vector field.

\vspace{3mm}

Consider decomposition (6.8) of anti-Hermitian 2-forms.

\vspace{3mm}

{\bf Lemma 6.1.} {\it If $h\in S_{2A}^{s}$, then in decomposition
(6.8) the second component $L_Yg$ is unique defined on $h^0$.}

\vspace{3mm}

{\bf Proof.} We can think, that $h=h^0+L_Yg\in S_{2A}^{\omega s}$.
Assume, that there are two forms $h_1,h_2\in S_{2A}^{\omega s}$,
such, that $h_1\ne h_2$ and $h_1=h^0+L_{Y_1}g$,
$h_2=h^0+L_{Y_2}g$. Then $L_{Y_1-Y_2}g=h_1-h_2\in S_{2A}^{\omega
s}\subset T_g\mathcal{AM}^s=S_{2A}^s$. On the other hand,
$L_{Y_1-Y_2}g\in\alpha_g(T_e\mathcal{D}_\omega^{s+1})$. It
follows from the M.Gromov's maximality theorem.
Indeed, let $\mathcal{D}(\mathcal{AM}^s)$ is group of all $H^s$ -
diffeomorphisms of a manifold $M$, saving the space
$\mathcal{AM}^s$. This is the closed subgroup in the group
$\mathcal{D}_\mu^{s+1}$ of diffeomorphisms saving the volume form
$\mu=\frac{1}{n!}\omega^n$ and
$\mathcal{D}_\omega^{s+1}\subset\mathcal{D}(\mathcal{AM}^s)
\subset\mathcal{D}_\mu^{s+1}$. Under the Gromov's theorem the
group $\mathcal{D}(\mathcal{AM}^s)$ coincides either with
$\mathcal{D}_\omega^{s+1}$, or with $\mathcal{D}_\mu^{s+1}$. If
we shall assume, that
$\mathcal{D}(\mathcal{AM}^s)=\mathcal{D}_\mu^{s+1}$, for any
$Y\in T_e\mathcal{D}_\mu^{s+1}$ we obtain, that $L_Yg\in
T_e\mathcal{AM}^s$, i.e. for any divergence-free field $Y$,
2-form $L_Yg$ is anti-Hermitian. It is obviously impossible.
Therefore $\mathcal{D}(\mathcal{AM}^s)=\mathcal{D}_\omega^{s+1}$.
Therefore, it follows from $L_{Y_1-Y_2}g\in T_g\mathcal{AM}^s$
that $L_{Y_1-Y_2}g\in\alpha_g(T_e\mathcal{D}_\omega^{s+1})$. The
simultaneous inclusion $L_{Y_1-Y_2}g\in S_{2A}^{\omega s}$ is
possible only at $L_{Y_1-Y_2}=0 $, i.e. at $L_{Y_1}=L_{Y_2}$.

\vspace{3mm}

In decomposition $S_{2A}^{s}=S_{2A}^{\omega
s}\oplus\alpha_g(T_e\mathcal{D}_\omega^{s+1})$ of the space
$S_{2A}^{s}$ of anti-Hermitian forms the first component
$S_{2A}^{\omega s}$ can be decomposed further in correspondence
with Berger-Ebin decomposition:
$$
S_{2A}^{\omega s}=AS_2^{0s}\oplus C^s,\quad
h=h^0+L_Yg.\eqno{(6.9)}
$$
Thus $C^s\subset\widetilde{B}^s$ and the second component $L_Yg$
is unique defined by $h_0$.

\vspace{3mm}

{\bf Corollary 3.} {\it In the decomposition (6.9) the space
$S_{2A}^{\omega s}$ is isomorphic to the closed subspace
$AS_{2}^{0 s}$ in the space $S_{2}^{0s}=\Ker\ \delta_g$. The
isomorphism is stated by projection $p(h)=p(h^0+L_Yg)=h^0$.}

\vspace{3mm}

Describe the space $AS_{2}^{0 s}$, i.e. describe those
divergence-free elements $h^0\in S_{2}^{0s}$, which are
components of decomposition of anti-Hermitian forms $h$.

Let $S^s_{2H}$ is the space of Hermitian symmetric 2-forms of
class $H^s$ on $M$. Natural orthogonal decomposition takes place
$$
S^s_2=S^s_{2A}\oplus S^s_{2H},
$$
$$
h(X,Y)=\frac 12\left(h(X,Y)-h(JX,JY)\right)+\frac
12\left((h(X,Y)+h(JX,JY) \right).
$$

\vspace{3mm}

{\bf Theorem 6.5.} {\it The element $h^0\in S_{2}^{0s}$ is
component of decomposition (6.8) of element $h\in S_{2A}^{s}$ iff
$h^0$ is orthogonal in $S_{2}^{0s}$ to the subspace
$S_{2H}^{0s}=S_{2}^{0s}\cap S_{2H}^{s}$ of the Hermitian
divergence-free forms.}

\vspace{3mm}

{\bf Proof.} Let $h=h^0+L_Yg$ is anti-Hermitian form. As $h\perp
S_{2H}^s$ and $L_Yg\perp S_2^{0s}$, $h^0=h-L_Y g$ is orthogonal
to intersection $S_{2H}^{0s}=S_{2H}^s\cap S_2^{0s}$. Conversely,
suppose, that $h^0\perp S_{2H}^{0s}=S_{2H}^s\cap S_2^{0s}$. Then
$h^0\in (S_{2H}^{0 s})^{\perp}=(S_{2H}^s\cap S_2^{0 s})^{\perp}
=(S_{2H}^s)^{\perp}\cup (S_2^{0 s})^{\perp}=S_{2A}^{s}\cup\alpha_g
(\Gamma^{s+1}(TM))$. Therefore $h^0$ is a linear combination of
an element $L_Vg$, $V\in\Gamma^{s+1}(TM)$ and anti-Hermitian form
$h$, $h^0=h-L_Vg$. We have obtained, that $h=h^0+L_Vg$ is the
anti-Hermitian form, at which $h^0$ is a divergence-free
component.

\vspace{3mm}

{\bf Corollary 4.}
$$
AS_{2}^{0s}=(S_{2H}^{0s})^{\perp},
$$
{\it where the orthogonal complement is taken in the space
$S_{2}^{0 s}$ of divergence-free forms.}

\vspace{3mm}

{\bf Corollary 5.} {\it Nonzero Hermitian form $h^0$ can not be a
divergence-free component of decomposition (6.7) of element $h\in
S_{2A}^{s}=T_g\mathcal{AM}^s$.}

\vspace{3mm}

Decompose component $h^0$ on Hermitian and anti-Hermitian parts,
$$
h^0=(h^0)_H+(h^0)_A.
$$
We have got, that the component $(h^0)_A$ of the nonzero form
$h^0$  is not equal to zero too. At the same Hermitian part
$(h^0)_H $, the second component $(h^0)_A$ is defined up to the
anti-Hermitian divergence-free form $h^0_A\in
S_{2A}^{0s}=S_{2A}^s\cap S_2^{0s}$. Therefore if $(h^0)_H\ne 0$,
$\delta_g ((h^0)_H)\ne 0 $, and so $\delta_g((h^0)_A)\ne 0$.

Consider Hermitian part $(h^0)_H$ of a divergence-free component
$h^0$ of the forms $h\in S_{2A}^{s}$.

Define mapping
$$
\mathbf{J}:\ \Gamma(T_2^0 M)\longrightarrow\Gamma(T_2^0M),
$$
which takes the 2-form $a(X,Y)$ on $M$ in the 2-form
$\mathbf{J}(a)(X,Y)=a(X,JY)$. The mapping $\mathbf{J}$ commutes
with an operator of taking of Hermitian part $a_H$ of the form
$a$:
$$
\mathbf{J}(a_H)=(\mathbf{J}a)_H.
$$

Note an obvious isomorphism
$$
\mathbf{J}:\ S_{2H}^{s}\oplus S_{2A}^{s}\longrightarrow
H^s(\Lambda^2_H(M)) \oplus S_{2A}^{s},
$$
where $\Lambda^2_H(M)$ is bundle of Hermitian skew-symmetric
2-forms on $M$. If $\varphi\in\Lambda^2_H(M)$, it is easy to see,
that $\mathbf{J}\varphi(X,Y)=\varphi(X,JY)$ is Hermitian
symmetric 2-form.

\vspace{3mm}

{\bf Theorem 6.6.} {\it The space $S_{2A}^{s}$ is isomorphic to a
direct sum of orthogonal subspaces
$$
S_{2A}^{s}\cong
AS_{2}^{0s}\oplus\alpha_g(T_e\mathcal{D}_\omega^{s+1}).
\eqno{(6.10)}
$$
Each tensor field $h\in S_{2A}^{s}$ is unique represented as
$$
h=h^0+L_Yg+L_Xg,\eqno{(6.11)}
$$
where $X$ is locally Hamilton vector field of class $H^s$, the
component $L_Yg$ is unique defined on $h_0$, a vector field $Y$
and 2-form $h^0$ have properties:

1) $\delta_g h^0=0$,

2) $h^0_H=-\mathbf{J}((L_Y\omega)_H)$,

3) $\gamma(Y)=\left(J(\Delta Y-\grad\ \div\
Y-2Ric(g)Y)\right)_\sharp$ is $d^*$-exact 1-form.}

\vspace{3mm}

{\bf Proof.} From the lemma 1 is obtained, that $L_Yg$ is unique
defined by $h^0$. The properties 1 and 3 follow from the theorem
6.4. Let $Y$ is vector field from decomposition (6.11) and
$\overline{h}=L_Yg$. Applying Lie derivative $L_Y$ to the left
and right parts of equality
$$
g(JV,W)-g(V,JW)=\omega(V,W)+\omega(JV,JW)
$$
and taking into account Leibniz rule, we obtain
$$
\overline{h}(JV,W)-\overline{h}(V,JW)=L_Y\omega(V,W)+L_Y\omega(JV,JW)=
2(L_Y\omega)_H(V,W).\eqno {(6.12)}
$$
Substituting instead of $W$ a vector field $JW$, we have
$$
2(\overline{h})_H=\overline{h}(JV,JW)+\overline{h}(V,W)=L_Y\omega(V,JW)
-L_Y\omega(JV,W)=2\mbox{\bf J}(L_Y\omega)_H(V,W).
$$
Thus, $(\overline{h})_H=\mbox{\bf J}(L_Y\omega)_H$. As the form
$h=h^0+L_Yg+L_Xg=h^0+\overline{h}+L_Xg$ belongs to the space
$S_{2A}^{s}=T_g\mathcal{AM}^s$ and $L_Xg\in T_g\mathcal{AM}^s$,
$(h)_H=0$ and $(L_Xg)_H=0$. Therefore $(h^0+\overline{h})_H=0$, so
$(h^0)_H=-(\overline{h})_H$. Finally we have $(h^0)_H=-\mbox{\bf
J}(L_Y\omega)_H$.

\vspace{3mm}

{\bf Remark.} We have obtained, that for any vector field $Y$ on
$M$ the equalities take place:
$$
L_Y\omega=-\mathbf{J}(L_Yg)-gL_YJ,\qquad\mathbf{J}(L_Y\omega)_H=(L_Yg)_H.
$$

\vspace{3mm}

{\bf Proposition 6.7.} {\it The space $AS_{2}^{0s}$ consists of
symmetric 2-forms $h^0\in S_{2}^{0s}$, having a property,
$$
(h^0)_H=\mathbf{J}(d\beta)_H,
$$
where $\beta$ is 1-form on $M$.}

\vspace{3mm}

{\bf Proof.} As $L_Y\omega$ is exact form, in one direction the
statement follows from the theorem 6.6. Conversely, let $h^0$ is
such that $\delta_gh^0=0$ and $(h^0)_H=\mbox{\bf J}(d\beta)_H$
for some exact 2-form $d\beta$. Find a vector field $Y$,
possessing properties 2,3 of theorems 6.6. For this purpose we
shall find a vector field $Z$ from equality
$\iota_Z\omega=-\beta$. From here $L_Z\omega=-d\beta$.
Decomposition (6.11) implies, that $L_Zg$ is orthogonal
$S_{2}^{0s}$, therefore under the theorem 6.4 $L_Z g$ is
decomposed as $L_Zg=L_Yg+L_Xg$, where $Y$ has the necessary
property 3. As $X\in T_e\mathcal{D}_{\omega}^{s+1}$, $L_X\omega=0
$. As a result $L_Y\omega=L_Z\omega=-d\beta$ and
$(h^0)_H=\mbox{\bf J}(L_Y\omega)_H$.

\vspace{3mm}

{\bf Remark.} One can point decomposition of the space $S_2$,
dual to the decomposition (6.2) in some sence. Namely, any
symmetric form $h\in S_2$ is decomposed in a sum of orthogonal
addends
$$
h=h^1+L_Xg,\eqno {(6.13)}
$$
where $X$ is exact divergence-free vector field on $M$, and $h^1$
has the property, that a vector field $J\delta_g h^1$ is locally
Hamilton. If one requires of a field $X$ to be only
divergence-free in (6.13), then vector field $J\delta_gh^1$
should be Hamilton.

Recall, that exact divergence-free is such vector field $X$, for
which the 1-form $X_\sharp=g^{-1}X$ is $d^*$ - exact.

The decomposition (6.13) turns out the same as (6.2), via using
of a differential operator $R:\
\Gamma(\Lambda^2TM)\longrightarrow S_2$, which is defined by
equality $R(\varphi)=\alpha_g(d^*\varphi)^{\sharp}$, where the
symbol $\sharp$ designates operation of a raising of an index.

\newpage

\vspace {10mm} \centerline {\bf \S 7. A curvature of a quotient
space $\mathcal{AM}/\mathcal{D}_{\omega}$.}

\vspace {7mm}

Let $\mathcal{AM}$ is the space of all smooth associated metrics
on a symplectic manifold $M,\omega $. The group
$\mathcal{D}_\omega$ of symplectic diffeomorphisms of a manifold
$M$ acts natural on this space $\mathcal{AM}$. Let
$\mathcal{D}_{\omega 0}$ is connected component of the unit of
group $\mathcal{D}_\omega$. Lie algebra of the group
$\mathcal{D}_\omega$ consists of locally Hamilton vector fields
on $M$.

Suppose that the symplectic structure $\omega$ allows integrable
associated almost complex structures. In this case commutator
$\mathcal{G} =[\mathcal{D}_{\omega 0},\mathcal{D}_{\omega 0}]$ is
\cite{Rat-Shm}, \cite{Ban2} the connected closed ILH-Lie subgroup
of the group $\mathcal{D}_{\omega 0}$ with a Lie algebra
$\mathcal{H}$, consisting of Hamilton vector fields on $M$. At
that $\mathcal{D}_{\omega 0}/\mathcal{G}\cong H^1(M,{\bf
R})/\Gamma$, where $\Gamma$ is some discrete subgroup of the
first group of cohomologies $H^1(M,{\bf R})$.

Consider the matter of a curvature of the space
$\mathcal{AM}/\mathcal{G}$ of classes of equivalent associated
metrics. The quotient space $\mathcal{AM}/\mathcal{G}$ is not a
manifold, in general, as it has singularities corresponding to
such of metrics $g$, that there isometry group $I(g)$ has
nontrivial intersection with the group $\mathcal{G}$. Calculate a
curvature $\mathcal{AM}/\mathcal{G}$ in the regular points
$[g]=g\mathcal{G}$, i.e. in such that $I(g)\cap\mathcal{G}={e}$.
The set of such classes $[g]$ is open in
$\mathcal{AM}/\mathcal{G}$. Indeed, it contains metrics with
discrete group $I(g)$. The set of metrics with discrete (and even
with trivial) isometry group is open and is everywhere dense in
the space of metrics $\mathcal{AM}$ \cite {Ebi1}.

The weak Riemannian structure on $\mathcal{AM}$ is invariant
\cite {Ber-Eb} with respect to action of the group
$\mathcal{D}_{\omega}$. Therefore on a regular part of the space
$\mathcal{AM}/\mathcal{G}$ a weak Riemannian structure such, that
the projection $p:\
\mathcal{AM}\longrightarrow\mathcal{AM}/\mathcal{G}$ is a
Riemannian submersion, is natural defined \cite {ONe}. The
vertical subbundle $V(\mathcal{AM})$ consists of subspaces, which
are tangents to orbits $g\mathcal{G}$. Recall, that
$V_g=T_g(g\mathcal{G})$ consists of 2-forms of the view $a=L_Xg$,
where $L_X$ is the Lie derivate along a Hamilton vector field $X$.

It is shown in \S 6, that the horizontal subspace
$H_g=V_g^{\perp}$ consists of the anti-Hermitian symmetric
2-forms $a$ such, that $\div J\delta_ga=0$, where $\delta_g$ is
covariant divergence: $(\delta_g a)^i =-\nabla_k a^{kj}$. Recall
that for an operator $\delta_g: S_2\longrightarrow\Gamma(TM)$ the
adjoint operator $\delta_g^{*}:\Gamma(TM)\longrightarrow S_2$ is
expressed via a Lie derivative: $\delta_g^{*}=\frac
12\alpha_g(X)=\frac 12L_Xg$.

Let $\overline{a},\overline{b}\in
T_{[g]}(\mathcal{AM}/\mathcal{G})$. As
$p:\mathcal{AM}\longrightarrow\mathcal{AM}/\mathcal{G}$ is
Riemannian submersion, the sectional curvature
$\overline{K}(\overline{a},\overline{b})$ of the space
$\mathcal{AM}/\mathcal{G}$ is under the formula \cite{ONe},
\cite{Bes2}:
$$
\overline{K}(\overline{a},\overline{b})= K(a,b)+\frac
34\frac{\left([a,b]^V,[a,b]^V\right)_g}{\|a\wedge b\|^2}, \eqno
{(7.1)}
$$
where $a,b\in T_g\mathcal{AM}$ are horizontal lifts of vectors
$\overline{a},\overline{b}$, $K(a,b)$ is sectional curvature of
the space $\mathcal{AM}$, $[a,b]^V$ is vertical part of a Lie
commutator of horizontal vector fields on $\mathcal{AM}$, which
continue $a, b\in T_g \mathcal {AM}$.

For an evaluation $[a,b]^V$ we will need two differential
operators:

1) Elliptic operator of the 4-th order:
$$
\mathrm{E}_g:C^\infty(M,{\bf R})\longrightarrow C^\infty (M,{\bf
R}),\qquad \mathrm{E}_g(f)=\div J\delta_g\ \alpha_gJ\grad f,\eqno
{(7.2)}
$$
It is obvious, that its kernel consists of a constant functions;

2) Differential operator of the 1-st order \cite {Ber-Eb}:
$$
\framebox(8,8){}:\ S_2\longrightarrow\Gamma(S_2 M\otimes TM),\quad
(\framebox(8,8){}\ a)_{ij}^k=a_{ij}^k=\frac 12\left(\nabla_i
a_j^k+ \nabla_ja_i^k-\nabla^ka_{ij}\right),\quad a\in S_2.
$$
The geometrical sense of an operator $\framebox(8,8){}$  is, that
it is a differential of mapping $\Gamma:\
\mathcal{M}\longrightarrow RConn$, which takes a Riemannian
metric $g\in\mathcal{M}$ to  Riemannian connection without a
torsion $\Gamma(g)$. If $g_t$ is curve on the space of metrics
with tangent vector $a=\left.\frac{d}{dt}\right|_{t=0}g_t$ and
$\Gamma_{ij}^k(g_t)$ are Christoffel symbols of a Riemannian
connection of the metric $g_t$,
$$
(\framebox(8,8){}\
a)_{ij}^k=\left.\frac{d}{dt}\right|_{t=0}\Gamma_{ij}^k(g_t).
$$

Define a contraction $\framebox(8,8){}\ a$ and $\delta_g a$ with
the symmetric 2-form $b$:
$$
(\framebox(8,8){}\ a,b)^k=a_{ij}^kb^{ij},\quad(\delta_g a,b)^k=
b_i^k(\delta_ga)^i.\eqno {(7.3)}
$$

Introduce the notation:
$$
\{a, b\}=(\delta_gb,a)+(\framebox(8,8){}\ b,a)-(\delta_g
a,b)-(\framebox(8,8){} \ a, b).\eqno {(7.4)}
$$

{\bf Theorem 7.1.} {\it If $a,b$ are horizontal vector fields on
the space $\mathcal {AM}$, then
$$
[a,b]^V=-\alpha_gJ\ \grad\ \mathrm{E}_g^{-1}\div
J\left(\{a,b\}\right). \eqno {(7.5)}
$$}

{\bf Theorem 7.2.} {\it Let $[g]$ is regular point of the
quotient space $\mathcal{AM}/\mathcal{G}$ and $g\in\mathcal{AM}$
is any metric from the class $[g]$. The sectional curvature
$K(a,b)$ of the space $\mathcal {AM}/\mathcal{G}$ in a plane
section, defined by pair of vectors $\overline{a},\overline{b}\in
T_{[g]}(\mathcal{AM} /\mathcal{G})$, is expressed by the formula
$$
\overline{K}(\overline{a},\overline{b})=K(a,b)+\frac{3}{4\|a\wedge
b\|^2} \int_M\left(\div
J\{a,b\}\right)\left(\mathrm{E}_g^{-1}(\div
J\{a,b\})\right)d\mu_g, \eqno {(7.6)}
$$
where $a,b\in T_g\mathcal{AM}$ are horizontal lifts of vectors
$\overline{a},\overline{b}$ and $K(a,b)$ is sectional curvature
of the space $\mathcal{AM}$.}

\vspace{3mm}

{\bf Proof of the theorem 7.1.} The vertical subspace $V_g$
consists of 2-forms of a view $h=L_Xg$, where $X$ is Hamiltonian
vector field, $X=J\grad f$. Thus, $h=\alpha_gJ\grad f$. For such
forms $h$ the equality $\div J\delta_gh=\div
J\delta_g\alpha_gJ\grad f=\mathrm{E}_g(f)$, where $\mathrm {E}_g$
is operator defined earlier by the formula (7.2), is held.
Therefore orthogonal projection of the tangent space
$T_g\mathcal{AM}$ on vertical $V_g$ can be given as follows:
$$
h^V=\alpha_gJ\ \grad\ \mathrm{E}_g^{-1}\div
J\delta_g(h).\eqno{(7.7)}
$$
Find a projection of Lie commutator $[a,b]$ of horizontal vector
fields $a,b$ on $\mathcal{AM}$. As
$\mathcal{AM}\subset\mathcal{M}\subset S_2$,
$[a,b]_g=db(a)-da(b)$, where $db(a)$ is usual derivative of a
field $b$ in direction $a$ in the vector space $S_2$.

At first we shall obtain $\div J \delta_g(db(a))$ in a considered
point $g$. Take a curve $g_t$ on $\mathcal{AM}$, going out in
direction $a$, $\left.\frac {d}{dt}\right|_{t=0}g_t=a$. Vector
field $b$ is horizontal, $\div_tJ_t\delta_{g_t}(b(g_t))=0$,
therefore
$$
\left.\frac{d}{dt}\right|_{t=0}\div_tJ_t\delta_{g_t}(b(g_t))=0.
\eqno {(7.8)}
$$
Calculate derivatives of operators $\div_t$, $J_t$,
$\delta_{g_t}$,
$$
\left.\frac{d}{dt}\right|_{t=0}\div_tX=\left.\frac{d}{dt}\right|_{t=0}
(\nabla_k (t)X^k)=\left.\frac{d}{dt}\right|_{t=0}\left(
\frac{\partial X^k}{\partial
x^k}+\Gamma^k_{ki}X^i\right)=a_{ki}^kX^i=0,
$$
where $a_{ij}^k=(\framebox(8,8){}\ a)_{ij}^k=\frac
12\left(\nabla_ia_j^k+ \nabla_ja_i^k-\nabla^ka_{ij}\right)$.
Taking into account, that $\tr\ g^{-1}a=a_k^k=0$, we obtain, that
$a_{ki}^k=\frac
12\left(\nabla_ia_k^k+\nabla_ka_i^k-\nabla^ka_{ik}\right)=0$.

Let $A=g^{-1}a$. Then, it is known (\S 4), that
$$
\left.\frac{d}{dt}\right|_{t=0}J_t=JA.
$$
Further, taking into account that
$$
\left.\frac{d}{dt}\right|_{t=0}g_t^{ki}=-a^{ki}\ \mbox{and}\
\left.\frac{d}{dt}
\right|_{t=0}(\nabla_kb_{is})=-a_{ki}^qb_{qs}-a_{ks}^qb_{iq},
$$
we obtain:
$$
\left.\frac{d}{dt}\right|_{t=0}\delta_{g_t}(b)=\left.\frac{d}{dt}\right|_{t=0}
(-g^{kl}g^{si}\nabla_kb_{ls})=
$$
$$
=a^{kl}g^{si}\nabla_kb_{ls}+g^{kl}a^{si}\nabla_kb_{ls}+
g^{kl}g^{si}a_{kl}^qb_{qs}+g^{kl}g^{si}a_{ks}^qb_{lq} =
$$
$$
=a^{kl}\nabla_kb_l^{i}-(A(\delta_gb))^i-(B(\delta_ga))^i+g^{si}a_{ks}^qb_q^k.
$$
Using these expressions of derivatives, equality (7.8) can be
written as
$$
\div J\left((A(\delta_gb))^i+a^{kl}\nabla_kb_l^{i}(A(\delta_gb))^i
(B(\delta_g a))^i+g^{si}a_{ks}^qb_q^k\right)+\div
J\delta_g(db(a))=0.
$$
From here we obtain, that
$$
\div
J\delta_g(db(a))=-\nabla_p\left(J_i^p(a^{kl}\nabla_kb_l^i-b_k^i
(\delta_g a)^k+\frac 12\nabla^i(a_{kl})b^{kl})\right).
$$
$\div J \delta_g (da (b))$ is similarly calculated. Using operator
$\framebox(8,8){}\ $, we see, that
$$
\div J\delta_g[a,b]_g=-\div J((\delta_gb,a)+(\framebox(8,8){}\
b,a)- (\delta_ga,b)-(\framebox(8,8){}\ a,b)).
$$
This equality and  relation (7.7) imply tha formula (7.5).

\vspace{3mm}

{\bf Proof of the theorem 7.2} follows at once from the theorem
7.1 and formulas (7.1), (7.5). Indeed, using conjugacy of
operators $\delta_g^*=\alpha_g$, $J^*=-J$, $\grad^*=-\div$, one
can just see, that
$$
([a, b]_g^V,[a,b]_g^V)_g=\left(\alpha_gJ\ \grad\
\mathrm{E}_g^{-1}\div J\{a,b\}, \ \alpha_gJ\ \grad\
\mathrm{E}_g^{-1}\div J\{a,b\}\right)=
$$
$$
=\int_M\left(\div J\{a,b\}\right)\left(\mathrm {E}_g^{-1}(\div
J\{a,b\}) \right)d\mu_g.
$$

\newpage

\vspace {10mm}

\centerline {\bf \S 8. The Spaces of the associated metrics on
sphere and torus.}

\vspace {7mm}

\centerline {\bf 8.1. Associated metrics on a sphere $S^2$.}

\vspace {3mm}

In this section we shall find sectional curvatures of the space
$\mathcal{AM}$ of associated metrics on a sphere $S^2$ with a
natural symplectic structure. In particular, it will be shown,
that the sectional curvature of the space $\mathcal{AM}$ is
negative.

\vspace{3mm} {\bf 8.1.1. A view of associated metrics.} Consider
unit sphere $S^2$ in ${\bf R}^3$, with spherical coordinates
$$
\left.
\begin{array}{l}
x=\cos\varphi\cos\theta\\
y=\sin\varphi\cos\theta\\
z=\sin\theta
\end {array}
\right\},\qquad\varphi\in (0,2\pi),\ \theta\in (-\pi/2,\pi/2).
$$
Then the canonical metric $g_0$ on $S^2$ and volume form look
like:
$$
g_0=\cos^2\theta d\varphi^2+d\theta^2,\quad \mu=\omega=\cos\theta
d\varphi\wedge d\theta,\quad J_0=\left(
\begin{array}{cc}
0 & -\cos^{-1} \theta \\
\cos \theta & 0
\end {array}
\right).
$$
Let $z=\sin\theta$, then
$$
g_0=(1-z^2)d\varphi^2+\frac{1}{1-z^2}dz^2,\quad\mu=\omega=d\varphi\wedge
dz, \quad J_0=\left(
\begin{array}{cc}
0 & - (1-z^2)^{-1} \\
( 1-z^2) & 0
\end {array}
\right).
$$

Choose an orthonormalized frame
$$
{\bf
e}_1=\frac{1}{\sqrt{1-z^2}}\frac{\partial}{\partial\varphi},\quad
{\bf e}_2=\sqrt{1-z^2}\frac{\partial}{\partial z}
$$
and dual coframe
$$
{\bf e}_1^*=\sqrt {1-z^2}\ d\varphi,\quad {\bf e}_2^*
=\frac{1}{\sqrt{1-z^2}}\ dz.
$$
Then $g_0=({\bf e}_1^*)^2+({\bf e}_2^*)^2$, \
$J_0=\left(\begin{array}{cc}
0 & -1 \\
1 & 0
\end {array}
\right) $.

In a corresponding complex frame
$$
\partial=\frac 12({\bf e}_1-i{\bf e}_2),\quad\overline{\partial}=
\frac 12({\bf e}_1+i{\bf e}_2),
$$
$$
d w={\bf e}_1^*+i{\bf e}_1^*,\qquad d \overline{w}={\bf
e}_1^*-i{\bf e}_1^*
$$
matrixes of tensors $g_0=dw\ d\overline{w}$ and $J_0$ look like:
$$
g_0=\frac 12\left(\begin{array}{cc}
0 & 1 \\
1 & 0
\end {array}
\right), \qquad J_0 =\left(\begin{array}{cc}
i & 0 \\
0 & -i
\end{array}
\right).
$$
The endomorphism $P$ in this case is defined by one
complex-valued function $p(z,\varphi)$ on a coordinate map
$U=(-1,1)\times (0,2\pi)$. We obtain the following values,
$$
P=\left(\begin{array}{cc}
0 & \overline {p} \\
p & 0
\end {array}
\right),\qquad P^2=p\overline{p}\left(\begin{array}{cc}
1 & 0 \\
0 & 1
\end {array}
\right),\quad 1-P^2=(1-|p|^2)\left(\begin{array}{cc}
1 & 0 \\
0 & 1
\end{array}
\right) > 0,\quad\mbox{if}\ |p|<1,
$$
$$
D=(1-P^2)^{-1}=\frac{1}{1-|p|^2}\left(\begin{array}{cc}
1 & 0 \\
0 & 1
\end{array}
\right),\quad (1+P)^2=\left(\begin{array}{cc}
1 + p\overline{p} & 2\overline{p}\\
2p & 1+p\overline{p}
\end{array}\right).
$$
Associated metric and a.c.s. $J$, corresponding to an operator
$P$ look like:
$$
g(P)=D^{-1}(1+P)g_0(1+P) =
\frac{1}{2(1-p\overline{p})}\left(\begin{array}{cc}
2p & 1 + p\overline {p} \\
1+p\overline{p} & 2\overline{p}
\end{array}\right),\eqno{(8.1)}
$$
$$
J=J_0D(1+P)^2=\frac{i}{1-p\overline{p}}\left(\begin{array}{cc}
1+p\overline{p} & 2\overline{p}\\
2p & 1+p\overline{p}
\end{array}\right).\eqno{(8.2)}
$$
$$
\partial(J)=\partial-p\overline{\partial},\qquad
\overline{\partial}(J)=\overline{\partial}-\overline{p}\partial.
$$

The anti-Hermitian symmetric 2-form $a$ is also set by one
complex-valued function:
$$
a=\left(\begin{array}{cc}
\alpha & 0 \\
0 & \overline{\alpha}
\end{array}\right),\qquad
A=g_0^{-1}a=2\left(\begin{array}{cc}
0 & \overline{\alpha}\\
\alpha & 0
\end{array}\right).
$$
Let
$$
\Psi_{AM}:\ \mathcal{P}_{J_0}\longrightarrow\mathcal{AM},\quad
P\to g=g_0(1+P)(1-P)^{-1},
$$
$$
g(X,Y)=g_0(X,(1+P)(1-P)^{-1}Y).
$$
is global parametrization of the space $\mathcal{AM}$ of
associated metrics and
$$
d\ \Psi_{AM}:\ \mathrm{End}_{SJ_0}(TM)\longrightarrow
T_g\mathcal{AM},
$$
$$
d \Psi_{AM}(A)=h_A=2g_0(1-P)^{-1}A(1-P)^{-1}.
$$
is differential of mapping $\Psi_{AM}$.

Then the element $h\in T_g\mathcal{AM}$ has the following
expression via function $p$ of an operator $P$ and function
$\alpha$ of an operator $A$:
$$
A=\left(\begin{array}{cc}
0 & \overline{\alpha}\\
\alpha & 0
\end{array}\right)
\ \mapsto\ h_A=2g_0(1-P)^{-1}A(1-P)^{-1}=\frac{1}{(1-|p|^2)^2}
\left(\begin{array}{cc}
\alpha+p^2\overline{\alpha} & p\overline{\alpha}+\overline{p}\alpha\\
p\overline{\alpha}+\overline{p}\alpha &
\overline{\alpha}+\overline{p}^2\alpha
\end{array}\right).\eqno{(8.3)}
$$
Operator $H_A=g^{-1}h_A$:
$$
H_A=\frac{2}{(1-|p|^2)^2}\left(\begin{array}{cc}
p\overline{\alpha}-\overline{p}\alpha & \overline{\alpha}-\overline{p}^2\alpha \\
\alpha-p^2\overline{\alpha} &
\overline{p}\alpha-p\overline{\alpha}
\end{array}\right).\eqno{(8.4)}
$$

Write out expressions of associated metric $g(P)$, a.c.s. $J$ and
2-form $h=2g_0P$ in real basis $\frac{\partial}{\partial z},\
\frac{\partial}{\partial\varphi}$. We shall write the function
$p(z,\varphi)$  as $p=r\ e^{i\psi}=r(\cos\psi+ i\sin\psi)$, where
$r=r(z,\varphi)$, $\psi=\psi(z,\varphi)$ are functions on domain
$U\subset S^2$.
$$
g(P)=\frac {1}{1-r^2}\left(\begin{array}{cc}
(1-z^2)(1+r^2+2r\cos\psi) & -2r\sin\psi\\
-2r\sin\psi & \frac{1+r^2-2r\cos\psi}{1-z^2}
\end{array}\right),\eqno{(8.5)}
$$
$$
J(P)=\frac {1}{1-r^2}\left(\begin{array}{cc}
2r\sin\psi & -\frac{1+r^2-2r\cos\psi}{1-z^2}\\
(1-z^2)(1+r^2+2r\cos\psi) & -2r\sin\psi
\end{array}\right),\eqno{(8.6)}
$$
If $\alpha(z,\varphi)=u(z,\varphi)+iv(z,\varphi)$,
$$
a=g_0A=\left(\begin{array}{cc}
(1-z^2)u & -v \\
-v & -\frac{u}{1-z^2}
\end{array}\right).\eqno {(8.7)}
$$

\vspace{3mm}

{\bf 8.1.2. A curvature of the space $\mathcal{AM}(S^2)$}. The
sectional curvature $K(a,b)$ in the plane section, given by the
elements $a,b\in T_g\mathcal{AM}(S^2)$ find by the formula:
$$
K(a,b)=\frac {1}{4\|a\wedge b\|^2}\int_M\tr([A, B]^2)d\mu,
$$
where $M=S^2$, $A=g_{-1}a$, $\mu=\omega=d\varphi\wedge dz$.

Find expression of a sectional curvature in a coordinate map
$\Psi_{AM}$. Let $g_0$ is canonical metric,
$g=g_0(1+P)(1-P)^{-1}$ is any associated metric and $J$ is
complex structure, corresponding to it.

Take an arbitrary elements $A,B\in\mathrm{End}_{SJ_0}(TM)$. In a
complex frame $\partial_1$, $\overline{\partial}_1$ they are set
by matrixes of a view
$$
A=\left(\begin{array}{cc}
0 & \overline {\alpha}\\
\alpha & 0
\end{array}\right),\qquad
B=\left(\begin{array}{cc}
0 & \overline {\beta}\\
\beta & 0
\end{array}\right),
$$
where $\alpha$, $\beta $ are complex functions of variables $z,\
\varphi$.

Tangent elements $h_A$, $h_B$ of the space $T_g
\mathcal{AM}(S^2)$ correspond to the operators $A$ and $B$:
$$
h_A=2g(1-P)(1-P^2)^{-1}A(1-P)^{-1},\quad
h_B=2g(1-P)(1-P^2)^{-1}B(1-P)^{-1}
$$
and operators $H_A=g^{-1}h_A$, $H_B=g^{-1}h_B$:
$$
H_A=2(1-P)(1-P^2)^{-1}A(1-P)^{-1},\quad
H_B=2(1-P)(1-P^2)^{-1}B(1-P)^{-1}.
$$
Operator bracket has a view:
$$
[H_A,H_B]=4(1-P)[(1-P^2)^{-1}A,1-P^2)^{-1}B](1-P)^{-1}.
$$

In a two-dimensional case we have,
$(1-P^2)^{-1}=\frac{1}{1-|p|^2}\left(\begin{array}{cc}
1 & 0 \\
0 & 1
\end{array}
\right)$, therefore operators $A$ and $B$ commute with
$(1-P^2)^{-1}$.

Begin to calculate a sectional curvature on the plane $\sigma $,
generated by the elements $h_A, h_B\in T_g\mathcal{AM}(S^2)$.
$$
\|h_A\|^2=(h_A,h_A)_g=4\int_M\
\tr\left((1-P^2)^{-1}A(1-P^2)^{-1}A\right)d\mu
=4\int_M\frac{\tr\left(A^2\right)}{(1-|p|^2)^2}d\mu.
$$
Therefore
$$
\|h_A\wedge h_B\|^2=\|h_A\|^2\|h_B\|^2-(h_A,h_B)^2_g=
$$
$$
=16\left(\int_M\frac{\tr\left(A^2\right)}{(1-|p|^2)^2}d\mu \
\int_M\frac{\tr\left(B^2\right)}{(1-|p|^2)^2}d\mu -
\left(\int_M\frac{\tr\left(AB\right)}{(1-|p|^2)^2}d\mu\right)^2\right).
$$
Similarly,
$$
\int_M\tr([H_A,H_B]^2)d\mu=16\int_M\frac{\tr\left([A,B]^2\right)}{(1-|p|^2)^4}
d\mu.
$$
Therefore formula for calculation of a sectional curvature
becomes:
$$
K_\sigma= \frac
14\frac{\int_M\frac{\tr\left([A,B]^2\right)}{(1-|p|^2)^4}d\mu}
{\int_M\frac{\tr\left(A^2\right)}{(1-|p|^2)^2}d\mu \
\int_M\frac{\tr\left(B^2\right)}{(1-|p|^2)^2}d\mu-\left(\int_M
\frac{\tr\left(AB\right)}{(1-|p|^2)^2}d\mu\right)^2}.\eqno {(8.9)}
$$

As the operators $A$ and $B$ are set by functions $\alpha $ and
$\beta $, then values $\tr\left([A,B]^2\right)$, \ $\tr\left
(A^2\right)$, \ $\tr\left (A^2\right)$, \ $\tr\left (AB\right)$
in this formula can also be expressed through $\alpha $ and
$\beta $.
$$
A=\left(\begin{array}{cc}
0 & \overline{\alpha}\\
\alpha & 0
\end{array}\right),\quad
A^2=\left(\begin{array}{cc}
\alpha\overline {\alpha} & 0 \\
0 & \alpha\overline{\alpha}
\end{array}\right)=|\alpha |^2E,\quad\tr A^2=2|\alpha |^2.
$$
$$
AB=\left(\begin{array}{cc}
0 & \overline{\alpha}\\
\alpha & 0
\end{array}\right)\left(\begin{array}{cc}
0 & \overline{\beta}\\
\beta & 0
\end{array}\right)=\left(\begin{array}{cc}
\overline{\alpha} \beta & 0 \\
0 & \alpha\overline{\beta}
\end{array}\right),\quad
\tr AB=\overline{\alpha}\beta+\alpha\overline{\beta}= 2{\rm Re }
(\overline{\alpha}\beta).
$$
$$
[A,B]=AB-BA =\left(\begin{array}{cc}
\overline{\alpha} \beta - \alpha\overline{\beta} & 0 \\
0 & \alpha\overline{\beta} -\overline{\alpha} \beta
\end{array}\right)=
2i\Im (\alpha\overline{\beta})\left(\begin{array}{cc}
-1 & 0 \\
0 & 1
\end{array}\right).
$$
$$
[A,B]^2=
-4\left(\Im(\alpha\overline{\beta})\right)^2\left(\begin{array}{cc}
1 & 0 \\
0 & 1
\end{array}\right),\quad\tr [A,B]^2=
-8\left(\Im(\alpha\overline{\beta})\right)^2.
$$

{\bf Theorem 8.1.} {\it Let $A,B\in\mathrm{End}_{SJ_0}(TM)$ are
operators given by functions $\alpha$ and $\beta$ and $h_A,h_B$
are elements of the tangent space $T_g\mathcal{AM}(S^2)$,
corresponding to them. Then the sectional curvature $K_\sigma$ of
the space $\mathcal{AM}(S^2)$ in the plane section $\sigma$,
generated by the elements $h_A,h_B\in T_g\mathcal{AM}(S^2)$ is
expressed by the formula
$$
K_\sigma=-\frac
12\frac{\int_M\frac{\left(\Im(\alpha\overline{\beta})\right)^2}
{(1-|p|^2)^4}d\mu}{\int_M\frac{|\alpha |^2}{(1-|p|^2)^2}d\mu\
\int_M \frac{|\alpha |^2}{(1-|p|^2)^2}d\mu-\left(\int_M \frac{{\rm
Re}(\overline{\alpha}\beta)}{(1-|p|^2)^2}d\mu\right)^2}.\eqno{(8.10)}.
$$
In particular:

1) If functions $\alpha $ and $\beta $ are simultaneously either
real, or pure imaginary, then $K_\sigma=0$.

2) If one of functions $\alpha$ and $\beta$ is real, and other is
pure imaginary, then
$$
K_\sigma= -\frac
12\frac{\int_M\frac{\left(\Im(\alpha\beta)\right)^2}
{(1-|p|^2)^4}d\mu}{\int_M\frac{|\alpha |^2}{(1-|p|^2)^2}d\mu \
\int_M\frac{|\beta |^2}{(1-|p|^2)^2}d\mu}.\eqno{(8.11)}
$$

3) The holomorphic sectional curvature is limited from above by
negative constant:
$$
K(h_A,\mathbf{J}h_A)\leq-\frac{1}{8\pi}.\eqno{(8.12)}
$$
}

{\bf Proof.} If functions $\alpha$ and $\beta$ are simultaneously
either real, or pure imaginary, then $\tr
[A,B]^2=-8\left(\Im(\alpha\overline{\beta})\right)^2=0$. In case,
when one of functions $\alpha$ and $\beta$ is real, and other is
pure imaginary, $\tr AB=2{\rm Re}(\overline{\alpha}\beta)=0$ and
$\left(\Im(\alpha\overline{\beta})\right)^2=\left(\Im(\alpha\beta)\right)^2$.

Prove the last statement. In a complex frame, when one gives
elements $A\in\mathrm{End}_{SJ_0}(TM)$ by complex functions,
operator of a complex structure on $\mathrm{End}_{SJ_0}(TM)$
$$
\mathbf{J}:\
\mathrm{End}_{SJ_0}(TM)\longrightarrow\mathrm{End}_{SJ_0}(TM),
\quad A\mapsto AJ_0
$$
acts as multiplication of function $\alpha$ on the number $i$.
Therefore if $B=\mathbf{J}A$, then $\beta=i\alpha$ and
$\tr[A,B]^2=-8\left(\Im(\alpha\overline{i\alpha})\right)^2=
-8\left(|\alpha|^2\right)^2$, $\tr AB=2{\rm
Re}(\overline{\alpha}i\alpha)=0$. Therefore the formula (8.10)
becomes
$$
K(h_A,\mathbf{J}h_A)= -\frac
12\frac{\int_M\frac{|\alpha|^4}{(1-|p|^2)^4}d\mu}
{\left(\int_M\frac{|\alpha|^2}{(1-|p|^2)^2}d\mu\
\right)^2}.\eqno{(8.12)}
$$
Now we use a Cauchy-Bunyakovskii inequality $\left(\int_Mf\
d\mu\right)^2\leq\int_Mf^2d\mu\int_Md\mu$ for function
$f=\frac{|\alpha|^2}{(1-|p|^2)^2}$. As $\int_Md\mu=4\pi$ for a
two-dimensional unit sphere $M=S^2$, then inequality
$$
\frac {1}{4\pi}\leq \frac{\int_M f^2 d\mu}{\left(\int_M f\
d\mu\right)^2}.
$$
implies the required evaluation
$K(h_A,\mathbf{J}h_A)\leq-\frac{1}{8\pi}$.


Find sectional curvatures without using of a parametrization of
the space $\mathcal{AM}$. Let $g\in\mathcal {AM}$ is any
associated metric and $J$ is complex structure, corresponding to
it. The volume form $\mu(g)$ coincides with the symplectic form
$\omega=d\varphi\wedge dz$.

Choose on a coordinate map $(\varphi,z)$ on a sphere $S^2$ a field
of orthonormalized (with respect to $g$) frames $e_1$, $e_2$. Let
$\partial_1=\frac 12(e_1-ie_2)$, $\overline{\partial}_1=\frac
12(e_1+ie_2)$ is corresponding field of complex frames.

We already marked, that in a complex frame the element $a\in T_g
\mathcal{AM}$ has a matrix $a=\left(\begin{array}{cc}
\alpha & 0\\
0 & \overline{\alpha}
\end{array}\right),$
where $\alpha(z,\varphi)$ is complex function on a map $U\subset
S^2$. When one gives elements $a\in T_g \mathcal{AM}(S^2)$ by
complex functions $\alpha(z,\varphi)$ on $U\subset S^2$, the
operator of an almost complex structure
$$
\mathbf{J}:\ T_g\mathcal{AM}(S^2)\longrightarrow
T_g\mathcal{AM}(S^2)
$$
acts as multiplication of function $\alpha (z,\varphi)$ on $i$.

The functions $e^{i(kz+l\varphi)}$, $k\in\pi\mathbf{Z}$,
$l\in\mathbf{Z}$, form full orthogonal system of functions on
domain $U=\{(\varphi,z);\ \varphi\in(0,2\pi),\ z\in(-1,1)\}$. They
define complex basis in $T_g\mathcal{AM}$. The real basis of the
space $T_g\mathcal{AM}$ corresponds to functions
$$
\cos(kz+l\varphi),\ \sin(kz+l\varphi),\ i\cos(kz+l\varphi), \
i\sin(kz+l\varphi),\ k\in\pi\mathbf{Z},\
l\in\mathbf{Z}.\eqno{(8.13)}
$$

{\bf Remark.} Coordinate map $(U;\ \varphi,z)$ has singularities
in two points of a sphere $z=1$ and $z=-1$. The 2-form $a$,
defined by function $\alpha(z, \varphi)$ on $U$ can have
singularities at $z=\pm 1$ and $\varphi=2\pi$. At calculation of
sectional curvatures the integration on a sphere is used. As the
singularities of limited functions on a sphere concentrated on
set of zero measure, do not influence on values of integrals,
then the functions of the system (8.13) can be used for
calculation of sectional curvatures of the space
$\mathcal{AM}(S^2)$.

\vspace{3mm}

{\bf Theorem 8.2.} {\it If functions $\alpha(z,\varphi)$, $\beta
(z,\varphi)$ from system (8.13), which gives $a,b\in
T_g\mathcal{AM}$, are simultaneously real, or pure imaginary,
then $K(a,b)=0$. If the form $a$ is given by real function, and
form $b$ is given by pure imaginary function (8.13), then
$$
K(a,b)=\left\{\begin{array}{l}
-\frac{1}{8\pi},\quad\mbox{if}\quad |\alpha |\neq |\beta |,\\
-\frac{3}{16\pi},\quad\mbox{if}\quad |\alpha|=|\beta|\neq 1,\\
-\frac{1}{8\pi},\quad\mbox{if}\quad |\alpha |=|\beta |=1.
\end{array}\right.
$$
In particular, for any basis function $\alpha$ the holomorphic
sectional curvature $K(a,\mathbf{J}a)$ takes values
$$
K(a,\mathbf{J} a)=\left\{\begin{array}{c}
-\frac{3}{16\pi},\quad |\alpha|\neq 1,\\
-\frac{1}{8\pi},\quad |\alpha|=1.
\end{array}\right.
$$
}

{\bf Proof.} For operators $A=g^{-1}a$, $B=g^{-1}b$ corresponding
to the forms $a,b$ we have:
$$
A=2\left(\begin{array}{cc}
0 & \overline{\alpha}\\
\alpha & 0
\end{array}\right),\quad
B=2\left(\begin{array}{cc}
0 & \overline{\beta}\\
\beta & 0
\end{array}\right),\quad
[A,B]=8i\left(\begin{array}{cc}
\Im(\overline{\alpha}\beta) & 0\\
0 & -\Im(\overline{\alpha}\beta)
\end{array}\right).
$$
Thus, $[A,B]=0$ if and only if $\Im(\overline{\alpha}\beta)=0$.
In particular, so will be, if $\alpha $ and $\beta $ are
simultaneously real, or pure imaginary. The first statement of
the theorem follows from here. Further, in the second case we
have $\tr(AB)=8{\rm Re}(\overline{\alpha}\beta)=0$. From the
formulas
$$
\int_M\tr([A,B]^2)dzd\varphi=
-128\int_M\left(\Im(\overline{\alpha}\beta)\right)^2 dzd\varphi,
$$
$$
\|a\wedge b\|^2=\int_M\tr(A^2)d\mu\int_M\tr(B^2)\mu
\left(\int_M(\tr(AB)\mu)\right)^2=64\int_M |\alpha |^2
dzd\varphi\int_M |\beta |^2 dzd\varphi,
$$
we obtain expression of a sectional curvature:
$$
K_\sigma = -\frac
12\frac{\int_M\left(\Im(\alpha\beta)\right)^2d\mu} {\int_M
|\alpha |^2 d\mu\ \int_M |\beta |^2 d\mu}.
$$

If
$$
a=\cos(kz+l\varphi),\ \mbox{or}\ \sin(kz+l\varphi),
$$
$$
\beta=i\cos(pz+q\varphi),\ \mbox{or}\ i\sin(pz+q\varphi),
$$
then the direct calculations give values, indicated in the
theorem. Indeed, we consider all essentially various cases.

1) Let, for example,
$$
\alpha=\cos(kz+l\varphi),\ \beta=i\cos(pz+q\varphi),\quad
(k,l)\neq(0,0),\ (p,q)\neq (0,0),\ (k,l)\neq(p,q).
$$
Designate $\zeta=kz+l\varphi$,\quad $\psi=pz+q\varphi$.\quad Then
\quad $|\alpha |^2=\cos^2\zeta,\quad |\beta |^2=\cos^2\psi$\quad
and \quad
$\left(\Im(\overline{\alpha}\beta)\right)^2=(\cos\zeta\cos\psi)^2=
\frac 14\left(\cos(\zeta+\psi)+\cos(\zeta-\psi)\right)^2$.
$$
\int_{0}^{2\pi}\int_{-1}^{1}\left
(\Im(\overline{\alpha}\beta)\right)^2 d\varphi dz =
$$
$$
=\frac 14\int_{0}^{2\pi}d\varphi
\int_{-1}^{1}(\cos^2(\zeta+\psi)+2\cos(\zeta+\psi)\cos(\zeta-\psi)+
\cos^2(\zeta-\psi))dz=\frac 14(2\pi+2\pi)=\pi.
$$
$$
\|a\wedge b\|^2=\int_U\cos^2\zeta\ dzd\varphi\ \int_U\cos^2 \psi\
dzd\varphi = 4\pi^2.
$$
Therefore
$$
K(a,b)=-\frac{1}{2}\frac{\pi}{4\pi^2}=-\frac{1}{8\pi}.
$$

2) Let \quad $\alpha=1,\ \beta=i\cos(pz+q\varphi),\ (p,q)\neq
(0,0)$. Designate $\psi=pz+q\varphi$. Then
$\left(\Im(\overline{\alpha}\beta)\right)^2=(\cos\psi)^2$\quad
and \quad $|\alpha|^2=1,\quad |\beta|^2=\cos^2\psi$.
$$
\int_{0}^{2\pi}\int_{-1}^{1}\left(\Im(\overline{\alpha}\beta)\right)^2
d\varphi
dz=\int_{0}^{2\pi}d\varphi\int_{-1}^{1}\cos^2(\psi)dz=2\pi.
$$
$$
\|a\wedge b\|^2=\int_Udzd\varphi\ \int_U\cos^2\psi\
dzd\varphi=(4\pi)(2\pi).
$$
Therefore
$$
K(a,b)=-\frac 12\frac{2\pi}{(4\pi)(2\pi)}=-\frac{1}{8\pi}.
$$

3) Now let (holomorphic sectional curvature),
$$
\alpha=\cos(kz+l\varphi),\ \beta=i\cos(kz+l\varphi),\ (k,l)\neq
(0,0).
$$
Designate $\zeta=kz+l\varphi$. Then
$\left(\Im(\overline{\alpha}\beta)\right)^2=(\cos^2\zeta)^2=\frac
14 \left(1+\cos(2\zeta)\right)^2$ \quad and \quad $|\alpha
|^2=\cos^2\zeta,\quad |\beta |^2=\cos^2\zeta$.
$$
\int_{0}^{2\pi}\int_{-1}^{1}\left(\Im(\overline{\alpha}\beta)
\right)^2d\varphi dz = \frac 14 \int_{0}^{2\pi} d\varphi
\int_{-1}^{1}(1+2\cos(2\zeta)+\cos^2(2\zeta))dz=\frac
14(4\pi+2\pi)= \frac{3\pi}{2}.
$$
$$
\|a\wedge b\|^2=\int_U\cos^2\zeta\ dzd\varphi\ \int_U\cos^2\zeta\
dzd\varphi =(2\pi)(2\pi).
$$
Therefore
$$
K(a,ia)=-\frac 12\frac{3\pi}{2(2\pi)(2\pi)}=-\frac {3}{16\pi}.
$$

4) The last case (holomorphic sectional curvature): $\alpha=\pm
1,\ \beta=\pm i$. Then
$\left(\Im(\overline{\alpha}\beta)\right)^2=1$, \quad $|\alpha
|^2=1,\quad |\beta |^2=1$. Therefore
$$
\int_{0}^{2\pi} \int_{-1}^{1} \left(\Im (\overline{\alpha}\beta)
\right)^2d\varphi dz = 4\pi, \quad \|a\wedge b \|^2 = \int_U
dzd\varphi\ \int_U dzd\varphi = (4\pi)(4\pi).
$$

$$
K(1,i)=-\frac 12\frac{4\pi}{(4\pi)(4\pi)}=-\frac{1}{8\pi}.
$$
The theorem is proved.

\vspace{5mm}

Finally we shall write out the equation $\div J \delta_g a = 0$,
when $g=g_0$. The covariant divergence of the form $a$ is a
vector field
$$
\delta_g a = -2\left ((1-z^2)\frac {\partial u}{\partial z} -
\frac{\partial v}{\partial \varphi} - 2zu\right)
\frac{\partial}{\partial z} + 2\left(\frac{\partial v}{\partial
z} + \frac{1}{1-z^2}\frac{\partial u}{\partial\varphi} -
\frac{2z}{1-z^2} v \right)\frac{\partial}{\partial\varphi}.
$$
Calculating value $\div J\delta_ga$ we obtain the equation, which
pickes out horizontal directions $h$:
$$
-\frac 12\div J \delta_g a =
$$
$$
( 1-z^2) \frac{\partial^2 v}{\partial z^2} -
\frac{1}{1-z^2}\frac{\partial^2 v}{\partial \varphi^2} -
4z\frac{1}{1-z^2}\frac{\partial v}{\partial z} - 2v +
2\frac{\partial^2 u}{\partial z\partial \varphi} -
\frac{2z}{1-z^2}\frac{\partial u}{\partial\varphi}=0.\eqno{(8.14)}
$$

\newpage

\vspace {10mm}

\centerline {\bf 8.2. Associated metrics on a torus $T^2$.}

\vspace {3mm}

In this section we shall find sectional curvatures of the space
$\mathcal{AM}$ of associated metrics on a torus $T^2$ with a
natural symplectic structure. In particular, values of sectional
curvatures in regular points of the quotient space
$\mathcal{AM}(T^2)/\mathcal{G}$, where $\mathcal{G}$ is the group
of symplectic diffeomorphisms of a torus, which Lie algebra is
algebra of Hamilton vector fields on $T^2$, will be found.

\vspace {3mm} {\bf 8.2.1. A view of associated metrics}. Consider
a torus $T^2={\bf R}^2/2\pi{\bf Z}^2$ with coordinates $(x,y)mod\
2\pi$, the flat metric $g_0=dx^2+dy^2$ and complex structure
$J_0$: $z=x+iy$. A volume form: $\mu=dx\wedge dy$. In complex
basis
$$
\frac{\partial}{\partial z}=\frac
12\left(\frac{\partial}{\partial x}- i\frac{\partial}{\partial
y}\right),\quad\frac{\partial}{\partial\overline{z}}= \frac
12\left(\frac{\partial}{\partial x}+i\frac{\partial}{\partial
y}\right).
$$
matrixes of tensor $g_0=dz\ d\overline{z}$ and operator of a
complex structure $J_0$ look like the same ones in case of
sphere. The endomorphism $P$ is defined by one complex-valued
function $p(x,y)$ on a torus $T^2$. Associated metric and a.c.s.
$J$, corresponding to an operator, $P$ have the same form, as
well as in case of a sphere $S^2$.
$$
\partial (J)=\frac {\partial}{\partial z} -
\overline{p}\frac {\partial}{\partial \overline{z}},\quad
\overline{\partial}(J)=\frac {\partial}{\partial \overline{z}} -
p\frac {\partial}{\partial z}.
$$

Write out expressions of associated metric $g(P)$, a.c.s. $J$ and
2-form $h=2g_0P$ in real basis $\frac {\partial}{\partial x},\
\frac {\partial}{\partial y}$. Write function $p(x,y)$ as $p=r\
e^{i\psi}=r(\cos\psi+i\sin\psi)$, where $r=r(x,y)$,
$\psi=\psi(x,y)$ is function on a torus $T^2$.
$$
g(P)=\frac{1}{1-r^2}\left(\begin{array}{cc}
1+r^2+2r\cos\psi & -2r\sin\psi\\
-2r\sin\psi & 1+r^2-2r\cos\psi
\end{array}\right),
$$
$$
J(P)=\frac{1}{1-r^2}\left(\begin{array}{cc}
2r\sin\psi & -(1+r^2-2r\cos\psi)\\
1+r^2+2r\cos\psi & -2r\sin\psi
\end{array}\right).
$$

{\bf 8.2.2. A curvature of the space $\mathcal{AM}(T^2)$.} The
expression of sectional curvature in a coordinate map $\Psi_{AM}$
is the same, as in case of a sphere. Let $g_0$ is canonical
metric, $g=g_0(1+P)(1-P)^{-1}$ is any associated metric and $J$
is complex structure, corresponding to it.

Take the arbitrary elements $A,B\in\mathrm{End}_{SJ_0}(TM)$. In a
complex frame $\partial_1$, $\overline{\partial}_1$ they are set
by matrixes of the view
$$
A=\left(\begin{array}{cc}
0 & \overline{\alpha}\\
\alpha & 0
\end{array}\right),\qquad
B=\left(\begin{array}{cc}
0 & \overline{\beta}\\
\beta & 0
\end{array}\right),
$$
where $\alpha$, $\beta $ are complex functions of variables $x,y$.

\vspace {3mm} {\bf Theorem 8.3.} {\it Let
$A,B\in\mathrm{End}_{SJ_0}(TM)$ are operators, given by functions
$\alpha$ and $\beta$ and $h_A, h_B$ are the elements  of tangent
space $T_g\mathcal{AM}(T^2)$, corresponding to them. Then the
sectional curvature $K_\sigma$ of the space $\mathcal{AM}(S^2)$
in the plane section $\sigma$, generated by the elements
$h_A,h_B\in T_g\mathcal{AM}(T^2)$, is expressed by the formula
$$
K_\sigma=-\frac
12\frac{\int_M\frac{\left(\Im(\alpha\overline{\beta})
\right)^2}{(1-|p|^2)^4} d\mu}{\int_M \frac {|\alpha
|^2}{(1-|p|^2)^2} d\mu \ \int_M \frac {|\alpha |^2}{(1-|p|^2)^2}
d\mu - \left(\int_M \frac{{\rm Re }(\overline
{\alpha}\beta)}{(1-|p|^2)^2} d\mu\right)^2}.
$$
where $M=T^2$. In particular:

1) If functions $\alpha$ and $\beta$ are simultaneously either
real, or pure imaginary, then $K_\sigma=0$.

2) If one of functions $\alpha$ and $\beta$ is real, and other
pure imaginary, then
$$
K_\sigma = -\frac
12\frac{\int_M\frac{\left(\Im(\alpha\beta)\right)^2}
{(1-|p|^2)^4} d\mu}{\int_M \frac{|\alpha |^2}{(1-|p|^2)^2}d\mu \
\int_M\frac{|\beta |^2}{(1-|p|^2)^2} d\mu}.
$$

3) The holomorphic sectional curvature is limited from above by
negative constant:
$$
K(h_A,\mathbf{J} h_A)\leq-\frac{1}{8\pi^2}.
$$
}

{\bf Proof.} The reason of appearance of number $\pi^2$, instead
of $\pi$, is that the area of a sphere is equal to $4\pi$, and
the area of our torus is equal to $4\pi^2$.

\vspace {3mm}

Find sectional curvatures without using of a parametrization of
the space $\mathcal{AM}$. Let $g\in\mathcal{AM}$ is any
associated metric and $J$ is complex structure, corresponding to
it. The volume form $\mu(g)$ coincides with the symplectic form
$\omega=dx\wedge dy$.

Choose a field of orthonormalized (with respect to $g$) frames on
torus $T^2$ $e_1$, $e_2$. Let $\partial_1=\frac 12(e_1-ie_2)$,
$\overline{\partial}_1=\frac 12(e_1+ie_2)$ is field of
corresponding complex frames.

We already marked, that in a complex frame the element $a\in
T_g\mathcal{AM}$ has a matrix $a=\left(\begin{array}{cc}
\alpha & 0\\
0 & \overline{\alpha}
\end{array}\right),$
where $\alpha(x,y)$ is complex $2\pi$ -- periodic function of
variables $x,y$. When one gives elements $a\in
T_g\mathcal{AM}(T^2)$ by complex functions $\alpha(x,y)$,
operator of an almost complex structure
$$
\mathbf{J}:\ T_g\mathcal{AM}(T^2)\longrightarrow
T_g\mathcal{AM}(T^2)
$$
acts as multiplication of function $\alpha (x,y)$ on $i$.

The full orthogonal system of functions on a torus $T^2$ is formed
by functions
$$
p_{kl}=e^{i(kx+ly)},\qquad k,l\in{\bf Z},\ (x,y)\in T^2.
$$
One can find a curvature $R=\frac{1}{det\ g(p)}R_{1212}$ of
associated metric $g(p_{kl})$.

The direct count in a system of analytical calculations MapleV
gives for a set of the metrics $g_t=g(tp_{kl})$ the following
expression of a Gaussian curvature:
$$
R=-\frac{t}{1-t^2}\left((k^2-l^2)\cos(kx+ly)-2kl\sin(kx+ly)\right).
$$

The functions $e^{i(kx+ly)}$ define the complex basis in
$T_g\mathcal {AM}$. The real basis of the space $T_g\mathcal{AM}$
corresponds to functions
$$
\cos(kx+ly),\ \sin(kx+ly),\ i\cos(kx+ly),\ i\sin(kx+ly),\quad
k,l\in{\bf Z}. \eqno{(8.15)}
$$

\vspace {3mm} {\bf Theorem 8.4.} {\it If functions $\alpha
(x,y)$, $\beta (x,y)$ from system (8.15), which give 2-forms
$a,b\in T_g\mathcal{AM}$ are simultaneously either real, or pure
imaginary, then $K(a,b)=0$. If the form $a$ is set by real
function, and form $b$ is pure imaginary function (8.15), then
$$
K(a,b)=\left\{\begin{array}{l}
-\frac{1}{8\pi^2},\quad\mbox{if}\quad|\alpha|\neq|\beta|,\\
-\frac{3}{16\pi^2},\quad\mbox{if}\quad|\alpha|=|\beta|\neq 1,\\
-\frac{1}{8\pi^2},\quad\mbox{if}\quad|\alpha|=|\beta|=1.
\end{array}\right.
$$
In particular, for any basis function $\alpha$ the holomorphic
sectional curvature $K(a,\mathbf{J}a)$ takes values
$$
K(a,\mathbf{J}a)=\left\{\begin{array}{c}
-\frac{3}{16\pi^2},\quad|\alpha|\neq 1,\\
-\frac{1}{8\pi^2},\quad|\alpha|=1.
\end{array}\right.
$$
In case, when $\alpha (x,y)=e^{i(kx+ly)}$,
$\beta(x,y)=e^{i(px+qy)}$,
$$
K(a,b)=\left\{\begin{array}{l}
-\frac{1}{16\pi},\quad\mbox{if}\quad\beta\neq\pm i\alpha,\\
-\frac{1}{8\pi},\quad\mbox{if}\quad\beta=\pm i\alpha.
\end{array}\right.
$$
}

{\bf Proof.} Let $\alpha=e^{i\varphi}$, $\beta=E^{i\psi}$, where
$\varphi,\psi$ are functions on $T^2$ of view $\varphi=kx+ly$,
$\psi=px+qy$, we assume, that $\varphi\neq\psi$. In this case we
have the following expressions:
$$
|\alpha|^2=1,\quad|\beta|^2=1,
$$
$$
\quad\left(\Im(\overline{\alpha}\beta)\right)^2=
\left(\Im(e^{i(\psi-\varphi)})\right)^2=\sin^2(\psi-\varphi)=
\frac{1-\sin 2(\psi-\varphi)}{2}.
$$

$$
\int_{0}^{2\pi}\int_{0}^{2\pi}\left(\Im(\overline{\alpha}\beta)\right)^2dxdy
= \int_{0}^{2\pi}\int_{0}^{2\pi}\frac{1-\sin 2(\psi-\varphi)}{2}
dxdy = \left\{\begin{array}{l}
2\pi^2,\quad\mbox{if}\quad\psi-\varphi\neq\pm\pi/2,\\
4\pi^2,\quad\mbox{if}\quad\psi-\varphi=\pm\pi/2.
\end{array}\right.
$$
$$
\|a\wedge b\|^2=\int_{0}^{2\pi}\int_{0}^{2\pi}dxdy
\int_{0}^{2\pi} \int_{0}^{2\pi} dxdy - \left(\int_{0}^{2\pi} \int
_{0}^{2\pi}\cos (\psi-\varphi)dxdy\right)^2 = (4\pi^2)(4\pi^2).
$$
Therefore
$$
K(a,b)=\left\{\begin{array}{l}
-\frac{1}{2}\frac{2\pi^2}{(4\pi^2)(4\pi^2)}=-\frac{1}{16\pi},\quad
\mbox{if}\quad\beta\neq\pm i\alpha,\\
-\frac{1}{2}\frac{4\pi^2}{(4\pi^2)(4\pi^2)}=-\frac{1}{8\pi},\quad
\mbox{if}\quad\beta=\pm i\alpha.
\end{array}\right.
$$
The theorem is completely proved.

\vspace {3mm} {\bf 8.2.3. A curvature of the space
$\mathcal{AM}/\mathcal{G}$.} Let $\overline{a},\overline{b}\in
T_{[g]}(\mathcal{AM}/\mathcal{G})$. As $p:\
\mathcal{AM}\longrightarrow\mathcal{AM}/\mathcal {G}$ is
Riemannian submersion, the sectional curvature
$\overline{K}(\overline{a},\overline{b})$ of the space
$\mathcal{AM}/\mathcal{G}$ is under the formula \cite{ONe},
\cite{Bes2}:
$$
\overline{K}(\overline{a},\overline{b})= K(a,b)+\frac
34\frac{\left([a,b]^V,[a,b]^V\right)_g}{\|a\wedge b\|^2}, \eqno
{(8.16)}
$$
where $a,b\in T_{g}\mathcal{AM}$ are horizontal lifts of vectors
$\overline{a},\overline{b}$, $K(a,b)$ is sectional curvature of
the space $\mathcal{AM}$, $[a,b]^V$ is vertical part of Lie
commutator of horizontal vector fields on $\mathcal{AM}$, which
continue $a,b\in T_g\mathcal{AM}$.

Consider the flat metric $g_0=dz\ d\overline{z}$ on $T^2$. Let
symmetric 2-form $h$ in real basis $dx,dy$ has a view
$$
h=u\ dx^2-2v\ dxdy-u\ dy^2,\eqno{(8.17)}
$$
Corresponding quadratic differential $h=pdzdz+\overline{p}\
d\overline{z}d\overline{z}$ is set by function $p=1/2(u+iv)$.

The covariant divergence $\delta_gh=-\nabla^ih_{ij}$ represents a
vector field
$$
\delta_gh=-2\left(\frac{\partial u}{\partial x}-\frac{\partial
v}{\partial y} \right)\frac{\partial}{\partial
x}+2\left(\frac{\partial v}{\partial x}+ \frac{\partial
u}{\partial y}\right)\frac{\partial}{\partial y}. \eqno {(8.18)}
$$
Mark incidentally, that the equality $\delta_gh=0$ means, that
$h=p\ dzdz$ is holomorphic quadratic differential. Find $\div
J\delta_gh$.
$$
J\delta_gh=-2\left(\frac{\partial v}{\partial x}+\frac{\partial
u}{\partial y} \right)\frac{\partial}{\partial
x}-2\left(\frac{\partial u}{\partial x}- \frac{\partial
v}{\partial y}\right)\frac{\partial}{\partial y}.
$$
$$
-\frac 12\div J\delta_gh=\frac{\partial^2v}{\partial x^2}-
\frac{\partial^2 v}{\partial y^2}+2\frac{\partial^2 u}{\partial
x\partial y}. \eqno {(8.19)}
$$
In the complex form:
$$
\div
J\delta_gh=-8\Im\left(\frac{\partial^2p}{\partial\overline{z}^2}\right).
\eqno {(8.20)}
$$
The horizontal space $S^*_{2A}$, orthogonal to an orbit
$g_0\mathcal{G}\subset\mathcal{AM}$ consists of quadratic
differentials $a=pdz^2$, satisfying to the equation,
$$
\frac{1}{4}\Im\left(\frac{\partial^2p}{\partial\overline{z}^2}\right)=
\frac{\partial^2 v}{\partial x^2}-\frac{\partial^2v}{\partial
y^2}+ 2\frac{\partial^2u}{\partial x\partial y}=0.\eqno {(8.21)}
$$

The basis of solutions of this equation is given by complex
functions $p(x,y)$ of a view
$$
p(x,y)=\left(f_1(x)+f_2(y)\right) + i\left(f_3(x+y) +
f_4(x-y)\right).
$$
Basis of real functions $p=u(x,y)$, defining horizontal forms,
i.e., satisfying to the equation,
$$
\frac{\partial^2u}{\partial x\partial y}=0,
$$
form the functions of a view
$$
\cos kx,\ \sin kx,\ \cos ly,\ \sin ly,\quad k,l\in{\bf
Z}.\eqno{(8.22)}
$$
Basis of pure imaginary functions $p=iv(x,y)$, defining
horizontal forms, i.e., satisfying to the equation,
$$
\frac{\partial^2 v}{\partial x^2}-\frac {\partial^2 v}{\partial
y^2} = 0,
$$
is formed by functions of a view
$$
i\cos p(x+y),\ i\sin p(x+y),\ i\cos q(x-y),\ i\sin q(x-y),\quad
p,q\in{\bf Z}. \eqno{(8.23)}
$$
Find a sectional curvature of the quotient space
$\mathcal{AM}/\mathcal {G}$ in case, when the horizontal forms
$a,b$, given by functions of view (8.22) and (8.23),  are taken
as the elements $a,b\in T_{[g]}\mathcal {AM}/\mathcal {G}$.

\vspace {3mm} {\bf Theorem 8.5.} {\it Sectional curvature
$\overline{K}(\overline{a},\overline{b})$ of the space
$\mathcal{AM}/\mathcal{G}$ at a point $[g_0]$ in direction of the
horizontal forms $a,b\in T_{g_0}\mathcal {AM}$ takes the
following values:

1) If $\alpha$ is $\cos kx$ or $\sin kx$, and $\beta$ is $\cos
ly$ or $\sin ly$, then
$$
\overline{K}(\overline{a},\overline{b})=\frac{3}{8\pi^2}
\frac{k^2l^2}{(k^2+l^2)^2},
$$
in remaining cases of real functions,
$\overline{K}(\overline{a},\overline{b})=0$,

2) If $\alpha$ is $i\cos p(x+y)$ or $i\sin p(x+y)$, and $\beta$ is
$i\cos q(x-y)$ or $i\sin q(x-y)$, then
$$
\overline{K}(\overline{a},\overline{b})=
\frac{3}{32\pi^2}\frac{p^2q^2}{(p^2+q^2)^2},
$$
in remaining cases of pure imaginary functions,
$\overline{K}(\overline{a},\overline{b})=0$,

3) If $\alpha$ is any real function from (8.22) and $\beta$ is
imaginary function of a view (8.23), then
$$
\overline{K}(\overline{a},\overline{b})=\left\{\begin{array}{l}
\frac{1}{4\pi^2},\quad\mbox{if}\quad|\alpha|\neq 1,|\beta|\neq 1,\\
-\frac{1}{8\pi},\quad\mbox{if}\quad|\alpha|=1\ \mbox{or}\
|\beta|=1.
\end{array}\right.
$$}

{\bf Proof.} We shall use the formula (7.6),
$$
\overline{K}(\overline{a},\overline{b})=K(a,b)+\frac{3}{4\|a\wedge
b\|^2} \int_M\left(\div J\{a,b\}\right)\left(\mathrm
{E}_g^{-1}(\div J\{a,b\})\right) d\mu_g,
$$
where:
$$
\mathrm{E}_g(f)=\div J\delta_g\alpha_gJ\grad f,
$$
$$
\{a, b\}=A(\delta_g b)+(\framebox(8,8){}\ (b),a)- B(\delta_g
a)(\framebox(8,8){}\ (a),b),
$$
$$
A(\delta_gb)= a_i^k(\delta_gb)^i,\quad(\delta_g
b)^i=-\nabla_jb^{ij},
$$
$$
(\framebox (8,8){}\ (b),a)=b^k_{ij}a^{ij},
$$
$$
(\framebox (8,8){}\ (b))^k_{ij}=b^k_{ij}= \frac
12\left(\nabla_ib_i^k+\nabla_jb_i^k-\nabla^kb_{ij}\right),\quad
b\in S^2.
$$
In our case of a two-dimensional torus $T^2$ and flat metric
$g=dx^2+dy^2$, the operator $\mathrm{E}_g$ has a view:
$$
\mathrm{E}_gf=\frac 12\left(\frac{\partial^4f}{\partial x^4}+
2\frac{\partial^4f}{\partial x^2\partial
y^2}+\frac{\partial^4f}{\partial y^4} \right)=\frac 12\Delta^2f.
$$

It is convenient now to set the symmetric 2-forms $a,b$ in real
base $dx,dy$. The 2-form $a$, corresponding to complex function
$\alpha=1/2(u+iv)$ has a matrix $a=\left(\begin{array}{cc}
u & -v \\
-v & -u
\end{array}\right)$
in real basis. Calculate $\delta_ga$, $\framebox(8,8){}\ (a)$ for
the following basic horizontal forms.

{\bf Type 1.} $a=\left(\begin{array}{cc}
\cos kx & 0 \\
0 & -\cos kx
\end{array}\right)$.

Covariant divergence: $\delta_ga=(k\sin kx,0)$.

Tensor $a^k_{ij}$:
$$
a^1_{11}=-\frac 12 k\sin kx,\quad a^1_{12}=0,\quad
a^1_{22}=-\frac 12k\sin kx,
$$
$$
a^2_{11}=0,\quad a^2_{12}=\frac 12k\sin kx,\quad a^2_{22}=0.
$$

{\bf Type 2.} $a=\left(\begin{array}{cc}
\cos ly & 0 \\
0 & -\cos ly
\end{array}\right)$.

Covariant divergence: $\delta_ga=-(0,l\sin ly)$.

Tensor $a^k_{ij}$:
$$
a^1_{11}=0,\quad a^1_{12}=-\frac 12l\sin ly,\quad a^1_{22}=0,
$$
$$
a^2_{11}=\frac 12l\sin ly,\quad a^2_{12}=0,\quad a^2_{22}=\frac
12l\sin ly.
$$

{\bf Type 3.} $a=\left(\begin{array}{cc}
0 & -\cos p(x+y)\\
-\cos p(x+y) & 0
\end{array}\right)$.

Covariant divergence: $\delta_ga=-(p\sin p(x+y), p\sin p(x+y))$.

Tensor $a^k_{ij}$:
$$
a^1_{11}=0,\quad a^1_{12}=0,\quad a^1_{22}=p\sin p(x+y),
$$
$$
a^2_{11}=p\sin p(x+y),\quad a^2_{12}=0,\quad a^2_{22}=0.
$$

{\bf Type 4.} $a=\left(\begin{array}{cc}
0 & -\cos q(x-y)\\
-\cos q(x-y) & 0
\end{array}\right)$.

Covariant divergence: $\delta_ga=-(q\sin q(x-y),q\sin q(x-y))$.

Tensor $a^k_{ij}$:
$$
a^1_{11}=0,\quad a^1_{12}=0,\quad a^1_{22}=-q\sin q(x-y),
$$
$$
a^2_{11}=q\sin q(x-y),\quad a^2_{12}=0,\quad a^2_{22}=0.
$$

{\bf Remark.} In case, when for definition of the form $a$ the
function $\sin$ is taken, the values $\delta_ga$ and
$\framebox(8,8){}\ (a)$ are similarly found with elementary
replacement of functions $\sin$ and $\cos $.

Begin calculation of a sectional curvature
$\overline{K}(\overline{a},\overline{b})$ of the space
$\mathcal{AM}/\mathcal{G}$.

{\bf First part.} Let
$$
a =\left(\begin{array}{cc}
\cos kx & 0 \\
0 & -\cos kx
\end{array}\right),\quad
b=\left(\begin{array}{cc}
\cos ly & 0 \\
0 & -\cos ly
\end{array}\right),\quad k\neq 0,\ l\neq 0.
$$
In this case: $K(a,b)=0$,
$$
(\framebox(8,8){}\ (b),a)=0,\quad(\framebox(8,8){}\ (a),b)=0,
$$
$$
A(\delta_g b)=(0,l\cos kx \sin ly), \quad B(\delta_g a)=(k\sin
kx\cos ly,0),
$$
$$
\{a,b\}=(-k\sin kx\cos ly,l\cos kx\sin ly),
$$
$$
J\{a,b\}=-(l\cos kx\sin ly,k\sin kx\cos ly),
$$
$$
\div J \{a,b\}=2kl\sin kx\sin ly=kl(\cos(kx-ly)-\cos (kx+ly)),
$$
$$
\mathrm{E}_g^{-1}(\div J\{a,b\})= 2kl\frac
{1}{(k^2+l^2)^2}\left(\cos(kx-ly)-\cos (kx+ly)\right),
$$
$$
\int_{T^2}\left(\div J\{a,b\}\right)\mathrm{E}_g^{-1}\left(\div J
\{a, b\} \right)\ dx\
dy=\frac{2k^2l^2}{(k^2+l^2)^2}(2\pi^2+2\pi^2),
$$
$$
\|a\|^2=\int_{T^2}\tr(A^2)\ dx\ dy=2\int_{T^2}\cos^2kx\ dx\
dy=4\pi^2=\|b\|,
$$
$$
(a,b)=2\int_{T^2}\cos kx\cos ly\ dx\ dy=0,
$$
Therefore
$$
\overline{K}(\overline{a},\overline{b})=
\frac{3}{8\pi^2}\frac{k^2l^2}{(k^2+l^2)^2}.
$$
If
$$
a=\left(\begin{array}{cc}
\cos kx & 0\\
0 & -\cos kx
\end{array}\right),\quad
b=\left(\begin{array}{cc}
\sin ly & 0 \\
0 & -\sin ly
\end{array}\right),\quad k\neq 0,\ l\neq 0,
$$
then similarly,
$$
\overline{K}(\overline{a},\overline{b})=\frac{3}{8\pi^2}
\frac{k^2l^2}{(k^2+l^2)^2},
$$
Let now
$$
a=\left(\begin{array}{cc}
\cos kx & 0\\
0 & -\cos kx
\end{array}\right),\quad
b=\left(\begin{array}{cc}
\cos lx & 0\\
0 & -\cos lx
\end{array}\right).
$$
In this case: $K(a,b)=0$,
$$
(\framebox (8,8){}\ (b),a)=0,\quad(\framebox (8,8){}\ (a), b) = 0,
$$
$$
A(\delta_g b) = (\l\sin lx\cos kx, 0), \quad B (\delta_g a) =
(k\sin kx\cos lx, 0),
$$
$$
\{a,b\} = (l\sin lx\cos kx - k\sin kx\cos lx, 0),
$$
$$
\div J \{a, b \} = 0.
$$
Therefore $\overline{K}(\overline{a},\overline{b})=0$.

If
$$
a=\left(\begin{array}{cc}
\cos kx & 0 \\
0 & -\cos kx
\end{array}\right),\quad
b=\left(\begin{array}{cc}
\sin lx & 0 \\
0 & -\sin lx
\end{array}\right).
$$
then similarly, $\overline {K}(\overline {a}, \overline {b}) = 0$.

{\bf Second part.} Let
$$
a = \left (\begin {array}{cc}
0 & -\cos p (x + y) \\
-\cos p (x + y) & 0
\end {array} \right),
b = \left (\begin {array}{cc}
0 & -\cos q (x-y) \\
-\cos q (x-y) & 0
\end {array} \right), p^2 + q^2\neq 0.
$$
In this case: $K(a,b) = 0$,
$$
(\framebox (8,8){}\ (b), a) = 0, \quad (\framebox (8,8){}\ (a),
b) =0,
$$
$$
A(\delta_g b)=\left(q\cos p(x+y) \sin q(x-y), -q\cos p(x+y)\sin
q(x-y)\right),
$$
$$
B(\delta_g a)=\left (p\sin p(x+y)\cos q(x-y), p\sin p(x+y) \cos
q(x-y)\right),
$$
$$
\{a,b\}=\left (q\cos p(x+y) \sin q(x-y)-p\sin p(x+y) \cos q(x-y),
\qquad \qquad\qquad \right.
$$
$$
\left. \qquad \qquad \qquad -q\cos p(x+y) \sin q(x-y)- p\sin
p(x+y)\cos q(x-y)\right),
$$
$$
J \{a, b \}=\left (q\cos p(x+y)\sin q(x-y)+p\sin p(x+y) \cos
q(x-y), \qquad \qquad \qquad \right.
$$
$$
\left. \qquad \qquad \qquad q\cos p(x+y) \sin q(x-y) - p\sin
p(x+y) \cos q(x-y) \right),
$$
$$
\div J \{a,b \} = -2pq\sin p(x+y)\sin q(x-y) =
$$
$$
= pq\left (\cos ((p-q)x+(p+q)y)-\cos ((p + q) x + (p-q) y)
\right),
$$
$$
\mathrm {E}_g^{-1}(\div J \{a, b \}) =
$$
$$
= \frac {-2pq}{((p-q)^2 + (p + q)^2)^2} \left (\cos ((p-q)x +
(p+q)y)\cos ((p+q)x + (p-q)y) \right),
$$
$$
\int_{T^2} \left (\div J \{a, b \}\right) \mathrm {E}_g^{-1}
\left (\div J \{a, b \}\right)\ dx\ dy = \frac {2p^2 q^2}{(2p^2 +
2q^2)^2}(2\pi^2 + 2\pi^2),
$$
$$
\|a\|^2 = 4\pi^2 = \|b \|^2, \quad (a, b)_g = 0, \quad K(a,b)= 0.
$$
Therefore
$$
\overline {K}(\overline {a}, \overline {b}) = \frac {3}{32\pi^2}
\frac {p^2q^2}{(p^2 + q^2)^2}.
$$
If
$$
a = \left (\begin {array}{cc}
0 & -\cos p (x + y) \\
-\cos p (x + y) & 0
\end {array} \right), \quad
b = \left (\begin {array}{cc}
0 & -\sin q (x-y) \\
-\sin q (x-y) & 0
\end {array} \right), p^2 + q^2\neq 0.
$$
then similarly
$$
\overline {K}(\overline {a}, \overline {b}) = \frac {3}{32\pi^2}
\frac {p^2q^2}{(p^2 + q^2)^2}.
$$
Let now,
$$
a = \left (\begin {array}{cc}
0 & -\cos p (x + y) \\
-\cos p (x + y) & 0
\end {array} \right), \quad
b = \left (\begin {array}{cc}
0 & -\cos q (x + y) \\
-\cos q (x + y) & 0
\end {array} \right), p\neq q.
$$
In this case: $K(a,b)=0$, $(\framebox (8,8){}\ (b),a) = 0$, \
$(\framebox (8,8){}\ (a), b) = 0$, $\div J \{a, b \} = 0$ and
$\overline {K}(\overline {a}, \overline {b})$. If
$$
a = \left (\begin {array}{cc}
0 & -\cos p (x + y) \\
-\cos p (x + y) & 0
\end {array} \right), \quad
b = \left (\begin {array}{cc}
0 & -\sin q (x + y) \\
-\sin q (x + y) & 0
\end {array} \right), p\neq q.
$$
then also $\overline{K}(\overline{a},\overline{b})$.

{\bf Third part.} Let
$$
a=\left(\begin{array}{cc}
-\cos kx & 0 \\
0 & -\cos kx
\end{array}\right),\quad
b=\left(\begin{array}{cc}
0 & -\cos p(x+y)\\
-\cos p(x+y) & 0
\end{array}\right),\quad p^2+q^2\neq 0.
$$
In this case: $K(a,b) = -\frac {1}{8\pi^2}$,
$$
(\framebox (8,8){}\ (b), a) = \left (-p\cos kx\sin p (x + y),
p\cos kx\sin p(x+y) \right),
$$
$$
(\framebox (8,8){}\ (a), b) = \left (0, -k\sin kx\cos p(x+y)
\right),
$$
$$
A(\delta_g b) = \left (-p\cos kx\sin p(x+y), p\cos kx\sin p(x+y)
\right),
$$
$$
B(\delta_g a) = \left (0, -k\sin kx\cos p(x+y) \right),
$$
$$
\{a, b \} = 2\left (-p\cos kx\sin p(x+y), p\cos kx\sin p(x+y)
\right) + 2\left (0, k\sin kx\cos p(x+y) \right),
$$
$$
J \{a, b \} = -2p\left (\cos kx\sin p(x+y), \cos kx\sin p(x+y)
\right) - 2k\left (0, \sin kx\cos p(x+y) \right),
$$
$$
\div J \{a, b \} = 4kp\sin kx\sin p(x+y) - 2 (k^2 + 2p^2) \cos
kx\cos p(x+y) =
$$
$$
= (2kp-k^2-2p^2) \cos ((k-p)x-py)-(k^2 + 2kp + 2p^2) \cos
((k+p)x+py),
$$
$$
\mathrm {E}_g^{-1}(\div J \{a, b \}) =
2\frac{2kp-k^2-2p^2}{((k-p)^2 + p^2)^2}\cos ((k-p)x-py)- 2\frac
{k^2+2kp + 2p^2}{((k+p)^2 + p^2)^2}\cos ((k+p)x + py),
$$
$$
\int_{T^2} \left (\div J \{a, b \}\right)\ \mathrm {E}_g^{-1}%
\left (\div J \{a, b \}\right)\ dx\ dy =
$$
$$
= \left (2\frac {(2kp-k^2-2p^2)^2}{((k-p)^2 + p^2)^2} + 2\frac
{(k^2 + 2kp + 2p^2)^2}{((k + p)^2 + p^2)^2}\right) 2\pi^2 =
8\pi^2,
$$
$$
\|a\|^2 = 4\pi^2 = \|b\|^2, \quad (a, b)_g = 0,
$$
Therefore
$$
\overline {K}(\overline {a}, \overline {b}) = -\frac {1}{8\pi^2} +
\frac{3}{4} \frac {1}{2\pi^2} = \frac {1}{4\pi^2}.
$$

The same value of a sectional curvature turns out, when the
2-form $a$ is given by functions $\sin kx$, $\cos ly$, $\sin ly$
or when the 2-form $b$ is given by functions $\sin p(x+y)$, $\cos
q(x-y)$, $\sin q(x-y)$. If one of functions has a view $\pm 1$ or
$\pm i$, then $\{a, b\} = 0$. Therefore $\overline {K}(\overline
{a},\overline {b}) = K(a,b) = -\frac {1}{8\pi^2}$. The theorem is
proved.

\newpage

\vspace {10mm}

\centerline {\bf \S 9. Critical associated metrics} \centerline
{\bf on a symplectic manifold.}

\vspace {7mm}

Consider a functional on the space of Riemannian metrics
$$
R:\ \mathcal {M} \longrightarrow {\bf R}, \quad R(g)=\int_M r(g)\
d\mu (g),
$$
where $r(g)$ is scalar curvature of the metric $g\in \mathcal{M}$.
Functional $R$ is invariant with respect to an action of group of
diffeomorphisms $\mathcal {D}$ on $\mathcal {M}$. The gradient of
a functional $R$ is simple found \cite{Bes2}:
$$
dR(g;h) = \int_M \left(dr(g;h)+r(g)\frac 12 \tr_g h \right)\
d\mu(g) =
$$
$$
=\int_M \left(\triangle(\tr_g h)+\delta_g \delta_g h g(h,Ric(g))+
r(g)\frac 12\tr_g h\right)\ d\mu (g) =
$$
\centerline {(under the Stokes formula, taking into account
$\partial M =\emptyset$)}
$$
=\int_M g\left (\frac 12r(g)g-Ric(g),h\right)\ d\mu (g) =
\left(\frac 12r(g)g - Ric(g), h\right)_g,
$$
where $Ric(g)$ is Ricci tensor, $g(h,Ric(g))$ is pointwise inner
product of tensor fields on $M$, $g(h,Ric(g)) = \tr\left (g^{-1}
hg^{-1}Ric(g)\right)$, $\tr_g h = \tr\ g^{-1} h$. Thus,
$$
\grad\ R = \frac 12r(g)g - Ric (g).
$$
Recall, that a metric $g$ on $M$ is called Einsteinian, if its
Ricci tensor is proportional to a metric tensor $g$,
$$
Ric(g) = \lambda g, \quad \lambda \in {\bf R}.
$$

Einsteinian metrics are critical for the functional $R$ on the
manifold $\mathcal {M}_1$ of metrics with the same volume (equal
to unit). Recall, that
$$
T_g\mathcal {M}_1 =\{h\in S_2;\ \int_M\tr_g h\ d\mu (g)=0 \}.
$$
The orthogonal complement to $T_g \mathcal {M}_1$ in the space
$T_g \mathcal{M}=S_2$ consists of 2-forms proportional to $g$:\
$h = cg$, $c\in \mathbf {R}$. The metric $g\in \mathcal {M}_1$ is
critical for the functional $R$ on $\mathcal {M}_1$ if and only
if $\grad\ R(g)$ is orthogonal to $T_g \mathcal {M}_1$, i.e. if
for some number $c\in \mathbf {R}$,
$$
\frac 12 r(g)g - Ric(g) = cg.
$$
is held. We just obtain from here, that $r(g)=const$ and
$Ric(g)=\lambda g$, $\lambda\in \mathbf {R}$.

Consider the functional $R$ on the manifold $\mathcal {AM}$ of
associated metrics on a symplectic manifold $M, \omega$. In this
case a set of critical metrics is much wider. D.Blair has shown
\cite{Bla2}, \cite{Bla-Ian}, that a metric $g\in \mathcal {AM}$ is
critical for the functional $R$ on $\mathcal {AM}$ if and only if
Ricci tensor $Ric(g)$ is Hermitian with respect to almost complex
structure $J$, corresponding to $g$.

Indeed, any element $h\in T_g \mathcal {AM}$ represents symmetric
anti-Hermitian form, therefore $\tr_g h = 0$ and then
$$
dR(g;h)=\int_M \left(-g(h,Ric(g))+r(g)\frac 12 \tr_g h\right)\
d\mu(g) =
$$
$$
=-\int_M g(h,Ric(g))\ d\mu (g)=0.
$$

From an arbitrary of the anti-Hermitian form $h$ we obtain
$g(h,Ric(g))=0$, and it follows from a pointwise orthogonality of
the Hermitian and anti-Hermitian forms on $M$, that the tensor
$Ric(g)$ is Hermitian.

\vspace{3mm}

{\bf Conclusion.} In case of associated metrics analog of
Einstein metrics is metrics with Hermitian Ricci tensor.

\vspace{3mm}

The important property of the space $\mathcal {EM}$ of
einsteinian metrics on $M$ is the finite dimensionality of a
moduli space $\mathcal {EM}/\mathcal{D}$. Such question is
natural for asking for the space of critical metrics of the
functional $R$ on the space of associated metrics.

In this paragraph we show a finite dimensionality of the space of
classes of equivalent critical associated metrics of a constant
scalar curvature. Consider the map
$$
ARic:\ \mathcal {AM} \longrightarrow S_{2A}, \quad g
\longrightarrow Ric(g)_A, \eqno {(9.1)}
$$
which takes each associated metric $g$ to an anti-Hermitian part
$Ric(g)_A$ of Ricci tensor. Then the set of all critical metrics
coincides with a set $ARic^{-1}(0)$.

By symbol $\mathcal {CAM}_c$ we shall designate a set of the
critical associated metrics of  constant scalar curvature, which
is equal to $c$. The space $\mathcal {CAM}_c$ can be considered
as a level set $(ARic\times r)^{-1}(0,c)$ of mapping
$$
ARic\times r:\mathcal{AM}\longrightarrow S_{2A}\times
C^{\infty}(M,{\bf R}), \quad g \longrightarrow (Ric(g)_A, r(g)).
\eqno {(9.2)}
$$
This mapping extends to a smooth mapping
$$
ARic\times r:\mathcal{AM}^s\longrightarrow S^{s-2}_{2A}\times
H^{s-2}(M,{\bf R})
$$
from Hilbert manifold $\mathcal {AM}^s$ of associated metrics of
the Sobolev class $H^s$, $s\geq 2n+3$, into the Hilbert space
$S^{s-2}_{2A} \times H^{s-2}(M,{\bf R})$ of forms and functions of
class $H^{s-2}$. It was proven in \S 3 that the manifold
$\mathcal {AM}^s$ is analytical. In further the analyticity of
the continued mapping $ARic \times r$ will be shown. Let
$\{\mathcal {G}; \mathcal {G}^s,\ s \geq 2n+5 \}$ is ILH-Lie group
of exact symplectic diffeomorphisms \cite{Omo5}, \cite{Rat-Shm}
and let $I_{\omega}(g)=I(g)\cap \mathcal {G}$ is isometry group
being also symplectic transformations.

The following theorem takes place, it just follows from the slice
theorem, \cite{Ebi1}, \cite{Fuj-Sch}, stated in general case for
action of groups of diffeomorphisms on the space of metrics.

\vspace{3mm}

{\bf Slice theorem.} {\it Let $g\in \mathcal {AM}$. If $s\geq
2n+5$, then there exists a submanifold $\mathcal {S}^s_g$ in
$\mathcal {AM}^s$ and a local section $\chi^{s+1}:\ I_{\omega}(g)
\setminus \mathcal {G}^{s+1} \longrightarrow \mathcal {G}^{s+1}$
defined on an open neighbourhood $U^{s+1}$ of a coset
$[I_{\omega}(g)]$ which posses the following properties:

1) If $\gamma\in I_{\omega}(g)$, then $\gamma^*(\mathcal
{S}^s_g)=\mathcal {S}^s_g$.

2) Let $\gamma\in \mathcal {G}^{s+1}$. If $\gamma^* (\mathcal
{S}^s_g)\cap \mathcal {S}^s_g\neq \emptyset$, then $\gamma\in
I_{\omega}(g)$.

3) The mapping $F:\ \mathcal {S}^s_g \times U^{s+1}
\longrightarrow \mathcal {AM}^{s}$,
$F^s(g_1,u)=\chi^{s+1}(u)^*g_1$, is homeomorphism on an open
neighbourhood $V^s$ of the element $g$ from $\mathcal {AM}^s$.}

\vspace{3mm}

The general schema of construction of a slice $\mathcal {S}^s_g$
is applied to our case too. Let$ S^{*s}_{2A}$ is the space of
anti-Hermitian symmetric 2-forms $h$ on $M$ of a class $H^s$,
satisfying to the condition,
$$
\div J \delta_g h = 0.
$$
The slice $\mathcal {S}^s_g$ is an image of a neighbourhood of
zero $W^s\subset S^{*s}_{2A}$ at exponential mapping \cite{Ebi1}.
In case of the space of associated metrics the exponential map is
set by usual exponential mapping (see \S 3):
$$
Exp_g:\ S^{*s}_{2A} \longrightarrow \mathcal {AM}^s, \quad
Exp_g(h)=g\ e^H, \eqno {(9.3)}
$$
where $H=g^{-1}h$ and $e^H=1+H+\frac 12 H^2+\frac
{1}{3!}H^3+\dots$. As the mapping $Exp_g$ is real-analytic, then
the slice $\mathcal {S}^s_g$ is a real-analytic submanifold in
$\mathcal {AM}^s$. Note, that mapping $F$ of the slice theorem
$$
F:\ \mathcal {S}_g \times U \longrightarrow \mathcal {AM}
$$
is ILH-smooth, as for any $s\geq 2n+5$:
$$
\mathcal {S}^s_g = \mathcal {S}^{2n + 5}_g \cap \mathcal {AM}^s,
\qquad U^{s+1} = U^{2n+6}\cap \left (I_{\omega}(g)\setminus
\mathcal{G}^{s+1}\right),
$$
$$
V^s=V^{2n+5} \cap \mathcal{AM}^s, \qquad
\chi^{s+1}=\chi^{2n+6}|U^{s+1}
$$
is held, and for any $k\geq 0$ mappings
$$
F^{s+1}:\ \mathcal {S}^{s+k}_g \times U^{s+k+1} \longrightarrow
V^{s},
$$
$$
p^{s+k} \times q^{s+k}:\ V^{s+k} \longrightarrow \mathcal
{S}^{s}_g \times U^{s+1}
$$
are $C^k$-differentiable.

The quotient space $\mathcal {S}^{s}_g/I_{\omega}(g)$ describes a
local structure of a quotient space $\mathcal {AM}^s/\mathcal
{G}^{s+1}$ in a neighbourhood of class $[g]$.

Let $g\in \mathcal {CAM}_c$ is critical metric of a constant
scalar curvature, which is equal to $c$.

\vspace{3mm}

{\bf Definition 9.1.} {\it The set of critical associated metrics
of a constant scalar curvature $c$, which are in slice $\mathcal
{S}_g\subset \mathcal {AM}$ at a point $g$, is called as a
premoduli space of critical associated metrics of a constant
scalar curvature in a neighbourhood $g\in \mathcal {AM}$.}

\vspace{3mm}

Premoduli space will be denoted by a symbol $\mathcal {PM}(g)$,
$$
\mathcal {PM}(g)=\mathcal {S}_g \cap \mathcal {CAM}_c =
\left((ARic\times r)|_{\mathcal {S}_g}\right)^{-1}(0,c).
$$
The local moduli space is the quotient space $\mathcal
{PM}(g)/I_{\omega}(g)$. It describes a local structure of the
space $\mathcal {CAM}_c/\mathcal{G}$ in a neighbourhood of class
$[g] = g\mathcal {G}$.

Let $\mathcal {PM}^s(g)=\mathcal {S}^s_g\cap \mathcal {CAM}^s_c$
is premoduli space of critical metrics of a constant curvature of
Sobolev class $H^s$, $s\geq 2n+5$.

\vspace{3mm}

{\bf Theorem 9.1.} {\it Let $g\in \mathcal {CAM}_c$, then for any
$s\geq 2n+5$ there is a neighbourhood $W^s$ of the element $g$ in
slice $\mathcal {S}^s_g$ such, that the space $\mathcal
{PM}^s(g)\cap W^s$ is an analytical set of finite-dimensional
real-analytic submanifold $Z\subset W^s$, whose tangent space
$T_g Z$ has dimension independent on $s$ and consists of the
anti-Hermitian 2-forms of $C^{\infty}$-class.}

\vspace{3mm}

For a proof we will need some facts about analytic mappings of
Hilbert spaces. Recall, that the mapping $f$, defined on an open
set $U$ of Hilbert space $E_1$ in a Hilbert space $E_2$ is called
real-analytic, if in a neighbourhood of each point it is
represented by a convergent power series.

Let $E_1$, $E_2$ are complex Hilbert spaces and $U\subset E_1$ is
open set. The mapping $f:\ U \longrightarrow E_2$ is called
holomorphic, if $f$ is of $C^1$-class and in each point $x \in U$
the differential $df(x)$ commutes with complex structures $J_1$
and $J_2$ on $E_1$ and $E_2$.

We cite several known statements \cite{Fuc-Ne-So-So}.

\vspace{3mm}

{\bf Proposition 1.} {\it Let $E_1$ and $E_2$ are complex Hilbert
spaces and $U \subset E_1$ is open set. The holomorphic mapping
$f:\ U \longrightarrow E_2$ is real-analytic.}

\vspace{3mm}

{\bf Proposition 2.} {\it Let $E_1$ and $E_2$ are real Hilbert
spaces and $E^C_1, E^C_2$ are their \\ complexifications. Let $U
\subset E_1$ is open set and $f:\ U \longrightarrow E_2$ is
real-analytic mapping. Then there exists an open set $U^C \subset
E^C_1$ including $U$ and such that $f$ extends to a holomorphic
mapping $f^C:\ U^C \longrightarrow E^C_2$.}

\vspace{3mm}

{\bf Proposition 3} \cite{Koi}. {\it Let $E_1$ and $E_2$ are real
Hilbert spaces and $f$ is real-analytic mapping from $E_1$ in
$E_2$, defined on an open neighbourhood of zero $0\in E_2$.
Suppose, that $f(0)=0$ and image of a differential $df(0)$ is
closed in $E_2$. Then there is an open neighbourhood $U$ of zero
in $E_1$ such, that the set $f^{-1}(0)\cap U$ is a real-analytic
set in a real-analytic submanifold $Z$ from $U$, the tangent
space $T_0 Z$ of which coincides with $\Ker df(0)$.}

\vspace{3mm}

{\bf Lemma 9.1.} {\it Let $E$, $F$ are vector bundles over a
manifold $M$ and $E^C$, $F^C$ are their complexifications and
$\Gamma^ s (E)$, $\Gamma^s (F)$ are spaces of sections of the
Sobolev class $H^s$. Let $U^s \subset \Gamma^s(E)$ is open set and
$\Psi:\ U^s \longrightarrow \Gamma^{s-k}(F)$, $s\geq \frac n2
+k+1$ is $C^{\infty}$-smooth mapping possessing a property: for
any point $x\in M$ there is a neighbourhood $V_x \subset M$ of
this point, such, that $\Psi$ defines mapping $\Psi|V_x: U^s \cap
\Gamma^s (E|V_x) \longrightarrow \Gamma^{s-k}(F|V_x)$ of a class
$C^{\infty}$, which extends to holomorphic mapping
$$
\Psi^C|V_x :\ (U^s)^C \cap \Gamma^s(E^C|V_x) \longrightarrow
\Gamma^{s-k}(F^C|V_x),
$$
where $(U^s)^C$ is open set in $\Gamma^s (E^C)$, containing $U^s$.
Then mapping $\Psi:\ U^s \longrightarrow \Gamma^{s-k} (F)$ is
real-analytic.}

\vspace{3mm}

{\bf Proof.} Let $J_E$ and $J_F$ are operators of complex
structures in complexified fiber bundles $E^C$ and $F^C$. They
define complex structures in Hilbert spaces of sections $\Gamma^s
(E^C)$ and $\Gamma^{s-k} (F^C)$, which we shall designate by the
same symbols $J_E$ and $J_F$. Show, that the mapping $\Psi:\
(U^s)^C \longrightarrow \Gamma^{s-k}(F^C)$ is holomorphic, i.e.
for any section $u\in (U^s)^C$ the equality $d\Psi
(J_E(u))=J_F(d\Psi (u))$ is held. The operators $J_E$ and $J_F$
on the spaces $\Gamma^s (E^C)$ and $\Gamma^s (F^C)$ acts
pointwise:
$$
J_E(u)(x)=J_E(u(x)), \quad J_F(v)(x)=J_F (v(x)),
$$
therefore equality $d\Psi (J_E (u)) = J_F (d\Psi (u))$ needs to be
only checked up at each point $x\in M$:
$$
d\Psi (J_E (u))(x)=J_F (d\Psi (u)(x)) \quad x\in M.
$$
But the last is held on a condition in a neighbourhood $V_x$ of
each point $x\in M$.

\vspace{3mm}

{\bf Corollary 1.} {\it If the mapping $\Psi:\ U^s\longrightarrow
\Gamma^s(F^C)$, $s\geq \frac n2 +1$, has a view
$$
\Psi (u) = \psi\circ u,
$$
where $\psi:\ E \longrightarrow F$ is  saving fibers map of class
$C^{\infty}$, defined on an open set in $E$ and extendible to a
mapping $\psi^C:\ E^C \longrightarrow F^C$, such that its
restriction $\psi^C$ on each fiber $E^C_x$ is holomorphic, then
$\Psi$ is real-analytic mapping.}

\vspace{3mm}

It follows from here, that the tensor operations (convolution,
raising of an index etc.) determine analytical mappings of spaces
of sections.

The partial derivative of a tensor field $u$ with respect to
coordinate $x^i$ on base of $M$ represents a linear operation in a
neighbourhood of a point $x\in M$. Therefore it extends to
holomorphic mapping of complex fields $u$ in a neighbourhood of a
point $x$.

\vspace{3mm}

{\bf Corollary 2.} {\it If the mapping $\Psi:\ U^s
\longrightarrow \Gamma^{s-k}(F)$ in local coordinates in a
neighbourhood of each point $x \in M$ is analytically expressed
through tensor operations and partial derivatives of sections $u
\in U^s$ to the order $k$, then $\Psi$ is real-analytic mapping.}

\vspace{3mm}

In particular, mapping
$$
ARic\times r :\ \mathcal {M}^s \longrightarrow S^{s-2}_{2A}\times
H^{s-2}(M,{\bf R})
$$
is real-analytic mapping at $s\geq \frac n2 +3$, since it is
locally analytically expressed through the second partial
derivatives of a metric tensor $g$ and operations of convolution,
raising of an index, taking of an anti-Hermitian part.

The slice $\mathcal {S}^s_g$ is a real-analytic submanifold in
$\mathcal {AM}^s$, as it is an image of a neighbourhood of zero of
analytical mapping
$$
Exp_g: \ S^{*s}_{2A} \longrightarrow \mathcal {AM}^s, \quad Exp_g
(h) = g\ e^H,
$$
where $H=g^{-1}h$.

\vspace{3mm}

{\bf Proof of the theorem 9.1.} Consider analytical mapping
$ARic\times r:\mathcal{AM}^s\longrightarrow S^{s-2}_{2A} \times
H^{s-2} (M,{\bf R})$ and take its restriction on an analytical
submanifold $\mathcal {S}^s_g$, $g\in \mathcal {AM}$,
$$
\left. ARic\times r \right|_{\mathcal {S}^s_g}:\ \mathcal {S}^s_g
\longrightarrow S^{s-2}_{2A} \times H^{s-2}(M,{\bf R}).
$$
The premoduli space $\mathcal {PM}^s(g)$ is a level set,
$$
\mathcal {PM}^s (g) = \left (\left. ARic\times r \right
|_{\mathcal {S}^s_g} \right)^{-1} (0, c).
$$
Therefore we can apply a proposition 3. It is necessary to show,
that an image of a differential $d\left (\left. ARic\times r
\right |_{\mathcal {S}^s_g} \right) (g)$ is closed. At first we
find differential $d(ARic)(g)$ of mapping
$$
ARic:\ \mathcal {AM}^s \longrightarrow S of^{s-2}_{2A}.
$$
Let $h\in T_g \mathcal {AM}^s$ and let $g_t$ is curve which is
going out from $g$ in direction $h$ and $J_t$ is its corresponding
set of associated almost complex structures. The tangent vector
$I=\left.\frac {d} {dt} \right|_{t=0} J_t$ is expressed through
$h$:\ $I=J\circ H$, where $H=g^{-1}h$. Let $X,Y\in \Gamma (TM)$.
Differentiating of equality with respect to $t$
$$
ARic(g_t)(X, Y) = \frac 12\left (Ric(g_t)(X,Y)-
Ric(g_t)(J_tX,J_tY)\right),
$$
we obtain,
$$
\left (d(ARic)(g;h)\right)(X,Y)=A\left(dRic(g;h)(X,Y)\right)-
$$
$$
-\frac 12\left (Ric(g)(IX,JY)+Ric(g)(JX,IY) \right) =
$$
$$
= \frac 12 A\left(\Delta_L h - 2\delta^*_g \delta_g h\right)(X,Y)
- \frac 12\left (Ric(g)(IX,JY) + Ric(g)(JX,IY) \right),
$$
where $A\left (\Delta_L h - 2\delta^*_g \delta_g h\right)$ is
anti-Hermitian part, $\delta_g$ is covariant divergence,
$\delta^*_g$ is operator, which is adjoint to $\delta_g$,
$\Delta_L h=-\nabla_k \nabla^k h_{ij}+R_{ik} h^k_j + R_{jk} h^k_i
-2R_{i k j l}\ h^{kl}$ is Lichnerowicz Laplacian, we take into
account, that $h\in T_g\mathcal {AM}^s$ and, therefore, $\tr_g
h=0$, the expression for $dRic(g)$ is obtained in work
\cite{Ber-Eb}, see also \cite{Bes2}.

Consider, that the point $g\in \mathcal {AM}$, in which the
differential is calculated, is the critical metric. Then the Ricci
tensor $Ric(g)$ is Hermitian, therefore:
$$
d(ARic)(g; h)(X,Y) =
$$
$$
= \frac 12 A\left (\Delta_L h - 2\delta^*_g \delta_g h\right)(X,Y)
- \frac 12\left (Ric(g)(HX,Y) + Ric(g)(X,HY) \right).
$$
or, omitting arguments $X,Y$,
$$
d(ARic)(g;h) = \frac 12 A\left (\Delta_L h - 2\delta^*_g \delta_g
h\right) - \frac 12 \left (H^T\circ Ric(g) + Ric(g) \circ
H\right). \eqno {(9.4)}
$$
We have obtained, that $dARic(g)$ is a differential operator of
the second order.

The differential of the mapping $r:\ \mathcal {AM}^s
\longrightarrow H^{s-2}(M,{\bf R})$ is known \cite{Ber-Eb} to be
$$
dr(g,h) = \Delta\ \tr_g h+\delta_g \delta_g h -g(h,Ric (g)).
$$
Since $h$ is anti-Hermitian and  $Ric (g)$ is Hermitian, it
follows  that $\tr_g h = 0$ and $g (h, Ric(g)) = 0$. Therefore,
$$
dr(g;h)=\delta_g \delta_g h. \eqno {(9.5)}
$$
We have obtained, that the differential $d(ARic\times r)$ of
mapping $ARic\times r:\ \mathcal {AM}^s \longrightarrow
S^{s-2}_{2A}\times H^{s-2}(M,{\bf R})$ at a critical point $g\in
\mathcal {AM}$ is a differential operator of the second order
$$
d(ARic\times r)(g):\ S^{s}_{2A} \longrightarrow S^{s-2}_{2A}
\times H^{s-2} (M,{\bf R}),
$$
$$
h \longrightarrow \left (\frac 12 A\left (\Delta_L h -
2\delta^*_g \delta_g h\right) - \frac 12 \left (H^T\circ Ric(g) +
Ric(g) \circ H\right), \delta_g \delta_g h \right).
$$
Since the mapping $ARic$ is restricted on the slice $\mathcal
{S}^s_g$, we should to assume
$$
h\in T_g \mathcal {S}^s_g = \{h\in S^{s}_{2A}; \ \div\ J\
\delta_g h=0\}=S^{*s}_{2A},
$$
i.e. we should to impose the additional condition: $\delta_g J
\delta_g h = 0$.

The following lemma is a basis for the proof of the theorem 9.1.

\vspace{3mm}

{\bf Lemma 9.2.} {\it For any associated metric $g\in \mathcal
{AM}$ the differential operator
$$
\mathrm {F}_g:\ S_{2A} \longrightarrow S_{2A}\times C^{\infty}
(M,{\bf R}) \times C^{\infty}(M,{\bf R}),
$$
$$
h \longrightarrow \left (\frac 12 A\left (\Delta_L h -
2\delta^*_g \delta_g h\right) - \frac 12 \left (H^T\circ Ric(g) +
Ric (g) \circ H\right),\delta_g \delta_g h, \delta_g J \delta_g h
\right),
$$
has an injective symbol.}

\vspace{3mm}

{\bf Proof.} Recall the definition of a symbol of a differential
operator $D: \Gamma (E) \longrightarrow \Gamma (E)(F)$ of the
order $k$, where $E$, $F$ are vector bundles above $M$, and
$\Gamma(E)$, $\Gamma (F)$ are spaces of their sections \cite{Pal}.
Let $x\in M$ is an arbitrary point. For any covector $\xi \in
T^*_x M$ there is a function $f$ on $M$ such, that $f(x)=0$ and
$df(x)=\xi$. Let $h \in \Gamma (E)$ is section, then expression
$$
\sigma_{\xi}(D) h_x=\frac {1}{k!} D\left (f^k h\right)(x)
$$
depends only on section $h(x)$ of $h$ at a point $x$ and, thus
defines linear mapping $\sigma_{\xi}(D):\ E_x \longrightarrow F_x$
of fibers above a point $x$ of bundles $E$ and $F$, which is
named as a symbol of an operator $D$.

The symbol $\sigma (D)$ is called injective, if $\sigma_{\xi}
(D):E_x \longrightarrow F_x$ is injective for everyone $x\in M$
and everyone nonzero $\xi\in T_x^* M$. Let $x \in M$ is any fixed
point. Show an injectivity of a symbol $\sigma_{\xi} (\mathrm
{F}_g)$ of an operator $\mathrm {F}_g$.

To the associated metric $g$ on $M$ there corresponds an
associated almost complex structure $J$, and it defines
decomposition $TM^C = T^{10} \oplus T^{01}$ of complexification
$TM^C$ of a tangent bundle $TM$.

Let $\partial_1, \ \dots, \ \partial_n$ is basis of sections of a
bundle $T^{10}$ in a neighbourhood of a point $x$ and
$\overline{\partial}_1, \ \dots, \ \overline {\partial}_n$ is its
corresponding basis of sections of $T^{01}$. Choose dual basis
$dz^1, \ \dots, \ dz^n$, $d\overline {z}^1,\ \dots,\ d\overline
{z}^n$ of bundle $T^*M^{C} = T^{*10} \oplus T^{*01}$ (we pay
attention, that $dz^k$ is simple notation, complex coordinates
$z^1,\ \dots,\ z^n$ on $M$ in a neighbourhood of a point $x\in M$
can be not defined).

If $g_{\alpha\overline {\beta}} = g (\partial_{\alpha},
\overline{\partial}_{\beta})$, for our Hermitian form $g$ on $M$,
continued on a complexification $TM^C$ we have, $g =
2g_{\alpha\overline {\beta}} dz^{\alpha} d\overline {z}^{\beta}$
and $g_{\alpha\overline {\beta}} =\overline {g}_{\beta\overline
{\alpha}}$. Note, that if $\widetilde {g} = g_{\alpha\overline
{\beta}}dz^{\alpha} \otimes d\overline {z}^{\beta}$, then $g=
2{\rm Re}\ \widetilde {g}$. Take an arbitrary anti-Hermitian form
$h$. In a local coframe it has a view:
$$
h = h_{\alpha\beta} dz^{\alpha} dz^{\beta} + \overline {h} _
{\alpha\beta} d\overline {z}^{\alpha} d\overline {z}_{\beta},
$$
where $h_{\alpha\beta} = h (\partial_{\alpha}, \partial_{\beta})$,
$h_{\alpha\beta}=h_{\beta\alpha}$. Let $\widetilde {h} =
h_{\alpha\beta}dz^{\alpha} \otimes dz^{\beta}$, then $h = 2 {\rm
Re} \widetilde {h}$.

Note obvious equalities:
$$
2\widetilde {g} = g - i\omega, \quad 2\widetilde {h} = h + i\Im
(2\widetilde {h}), \quad \Im (2\widetilde {h}) (u, v) = - h (u,
Jv).
$$
Remark, that a symbol of an operator
$$
\mathrm {F}_g (h) = \left (\frac 12 A\left (\Delta_L h -
2\delta^*_g \delta_g h\right) - \frac 12 \left (H^T\circ Ric (g) +
Ric(g) \circ H\right),\delta_g \delta_g h,\delta_g J \delta_g h
\right)
$$
coincides with a symbol of an operator
$$
h \longrightarrow \left (\frac 12 A\left (\overline {\Delta} h -
2\delta^*_g \delta_g h\right),\delta_g \delta_g h,\delta_g
J\delta_g h\right), \eqno {(9.6)}
$$
where $\overline {\Delta}$ is Laplacian.

Let $\nabla_{\alpha} = \nabla_{\partial_{\alpha}}$ and
$\nabla_{\overline {\alpha}} = \nabla_{\overline
{\partial}_{\alpha}}$, then for the anti-Hermitian form $h =
h_{\alpha\beta} dz^{\alpha} dz^{\beta} + \overline
{h_{\alpha\beta}} d\overline {z}^{\alpha} d\overline {z}^{\beta}$
we have:
$$
(\overline {\Delta} h)_{\alpha\beta} = -\nabla_{\gamma}
\nabla^{\gamma}h_{\alpha\beta} + \nabla_{\overline {\gamma}}
\nabla^{\overline {\gamma}}h_{\alpha\beta},
$$
$$
\delta_g h = -\left (\nabla^{\gamma} h_{\alpha\gamma},
\nabla^{\overline{\gamma}}h_{\overline {\alpha} \overline
{\gamma}} \right),
$$
$$
J\delta_g h = -i\left (\nabla^{\gamma} h_{\alpha\gamma},
-\nabla^{\overline{\gamma}}h_{\overline{\alpha}\overline{\gamma}}\right),
$$
$$
\delta_g \delta_g = \nabla^{\alpha} \nabla^{\gamma}
h_{\alpha\gamma} + \nabla^{\overline{\alpha}}\nabla^{\overline
{\gamma}} h_{\overline {\alpha} \overline{\gamma}},
$$
$$
\delta_g J \delta_g = i\left (\nabla^{\alpha} \nabla^{\gamma}
h_{\alpha\gamma} - \nabla^{\overline {\alpha}} \nabla^{\overline
{\gamma}} h_{\overline {\alpha}\overline{\gamma}} \right),
$$
$$
A(\delta^*_g \delta_g h) = -\frac 12\left (\begin {array} {cc}
\nabla_{\alpha} \nabla^{\gamma} h_{\beta\gamma} +
\nabla_{\beta} \nabla^{\gamma} h_{\alpha\gamma} & 0 \\
0 & \nabla_{\overline {\alpha}}\nabla^{\overline
{\gamma}}h_{\overline {\beta} \overline {\gamma}} +
\nabla_{\overline {\beta}} \nabla^{\overline {\gamma}}
h_{\overline{\alpha} \overline {\gamma}} \end {array} \right),
$$
Therefore for $\xi = \xi_{\alpha} dz^{\alpha} + \xi_{\overline
{\alpha}} d\overline {z}^{\alpha}$ we obtain,
$$
\sigma_{\xi} (\overline {\Delta}) (h) = -\xi^{\gamma}
\xi_{\gamma} h_{\alpha\beta}- \xi^{\overline {\gamma}}
\xi_{\overline {\gamma}} h_{\alpha\beta}.
$$
$$
\sigma_{\xi}\left (A\delta^*_g \delta_g \right) (h) = -\frac 12
\left(\begin {array} {cc} \xi_{\alpha} \xi^{\gamma}
h_{\beta\gamma} +
\xi_{\beta} \xi^{\gamma} h_{\alpha\gamma} & 0 \\
0 & \xi_{\overline{\alpha}} \xi^{\overline {\gamma}} h_{\overline
{\beta} \overline {\gamma}} + \xi_{\overline {\beta}}
\xi^{\overline {\gamma}} h_{\overline{\alpha} \overline {\gamma}}
\end {array} \right),
$$
$$
\sigma_{\xi} \left (\delta_g \delta_g \right) (h) = \xi^{\alpha}
\xi^{\gamma} h_{\alpha\gamma} + \xi^{\overline {\alpha}}
\xi^{\overline {\gamma}} h_{\overline {\alpha} \overline
{\gamma}} = h (\xi,\xi),
$$
$$
\sigma_{\xi} \left (\delta_g J\ \delta_g \right) (h) = i\left
(\xi^{\alpha} \xi^{\gamma} h_{\alpha\gamma} - \xi^{\overline
{\alpha}} \xi^{\overline {\gamma}} h_{\overline {\alpha}
\overline {\gamma}} \right) = -2\Im \widetilde
{h}(\xi,\xi)=h(\xi,J\xi),
$$
We identify a covector $\xi\in T^*_x M$ to a vector from $T_x M$
with the help of metric tensor $g$.

Our task is to show, that if $\sigma_{\xi} (\mathrm {F}_g)(h)=0$
for nonzero $\xi$, then $h=0$. Assume, that $h$ satisfies to a
condition $\sigma_{\xi} (\mathrm {F}_g)(h)=0$ for nonzero $\xi\in
T^*_x M$. In this case,
$$
h(\xi,\xi) = 0, \quad h (\xi, J\xi) = 0
$$
and for any vector $u\in T_x M$ the equality is held
$$
\sigma_{\xi} A\left (\overline {\Delta}h-2\delta^*_g \delta_g
h\right)(u,u)=0.
$$
Calculate the left part.
$$
\sigma_{\xi} A\left(\overline {\Delta}h - 2\delta^*_g
\delta_gh\right)(u,u)=
$$
$$
=-g(\xi,\xi)h(u,u)+ \xi_{\alpha} u^{\alpha}
h_{\beta\gamma}u^{\beta} \xi of^{\gamma} + \xi_{\beta} u^{\beta}
h_{\alpha\gamma}u^{\alpha} \xi of^{\gamma} + \xi_{\overline
{\alpha}}u^{\overline {\alpha}}h_{\overline {\beta} \overline
{\gamma}} u^{\overline {\beta}} \xi^{\overline {\gamma}} +
\xi_{\overline {\beta}} u^{\overline {\beta}} h_{\overline
{\alpha} \overline {\gamma}} u^{\overline {\alpha}}
\xi^{\overline {\gamma}} =
$$
$$
= -g(\xi, \xi) h(u,u)+2\widetilde {g}(\xi,u) \widetilde {h} (u,
\xi) + 2\overline {\widetilde {g}(\xi, u) \widetilde {h}(u,\xi)} =
$$
$$
= -g(\xi,\xi)h(u,u) + 2 {\rm Re} \left (2\widetilde{g}(\xi,u)
\widetilde{h} (u, \xi) \right) =
$$
$$
= -g(\xi,\xi)h(u,u) + g(\xi,u)h(u,\xi)+\omega (\xi,u)\Im
(\widetilde {h}(u,\xi)) =
$$
$$
= -g(\xi,\xi)h(u,u)+g(\xi, u)h(u,\xi) - \omega (\xi,u)h(u,J\xi).
$$
Therefore for any vector $u\in T_x M$, we obtain:
$$
-g(\xi,\xi)h(u,u)+g(\xi,u)h(u,\xi) - \omega (\xi,u)h(u,J\xi) = 0.
$$
Taking into account, that $- \omega (\xi, u)=-g(J\xi,u)$, we
rewrite the last equality as:
$$
-g(\xi,\xi)h(u,u) + g(\xi,u)h(u,\xi) - g(J\xi,u)h(u,J\xi) = 0.
\eqno {(9.7)}
$$
It is held for any $u\in T_x M$. If the vector $u$ is orthogonal
to vectors $\xi$ and $J\xi$, two last addends vanish and we
obtain,
$$
-g(\xi,\xi)h(u,u) = 0,
$$
so $h(u,u)=0$ for any vector $u$, orthogonal to two-dimensional
subspace $\mathbf{R}\{\xi,J\xi \}$ spanned by the vectors $\xi$
and $J\xi$. If $u=\xi$ or $u=J\xi$, the equality (9.7) is held
for any form $h$. However we have two equalities: $h(\xi,\xi)=0$,\
$h(\xi,J\xi)=0$. Taking into account, that
$h(J\xi,J\xi)=-h(\xi,\xi)$, we obtain, that the 2-form $h$
vanishes on a two-dimensional subspace $\mathbf{R}\{\xi,J\xi \}$.
Thus, the 2-form $h$ vanishes on all the space $T_x M$. Therefore
$h=0$. The lemma is proved.

\vspace{3mm}

{\bf Ending of proof of the theorem 9.1.} The differential
operator
$$
\mathrm {F}_g: \ S_{2A} \longrightarrow S \times C^{\infty}
(M,{\bf R}) \times C^{\infty} (M,{\bf R}),
$$
$$
h\longrightarrow \left(d(ARic)(g;h),dr(g;h), \delta_g J \delta_g
h \right)
$$
has an injective symbol. Therefore its kernel $\Ker\ \mathrm
{F}_g\subset S^s_{2A}$ is finite-dimensional and consists of
forms $h$ of class $C^{\infty}$. It follows from an ellipticity
of an operator $\mathrm {F}_g^*\circ \mathrm {F}_g: \ S^s_{2A}
\longrightarrow S^s_{2A}$. Besides the following Berger-Ebin
decomposition in a direct sum of the closed orthogonal subspaces
\cite{Ber-Eb}, takes place:
$$
S^{s-2}_{2A} \times H^{s-2} (M,{\bf R}) \times H^{s-2} (M,{\bf
R}) = \Im (\mathrm {F}_g) \oplus \Ker (\mathrm {F}^{*}_g). \eqno
{(9.8)}
$$

Consider mapping
$$
ARic\times r: \ \mathcal {S}^s_g \longrightarrow S^{s-2}_{2A}
\times H^{s-2} (M,{\bf R}).
$$
The tangent space $T_g \mathcal {S}^s_g$ consists of 2-forms
$h\in S^s_{2A}$, satisfying to a condition: $\div J \delta_g h=0$.
Since
$$
\mathrm {F}_g(h)=\left(d(ARic\times r)(g;h),\div J\delta_g
h\right),
$$
then the image of a differential $d (ARic\times r)_g (T_g
\mathcal {S}^s_g)$ coincides with an image of a subspace $\{h \in
S^s_{2A}; \ \div J \delta_g h=0 \}$ under the action of the
operator $\mathrm {F}_g$. Therefore
$$
d\left (ARic\times r\right)_g \left (T_g \mathcal {S}^s_g\right) =
\Im (\mathrm {F}_g) \cap \left (S^{s-2}_{2A}\times H^{s-2}
(M,{\bf R}) \times \{0\} \right).
$$
So the image is closed as intersection of two closed subspaces.
Designate $F^{s-2}=d(ARic\times r)_g (T_g\mathcal {S}^s_g)$. Let
$$
p:\ S^{s-2}_{2A} \times H^{s-2} (M,{\bf R})\times \{0\}
\longrightarrow F^{s-2}
$$
is orthogonal projection on the closed subspace (in
correspondence with decomposition (9.8)), assume $q = {\rm id}
-p$. The mapping $p\circ(ARic\times r):\ \mathcal {S}^s_g
\longrightarrow F^{s-2}$ is analytical and its differential at a
point $g$ maps tangent space $T_g \mathcal{S}^s_g$ onto the whole
space $F^{s-2}$. By the implicit function theorem for analytical
mapping of Hilbert manifolds there is a neighbourhood $W^s$ of
the element $g$ in the slice $\mathcal {S}^s_g$ such, that the
set $Z=\left (P\circ (ARic\times r) \right)^{-1}(0) \cap W^s$
represents a real-analytic submanifold in $W^s$. Moreover the
tangent space $T_g Z$ coincides with the kernel
$\Ker\left(d(ARic\times r)_g\right)=\Ker\ \mathrm {F}_g$ of an
operator $\mathrm{F}_g$, which is finite-dimensional and consists
of the forms $h$ of class $C^{\infty}$. Therefore $Z$ is
finite-dimensional analytical submanifold in $\mathcal {S}^s_g$.
Now
$$
\mathcal {PM}^s (g) \cap W of^s = \left (q\circ (ARic\times
r)|Z\right)^{-1}(0)
$$
is analytical set in $Z$ as the inverse image of zero under the
analytical mapping $q\circ(ARic\times r):\ Z \longrightarrow
(F)^{\bot}$. The theorem is proved.

\end{document}